\documentclass[]{article}

\usepackage[utf8]{inputenc}
\usepackage[T1]{fontenc}

\usepackage[hmarginratio=1:1,top=32mm,columnsep=20pt]{geometry} 
\usepackage[hang, small,labelfont=bf,up,textfont=it,up]{caption} 
\usepackage[table]{xcolor}

\usepackage{graphicx}%
\usepackage[font=scriptsize,labelfont=bf]{caption}
\usepackage{amsmath,amssymb,amsfonts}%
\usepackage{amsthm}%
\usepackage{mathrsfs}%
\usepackage{textcomp}%
\usepackage{manyfoot}%

\usepackage{algorithm}%
\usepackage{algorithmicx}%
\usepackage{algpseudocode}%
\usepackage{comment}
\usepackage{listings}%
\usepackage{bm}
\usepackage{hyperref}
\usepackage{siunitx}

\usepackage{subfig}
\usepackage{array, multirow}
\usepackage{supertabular}
\usepackage[title]{appendix}
\usepackage{booktabs}%

\usepackage{amsmath}

\DeclareMathOperator*{\argmin}{arg\,min}

\newcommand{\ndof}{\ensuremath{\text{n}_{\text{dof}}}}
\newcommand{\nsnap}{\ensuremath{\text{N}_{\text{train}}}}
\newcommand{\ntot}{\ensuremath{\text{N}_{\text{tot}}}}
\newcommand{\neval}{\ensuremath{\text{N}_{\text{eval}}}}
\newcommand{\ntest}{\ensuremath{\text{N}_{\text{test}}}}
\newcommand{\mueval}{\ensuremath{\bm{\mu}_{\text{eval}}}}
\newcommand{\etaeval}{\ensuremath{\bm{\eta}_{\text{eval}}}}
\newcommand{\snapeval}{\ensuremath{\bm{s}_{\text{eval}}}}

\newcommand{\etatest}{\ensuremath{\bm{\eta}_{\text{test}}}}

\newcommand{\snappredeval}{\ensuremath{\tilde{\bm{s}}_{\text{eval}}}}
\newcommand{\reddim}{\ensuremath{\text{r}}}

\newcommand{\mediumdim}{\ensuremath{\text{r}_{\text{med}}}}

\usepackage{amssymb, latexsym, amsmath}
\usepackage{xcolor}
\usepackage{colortbl}

\usepackage{titlesec}
\titleformat{\section}[block]{\large\scshape\centering}{\thesection.}{1em}{} 
\titleformat{\subsection}[block]{\large}{\thesubsection.}{1em}{} 

\usepackage{abstract} 

\title{\textbf{Enhancing non-intrusive Reduced Order Models
with space-dependent aggregation methods}}

\date{ }
\author{Anna Ivagnes  \\ \small SISSA, International School for Advanced Studies, \\ \small Mathematics Area, mathLab, Trieste, Italy. \\ \small  \href{mailto:anna.ivagnes@sissa.it}{anna.ivagnes@sissa.it} \normalsize \and Niccolò Tonicello \\ \small SISSA, International School for Advanced Studies, \\ \small Mathematics Area, mathLab, Trieste, Italy. \\ \small  \href{mailto:niccolo.tonicello@sissa.it}{niccolo.tonicello@sissa.it} \\ \normalsize
Paola Cinnella \\ \small Institut Jean Le Rond D'Alembert, Sorbonne Universit\'e, Boîte 162, Tour 55-65, Paris, France \\ \small Mathematics Area, mathLab, Trieste, Italy. \\ \small  \href{mailto:paola.cinnella@sorbonne-universite.fr}{paola.cinnella@sorbonne-universite.fr} \normalsize  \and Gianluigi Rozza \\ \small  SISSA, International School for Advanced Studies, \\ \small Mathematics Area, mathLab, Trieste, Italy. \\ \small \href{mailto:gianluigi.rozza@sissa.it}{gianluigi.rozza@sissa.it} }
\begin{document}
\maketitle
\begin{abstract}
\noindent In this manuscript, we combine non-intrusive reduced order models (ROMs) with space-dependent aggregation techniques to build a \emph{mixed-ROM}. The prediction of the \emph{mixed} formulation is given by a convex linear combination of the predictions of some previously-trained ROMs, where we assign to each model a space-dependent weight. The ROMs taken into account to build the \emph{mixed} model exploit different \emph{reduction} techniques, such as Proper Orthogonal Decomposition (POD) and AutoEncoders (AE), and/or different \emph{approximation} techniques, namely a Radial Basis Function Interpolation (RBF), a Gaussian Process Regression (GPR) or a feed-forward Artificial Neural Network (ANN). The contribution of each model is retained with higher weights in the regions where the model performs best, and, vice versa, with smaller weights where the model has a lower accuracy with respect to the other models. Finally, a regression technique, namely a Random Forest, is exploited to evaluate the weights for unseen conditions.
The performance of the aggregated model is evaluated on two different test cases: the 2D flow past a NACA 4412 airfoil, with an angle of attack of $5$ degrees, having as parameter the Reynolds number varying between $\num{1e5}$  and $\num{1e6}$ and a transonic flow over a NACA 0012 airfoil, considering as parameter the angle of attack.
In both cases, the \emph{mixed-ROM} has provided improved accuracy with respect to each individual ROM technique.

\end{abstract}

\maketitle

\section{Introduction}\label{sec:intro}
Thanks to the constantly increasing computational resources provided by modern hardware architectures, Computational Fluid Dynamics (CFD), has become a fundamental tool for the design processes of the aeronautical industry \cite{slotnick2014cfd}.

In the last few years great improvements have been made in using scale-resolving simulations such as Large-Eddy Simulations (LES) for real-world aeronautical applications of CFD \cite{goc2020wall,goc2021large}. However, for full-scale problems such as full wings or fuselage of operating airplanes, their computational cost is still prohibitive in terms of CPU time even with the intensive use of powerful supercomputers. Consequently, with the additional cost of a reduced accuracy in the smallest details of the flow field, the use of Reynolds-Averaged Navier-Stokes (RANS) equations still represents the most common technique to deal with turbulent flows in the aerospace industry. Even if RANS are much less expensive from a computational point of view with respect to LES, they may still require a high computational effort for design process, which require repeated queries of the CFD model for computing the performance of new geometries and configurations. This is even more true if such computations are performed within an optimisation algorithm to find optimal designs. In such case, the use of Reduced Order Modelling (ROM) can significantly improve the computational efficiency of the overall algorithm with a small loss in terms of accuracy.
In other words, LES, RANS e ROMs for RANS are in this way ordered for decreasing computational cost. Depending on the different tasks required by the industrial needs, one or the other can be used in order to simulate the aerodynamics of aeronautical vehicles. 

We focus here on the use of ROMs for RANS, which are suitable for applications where many queries in the parametric space are needed or for scenarios where real-time response is needed.
Whereas using ROMs on LES or DNS data is still quite limited due to the large amount of degrees of freedom involved and the strongly multi-scale character of the flow, predicting unseen configurations via ROMs for RANS is much more doable and common in the reduced-order modelling community \cite{hijazi2020data,hijazi2020effort,zancanaro2021hybrid,zancanaro2022segregated,zancanaro2022finite}. In this framework, in fact, the assumption of finding a low-dimensional manifold describing the parametric evolution of the system is much more reasonable due to the smoothness of averaged flow fields in turbulence.
Thus, we focus hereafter on the use of ROMs for RANS over two different airfoil configurations, namely, a NACA4412 in subsonic conditions and a transonic NACA0012 test case. Both cases represent well-documented benchmarks often used in the literature to validate RANS models \cite{nguyen2007rans,garcia2014steady,tober2018evaluation,iorio2014direct}. Specifically, we consider here a set of alternative non-intrusive reduced order models corresponding to different reduction methodologies and approximations in the latent space. It will be shown that such models succeed more or less well in approximating the reference RANS solution in different regions of the computational domain, and it is often difficult to select \emph{a priori} a single best-performing ROM for all configurations. Such uncertainty in the ROM choice may be critical in view of their use for further design processes, and should be quantified and, if possible, reduced.

For that purpose, we propose here a novel approach to further improve reduced modeling accuracy while quantifying uncertainty, by constructing a spatially-dependent convex linear combination of them depending on their agreement with the reference full-order model. In many problems of interest it might be difficult to know a-priori which combination of reduction and approximation technique might be the best in the whole spatial and parametric domain. What is usually consists in testing different models, different architectures, andtuning different hyperparameters until a sufficiently overall satisfyingly model is obtained. Such a task is however rather arbitrary and time consuming. The methodology proposed in the following is not only able to identify automatically a optimal model combination, but it can also further improve the overall accuracy by using the information obtained from all the alternative models. For example, non-linear reduction techniques might be a reasonable choice for advection-dominated problems. However, such techniques might be so finely trained in capturing sharp features such as shock waves and boundary layers that they end up lacking accuracy in smoother regions of the flow where, instead, linear techniques might be more suitable.
In practice, complex flows contain simultaneous various physical processes leading to competing requirements for the ROM to be used.

An efficient framework for combining a set of alternative models based on their merit is the Model Aggregation framework \cite{stoltz_agregation_2010,devaine_forecasting_2013,deswarte_sequential_2019},  which  aims at  combining multiple  predictions stemming  from
various models --also termed experts or forecasters--  to provide a global, enhanced, solution. A solution that provides a spatially-constant weighting of competing model predictions being not optimal, is not optimal itself because  the accuracy of the models varies according to  the local flow physics. Therefore, it becomes interesting to combine the Model Aggregation approach with so-called Mixture-of-local-Experts techniques \cite{jacobs_adaptive_1991},
also  referred-to  as Mixture-of-Experts  \cite{yuksel_twenty_2012}  or  Mixture Models. In such approaches,  the  input
feature space  (covariate space) is softly  split into partitions where  the locally best-performing
models are assigned  higher weights.  The soft partitioning is  accomplished through parametric gate
functions, or a  network of hierarchical gate functions \cite{JordanJacobs_1994expertregions} that
rank the  model outputs with probabilities. However, Mixture-of-Experts tend to promote a single best model in
every soft partition, thus  accounting for  the spatial variation  of the  best model  but 
neglecting the uncertainty in model choice.
A methodology combining the Model Aggregation and the Mixture-of-Experts ideas has been recently introduced by De Zordo Banliat et al. \cite{de2024space}, who applied it to optimally combining a set of competing RANS models depending on some local flow features. The approach, named space-dependent model aggregation (XMA) was specifically designed to improve the prediction over individual component models.
In this work, we build on the XMA framework, and we propose a space-dependent aggregation model for combining ROMs.


The new methodology is tested on the two above-mentioned airfoil configurations. A variety of different analyses such as varying the dimensions of the latent space, and the influence of the chosen reduction/approximation techniques are presented in order to assess the robustness and suitability of the proposed approach for the prediction of external turbulent aerodynamics problems.

\section{Numerical methods}\label{sec:methods}

The use of RANS is extremely widespread in the aeronautical industry. Here we simply considered well-established models and set-ups as full-order models. Consequently, we focus this section more on non-intrusive ROMs and the aggregation algorithm to combine them.
Non-intrusive ROMs, described in Section \ref{met-roms}, are employed to predict an approximated version of the full-order fields for unseen operating conditions with a considerably reduced computational cost.
The predictions of different ROMs are then combined together with a model mixture technique, deeply explained in Section \ref{subsec:met-aggregation}.
\subsection{Non-intrusive Reduced Order Models}
\label{met-roms}
This part of the manuscript is dedicated to the presentation of the theory behind the non-intrusive model order reduction technique. This approach is fully data-based, namely it only exploits the information provided by high-fidelity simulations without any a-priori knowledge on the equations or numerical scheme used for the high-fidelity simulations. Therefore, in this framework, the governing equations are only used at the full order level to perform the offline simulations.

The original full order solutions, named \emph{snapshots}, correspond to the fields of interest of our problem and each snapshot is related to a specific set of parameters, i.e. a specific simulation setup. For instance, in both test cases we will consider the pressure and wall shear stress fields acting on the airfoil, and the velocity and pressure fields evaluated in the 2D space around the airfoil.

Being $s$ a generic field, if we name $\bm{s}_i=\bm{s}(\bm{\mu_i})$ the $i$-th snapshot, having as $\bm{\mu}_i \in \mathbb{R}^p$ the corresponding set of parameters, we can assemble the following snapshots' matrix:
\[
\bm{S} = \begin{bmatrix}
    \vert & \vert & &\vert\\
    \bm{s}_1(\bm{x}) & \bm{s}_2(\bm{x}) & \dots & \bm{s}_\text{N}(\bm{x})\\
    \vert & \vert & &\vert
\end{bmatrix} \in \mathbb{R}^{\ndof \times \nsnap}
\]
where \ndof{}  is the number of degrees of freedom (equal for each snapshot), and \nsnap{}  is the number of snapshots.
In particular, \ndof{} is the number of cells on the airfoil 1D surface (when considering as snapshots the pressure and wall shear stress acting on it), and the number of cells in the 2D external mesh (when considering the pressure and velocity fields around the airfoil).
\medskip

The question which is addressed by ROMs is the following.
\textbf{Question 1}: 
``\emph{How can we find the approximation $\tilde{\bm{s}}(\bm{\mu^*}) \simeq \bm{s}(\bm{\mu^*})$ where $\bm{\mu}^*$ represents an unseen configuration?}"
\medskip

The procedure to address the above-mentioned task is composed of two steps: a \emph{reduction} step, performing a compression onto a space with reduced dimensionality, and an \emph{approximation} step, where interpolation or regression techniques are used to predict the reduced representation for unseen parameters.
Finally, the reduced representation is backmapped into the full-order space to find the approximated full-dimensional solution.

\subsubsection{Reduction techniques}
\label{subsec:reduction}
This first step consists in a linear or nonlinear mapping of the original matrix of snapshots onto a matrix with reduced dimensions.
In the following part we will refer to the reduction step with the mapping $\mathcal{R}$.
So, the problem here addressed is the following.
\medskip

\textbf{Question 2}:
``\emph{How can we find a reduced representation $\bm{a}(\bm{\mu_i})$ of each $\bm{s}(\bm{\mu_i})$?}"
\medskip

We consider here two different reduction techniques: a linear one, the Proper Orthogonal Decomposition (POD), and a non-linear one, the AutoEncoder (AE).

\paragraph{Reduction through Proper Orthogonal Decomposition}
The POD~\cite{chatterjee2000introduction, berkooz1993proper, kerschen2005method} consists in the projection of the snapshots' matrix onto a space spanned by a limited number of the so-called \emph{modes}, which can be computed either via the correlation matrix or through a Singular Value Decomposition (SVD) technique~\cite{wall2003singular, klema1980singular} in the offline stage. 
The main hypothesis of the POD is that each snapshot can be approximated as a linear combination of the modes:
\[\bm{s}_i \simeq \sum_{j=1}^{\reddim} a_j(\bm{\mu_i}) \bm{\phi_j},\]
where $\{\bm{\phi_j}\}_{j=1}^\reddim$ are the modes and $\reddim \ll \nsnap $. The parameter $\reddim$ is established \emph{a-priori} based on energy criteria and/or the singular values decay. The terms $\{a_j(\bm{\mu_i})\}_{j=1}^\reddim$ are the reduced coefficients associated with the modes. In the case of a SVD procedure, we will have $\bm{S} = \bm{U\Sigma V^T}$, where the columns of $\bm{U}$ are the POD modes and the operation $\bm{U^T s_i}=\bm{a}(\bm{\mu_i})$ provides the reduced representation.

\begin{figure}
    \centering
    \includegraphics[width=0.8\textwidth]{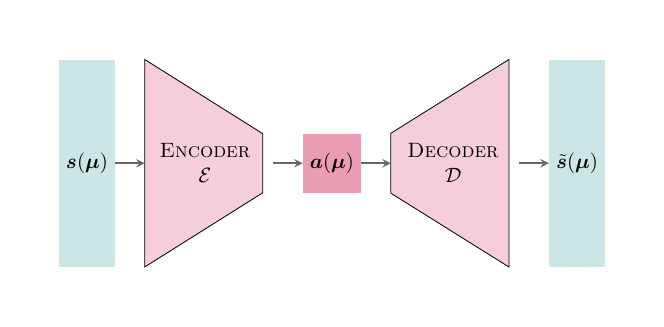}
    \caption{General structure of an autoencoder.}
    \label{fig:ae-scheme}
\end{figure}

\paragraph{Reduction through AutoEncoder}
Autoencoders are a family of neural networks that have become popular as dimensionality reduction technique~\cite{lee2020model, eivazi2020deep, ivagnes2023towards, halder2022non, kadeethum2022non, romor2022non} thanks to their peculiar architecture.
Indeed, AEs are composed of two components, as showed in Figure \ref{fig:ae-scheme}: an \emph{encoder} $\mathcal{E}:\mathbb{R}^{\ndof}\rightarrow \mathbb{R}^{\reddim}$, which maps the high-dimensional input into a latent space (the reduced space), and a \emph{decoder} $\mathcal{D}:\mathbb{R}^{\reddim}\rightarrow \mathbb{R}^{\ndof}$, which back-maps the latent representation onto the original full dimensionality.

Autoencoders can be composed by either convolutional or dense layers, but in this work we focus on dense feed-forward neural networks for both the encoder and the decoder. 

For instance, if we consider a single-layer structure for both the decoder and the encoder, we will have:
\[
\begin{split}
    &\mathcal{E}(\bm{s}(\bm{\mu})) = \sigma(\mathrm{W} \bm{s}(\bm{\mu}) + \mathrm{b}) = \bm{a}(\bm{\mu}),\\
    &\mathcal{D}(\bm{a}(\bm{\mu})) = \sigma^{\prime}(\mathrm{W}^{\prime} \bm{a}(\bm{\mu}) + \mathrm{b}^{\prime}) = \tilde{\bm{s}}(\bm{\mu}) \simeq {\bm{s}}(\bm{\mu}).
\end{split}
\]
The general structure with more hidden layers can be easily derived from the above formulation.

In the training of the AE, the loss function which is minimized is the following:
\begin{equation}
\frac{1}{\nsnap} \sum_{i=1}^{\nsnap}\left(\|\bm{s}(\bm{\mu_i})- \tilde{\bm{s}}(\bm{\mu_i})\|_{l^2}^2\right) = \frac{1}{\nsnap} \sum_{i=1}^{\nsnap}\left(\|\bm{s}(\bm{\mu_i})-  \sigma^{\prime}(\mathrm{W}^{\prime} (\sigma(\mathrm{W} \bm{s}(\bm{\mu}) + \mathrm{b})) + \mathrm{b}^{\prime})  \|_{l^2}^2\right),
    \label{eq:loss-ae}
\end{equation}
where $\tilde{\bm{s}}(\cdot)$ is the approximated high-fidelity field, and $\bm{s}(\cdot)$ is the ground truth.

\paragraph{Reduction through PODAE}
The mixed approach PODAE consists in a two-steps reduction, where the total mapping $\mathcal{R}:\mathbb{R}^{\ndof}\rightarrow \mathbb{R}^{\reddim}$ is expressed as $\mathcal{R}=\mathcal{R}_{\text{AE}} \circ \mathcal{R}_{\text{POD}}$.

The first reduction $\mathcal{R}_{\text{POD}}:\mathbb{R}^{\ndof}\rightarrow \mathbb{R}^{\mediumdim}$ is, in general, the most expensive from a computational point of view. In our case, it coincides to the POD retaining all the modes, namely $\mediumdim \equiv \nsnap$. On the other hand, the second reduction $\mathcal{R}_{\text{AE}}:\mathbb{R}^{\mediumdim}\rightarrow \mathbb{R}^{\reddim}$ is operated by an autoencoder.

The reason for employing such a reduction technique is the considerable gain in the computational time with respect to a pure nonlinear AE. The computational effort needed to train an autoencoder is, indeed, much higher than in the case of a mixed nonlinear approach.
However, since the PODAE includes a nonlinear mapping, it allows to capture complex behaviors that the POD is not able to capture.
In Section \ref{sec:results} we will employ this technique when $\ndof$ is of the order of magnitude of $10^6$.
\medskip

It is important to highlight that in all the numerical tests we considered the same latent dimension for both the POD, the AE and the PODAE techniques.

\subsubsection{Approximation techniques}
\label{subsec:approximation}
Once we have all the reduced representations of the snapshots, we can address the final task.
\medskip

\textbf{Question 3}:
``\emph{How can we find the reduced representation $\bm{a}(\bm{\mu^*})$, starting from the knowledge of $\{\bm{a}(\bm{\mu_i})\}_{i=1}^{\nsnap}$?}"
which can be also written as:
``\emph{Which is the mapping $\mathcal{M}:\mathbb{R}^p \rightarrow \mathbb{R}^r$ such that $\bm{a}(\bm{\mu_i}) = \mathcal{M}(\bm{\mu_i})$?}"
\medskip

We here consider different types of mapping $\mathcal{M}$: a Radial Basis Function (RBF) interpolation technique, a Gaussian Process Regression (GPR) technique, and a feed-forward Artificial Neural Network (ANN).
We provide in the following part a summary of the above-mentioned techniques, without going through all the details.

\paragraph{Radial Basis Function Interpolation}
Following the RBF~\cite{powell1987radial} technique, the unknown coefficients are computed as follows:
\begin{equation}
\mathcal{M}(\bm{\mu^*}) = \bm{a}(\bm{\mu^*}) = \sum_{i=1}^{\nsnap} \omega_i \varphi(\|\bm{\mu^*} - \bm{\mu_i}\|),
\label{eq:rbf_formula}
\end{equation}
where $\varphi(\|\bm{\mu^*} - \bm{\mu_i}\|)$ are the radial basis functions having center $\bm{\mu_i}$ and weight $\omega_i$. The weights are found considering known conditions, namely Equation \eqref{eq:rbf_formula} with $\bm{\mu^*}=\bm{\mu_j}$, $j=1, \dots, \nsnap$.
More details about the kernel are provided in the supplementary section.

\paragraph{Gaussian Process Regression}
The GPR~\cite{williams1995gaussian} is a supervised learning technique that exploits a stochastic model to build the regression function $\mathcal{M}$. In particular, the GPR model assumes that the regression function that relates the inputs $\bm{\mu}$ to the outputs $\bm{a}$ is drawn from a Gaussian process with mean function $\mathrm{m}$ and covariance function $\mathrm{K}$:
\begin{equation}
\mathcal{M}(\bm{\mu}) \sim \mathcal{N}(\mathrm{m}(\bm{\mu}), \mathrm{K}(\bm{\mu})),
    \label{eq:GP}
\end{equation}
where $\mathrm{K}_{ij}(\bm{\mu})=\mathcal{K}(\bm{\mu_i}, \bm{\mu_j})$. The shape of the covariance matrix, also called \emph{kernel}, is chosen depending on the particular problem of interest. The details of the kernel chosen for this problem are provided in the supplementary results' section.

\paragraph{Artificial Neural Network}
The last approximation technique considered in this work is a regression through a feed-forward neural network, which have been subject of a lot of research works in recent years in several fields, and in particular in reduced order modeling, such as ~\cite{hesthaven2018non, salvador2021non, bhattacharya2021model, san2019artificial, ivagnes2023towards}.

The architecture that we consider is a fully-connected neural network, defined as the concatenation of an input layer, multiple hidden layers, and a final output layer. If we isolate just an hidden layer, we can express the $i$-th output of the $h$-th layer as follows:

\begin{equation}
o_i^h = \sigma \left( \sum_{k=1}^{n_{h-1}} \mathrm{W}_{ik}^h o_k^{h-1} + b_i^h \right),
    \label{eq:output_1neuron}
\end{equation}
where $\sigma$ is the activation function representing the non-linearity of the model, $\mathrm{W}^h$ and $ b^h $ are the weights' matrix and the bias of the $h$-th layer.
All weights and bias are tuned during the training procedure in order to fit the input data, and our function can be here expressed as: $\mathcal{M}_{\theta}(\bm{\mu})$, where $\theta$ is the set of all the weights of the ANN.
The architecture of the ANN is reported in the supplementary results' part.

\subsection{Space-dependent aggregation methods and ROM mixture}
\label{subsec:met-aggregation}

The present section of the article is dedicated to the explanation of the logic of the Space-dependent Model Aggregation (XMA)~\cite{de2024space, cherroud2023space}. The main principle is to mix models into an ``hypermodel” (the model mixture) with improved predictive capabilities with respec to the mixture components themselves.

We call the set of all the models taken into account for the mixture $\mathcal{M} = \{M_1, M_2, \dots, M_{n_M}\}$, where $n_M$ is the number of models we are considering.

In our case, each $M_i$ correspond to a non-intrusive reduced order model, characterized by specific reduction technique $\mathcal{R}$ and approximation technique $\mathcal{A}$, i.e., 
$M_i = M(\mathcal{R}, \mathcal{A})$.

Moreover, we can write the prediction of each model $M_i$ as follows:

\[
\tilde{\bm{s}}^{(i)} = \tilde{\bm{s}}(\bm{\eta}, M_i), i=1, \dots, n_M,
\]
where $\bm{\eta}$ is the set of independent features on which the prediction depends, in our case the space coordinates and the ROM parameters: $\bm{\eta} = [\bm{x}, \bm{\mu}]$.

The model mixture is based on the assumption that the aggregated prediction can be written as a convex combination of the given models, namely:
\[
\tilde{\bm{s}}^{(\text{mix})} = \sum_{i=1}^{n_M} \omega_i(\bm{\eta}) \tilde{\bm{s}}^{(i)},
\]
where $\{\omega_i(\bm{\eta})\}_{i=1}^{n_M}$ are the weights' distribution associated with each model.

In particular, the weights have to satisfy the following conditions:
\[ \omega_i(\bm{\eta}) \in [0, 1], \quad \sum_{i=1}^{n_M}\omega_i(\bm{\eta}) = 1 \quad \forall \bm{\eta}\, . \]

At this point, it is important to specify that we split the original dataset of parameters and snapshots $\{(\mu_i, \bm{s})\}_{i=1}^{\ntot}$ into the following three parts:
\begin{itemize}
    \item a \emph{training} set, used to train all the reduced order models in $\mathcal{M}$, composed of $\nsnap$ components;
    \item an \emph{evaluation} set, used to compute the optimal weights $\{\omega_i\}_{i=1}^{n_M}$, composed of $\neval$ components;
    \item a \emph{test} set, composed of $\ntest$ components and used to compare the results among the single models and the mixture. In this set, the weights are computed with a \emph{regression} technique.
\end{itemize}

The key question that we have to address now is:

\textbf{Question 4}: 
``\emph{How can we compute the weights in the evaluation set in an optimal way?}"

We consider the following expression for the weights \cite{de2024space}:

\begin{equation}
\omega_i(\etaeval) = \dfrac{g_i(\etaeval)}{\sum_{j=1}^{n_M} g_j(\etaeval)},\, \; \text{where} \; \, g_j(\etaeval)=\exp{\left( -\dfrac{1}{2} \dfrac{(\snappredeval^{(j)}-\snapeval)^2}{\sigma^2} \right)}.
    \label{eq:weights}
\end{equation}
In the above formulation we consider $\etaeval=[\bm{x}, \mueval]$, where $(\mueval, \snapeval)$ is the evaluation set. Instead, $\snappredeval^{(j)}$ is the result predicted by the $j$-th model.
The parameter $\sigma$ is evaluated in each regression with an optimization algorithm, as in \cite{deswarte2019sequential} in the following way:
\begin{equation}
    \sigma_{\text{opt}} = \argmin_{\sigma \in \mathcal{I}} \sum_{k=1}^{\neval} \left( (\snapeval)_{k} - \sum_{l=1}^{n}\omega_l((\etaeval)_k;\sigma)
    \snappredeval(\etaeval)_k  \right)^2,
    \label{eq:opt-sigma}
\end{equation}
where $\mathcal{I}$ is the range where to search the optimal value of $\sigma$. As investigated in detail in \cite{deswarte2019sequential, de2024space}, the algorithm \eqref{eq:opt-sigma} can be used to establish the optimal \emph{order of magnitude} of $\sigma$, since it slightly influence the results if the order of magnitude is the optimal one.

Finally, one can train a regression technique for each model considering $(\etaeval, \omega_{j}(\etaeval))$ for each $j=1, \dots, n_M$ as data and then let the regression algorithm predict $\omega_{j}(\etatest)$ for each $j=1, \dots, n_M$. Being inspired from previous works such as \cite{cherroud2023space, de2024space}, a Random Forest Regression has here been used. However, any regression technique may be used to address this task, such that a GPR, or a neural network for example.

\paragraph{Random Forest regression}
The Random Forest (RF) regression~\cite{breiman2001random, liaw2002classification, rodriguez2015machine} is a machine learning algorithm used for supervised regression, which relies on different hyperparameters, namely the number of trees, the criterion for node splitting, and the minimum number of samples in a leaf. In our case the number of trees is set to 100, the criterion for splitting is the mean squared error, the minimum number of samples in a leaf is set to 2.

In general, at this point there are two operative choices:
\begin{itemize}
    \item train $n_M-1$ regressions, having as output $\{\omega_i(\bm{\eta}^{\star})\}_{i=1}^{n_M-1}$ and compute the remaining weights as $\omega_{n_M}(\bm{\eta}^{\star})=1-\sum_{i=1}^{n_M-1} \omega_i(\bm{\eta}^{\star}) $, for a fixed $\bm{\eta}^{\star}$;
    \item train separately $n_M$ regressions, having as outputs the Gaussians $\{g_j(\bm{\eta}^*)\}_{j=1}^{n_M}$ and then normalize the Gaussians as in \eqref{eq:weights}.
\end{itemize}

In this work, we adopt the second strategy.
It is also important to highlight that in some cases it may happen that the sum $\sum_{j=1}^{n_M} g_j (\etaeval)$ is considerably small, leading the division in Equation \eqref{eq:weights} to be ill-conditioned. To avoid this mathematical issue, we setup a minimum value for the range $\mathcal{I}$ in Equation \eqref{eq:opt-sigma} to search the optimal value of $\sigma$. In particular, we set $\mathcal{I}=[1e-3, 1]$.

\section{Numerical results}
\label{sec:results}
In the present Section, we present the numerical results obtained by aggregating different combinations of ROMs for two different test cases:
\begin{itemize}
    \item[\textbf{(i)}] steady and 2D flow past a NACA 4412 airfoil, having as parameter the inlet velocity, namely the Reynolds number. In particular, the Reynolds number is in the range $[\num{1e5}, \num{1e6}]$;
    \item[\textbf{(ii)}] transonic and 2D flow past the symmetric NACA 0012 airfoil, having as parameter the angle of attack, which varies within the range $[0^{\circ}, 10^{\circ}]$ at fixed Reynolds number $Re=\num{1e7}$.
\end{itemize}

The results' Section is divided into the following parts:
\begin{itemize}
    \item description of the RANS models here considered in both test cases, and validation of them with respect to the high-fidelity and/or experimental counterparts (section \ref{subsec:validation});
    \item analysis of the POD eigenvalues' decay for both cases, and, brief study of the ROMs' results (section \ref{subsec:rom-result});
    \item presentation of the results of the aggregated models for test case 1 (section \ref{subsubsec:aggr-result-1}) and 2 (section \ref{subsubsec:aggr-result-2}), and finally a brief discussion of the results, including a comparison between the performances in the two test cases.
\end{itemize}

\subsection{FOM description and validation}
\label{subsec:validation}

This Section is dedicated to the presentation of the full-order models considered to obtain the datasets for the ROMs in the two test cases. The software used to develop both the FOMs is OpenFOAM~\cite{ofsite}, a CFD open-source software, based on a finite-volume space discretization.
The following parts, namely \ref{subsubsec:fom-case1} and \ref{subsubsec:fom-case2} are dedicated to the description of the single FOMs for test case 1 and 2, respectively.

\subsubsection{Test case 1}
\label{subsubsec:fom-case1}
The first test case is the 2D steady flow past a NACA 4412 airfoil at fixed angle of attack $\alpha=5^{\circ}$. This case is simulated employing the Reynolds--Averaged Navier--Stokes (RANS) formulation, which consists of a time-averaged version of the Navier--Stokes equations.

The main hypothesis that characterizes the RANS approach is the \emph{Reynolds decomposition}~\cite{reynolds1895iv}.
This theory is based on the assumption that each flow field can be expressed as the sum of its mean and fluctuating parts. Such mean has different definitions depending on the case of application. Classical choices are, for example, time averaging, averaging along homogeneous directions or ensemble averaging.

We briefly recall here the standard RANS formulation for the incompressible Navier-Stokes equations:
\begin{equation}
\begin{cases}
\dfrac{\partial \overline{u}_{i}}{\partial x_i} = 0 ,\smallskip\\
    \overline{u}_{j} \dfrac{\partial \overline{u}_{i} }{\partial x_j}=-\dfrac{\partial \overline{p}}{\partial x_i}+ \dfrac{\partial( 2 \nu\overline{\mathbf{E}}_{ij}- \mathcal{R}_{ij})}{\partial x_j},
    \label{RANS}
\end{cases}
\end{equation}

where the Einstein notation has been adopted, $\mathcal{R}_{ij}=\overline{u'_i u'_j}$ is the Reynolds stress tensor, and $\overline{\mathbf{E}}_{ij}=\dfrac{1}{2}\left(\dfrac{\partial \overline{u}_{i}}{\partial x_j} + \dfrac{\partial \overline{u}_{j}}{\partial x_i}\right)$ is the averaged strain rate tensor.

The RANS formulation in \eqref{RANS} needs to be coupled with a turbulence model to close system \eqref{RANS}. In particular, we adopt the $\kappa-\omega$ \emph{Shear Stress Transport} (SST) model~\cite{menter1994two}.

This model belongs to the class of \emph{eddy viscosity models}, which are based on the Boussinesq hypothesis, i.e. the turbulent stresses are related to the mean velocity gradients as follows:
\[
 -\mathcal{R}_{ij}=2 \nu_t \overline{\mathbf{E}}_{ij} - \dfrac{2}{3} \kappa \delta_{ij},
    \label{bouss}
\]
where $\kappa=\frac{1}{2}\overline{u'_i u'_i}$ is the turbulent kinetic energy and $\nu_t$ is the eddy viscosity.
For the complete model we refer the reader to the original paper \cite{menter1994two}.

\begin{figure}[htpb!]
    \centering
    \subfloat[Domain and notation]{\includegraphics[width=0.5\textwidth, trim={5cm 0.5cm 3cm 0}, clip]{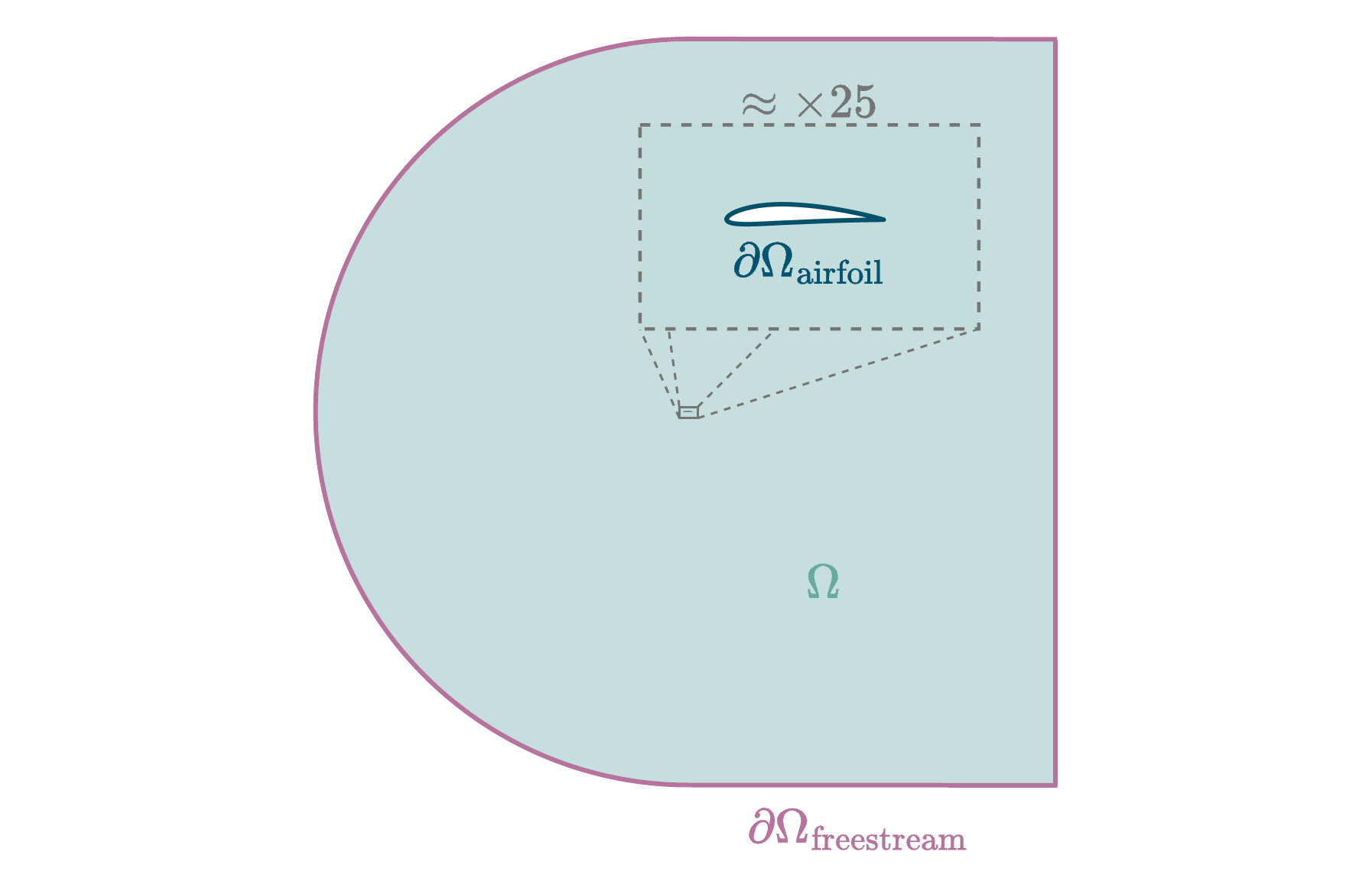}
    \label{case1domain}}
     \subfloat[Mesh representation]{\includegraphics[width=0.5\textwidth, trim={9cm 0cm 13cm 2cm}, clip]{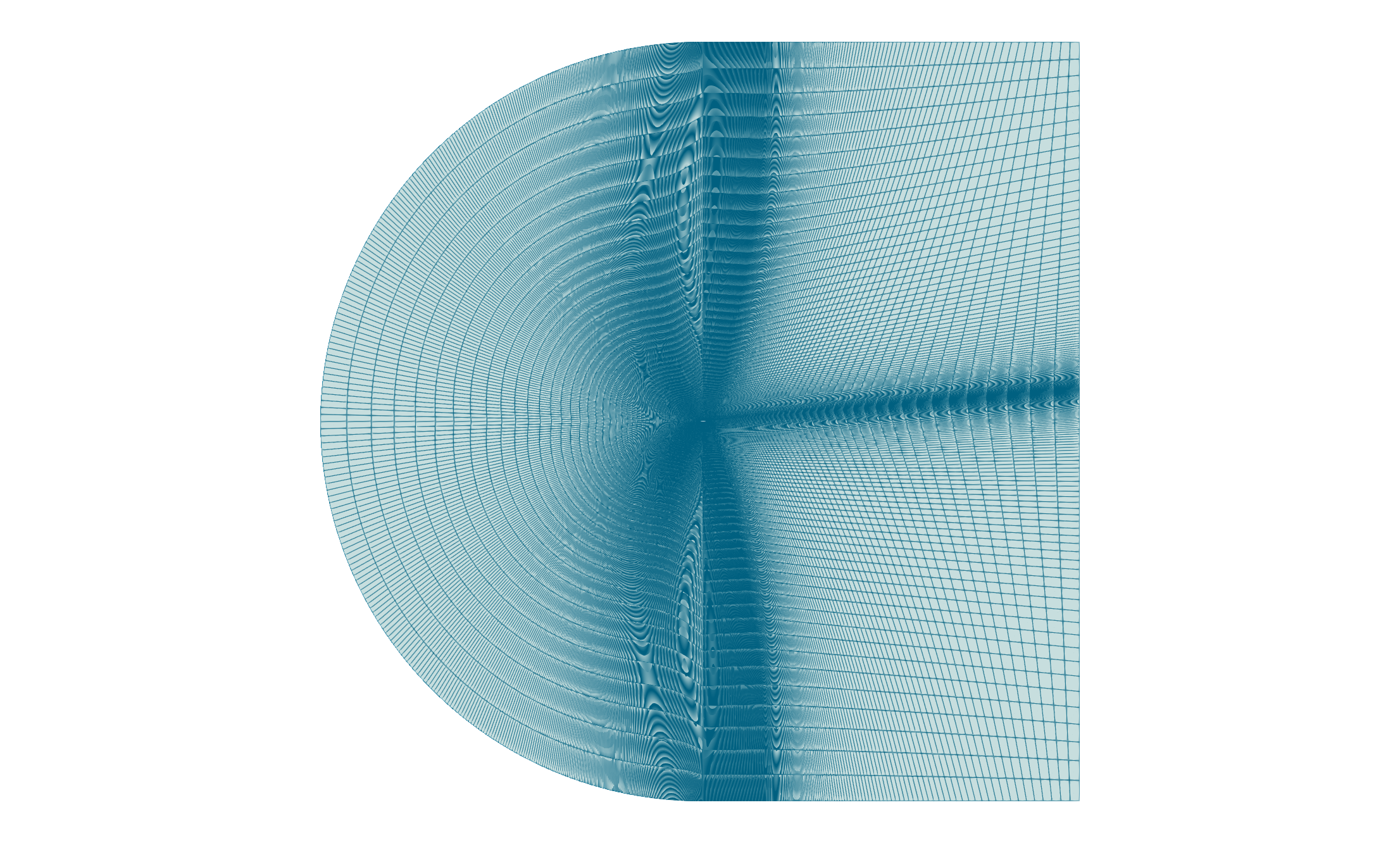}
     \label{case1mesh}}\\
     \subfloat[Zoomed mesh around the airfoil]{\includegraphics[width=0.35\textwidth]{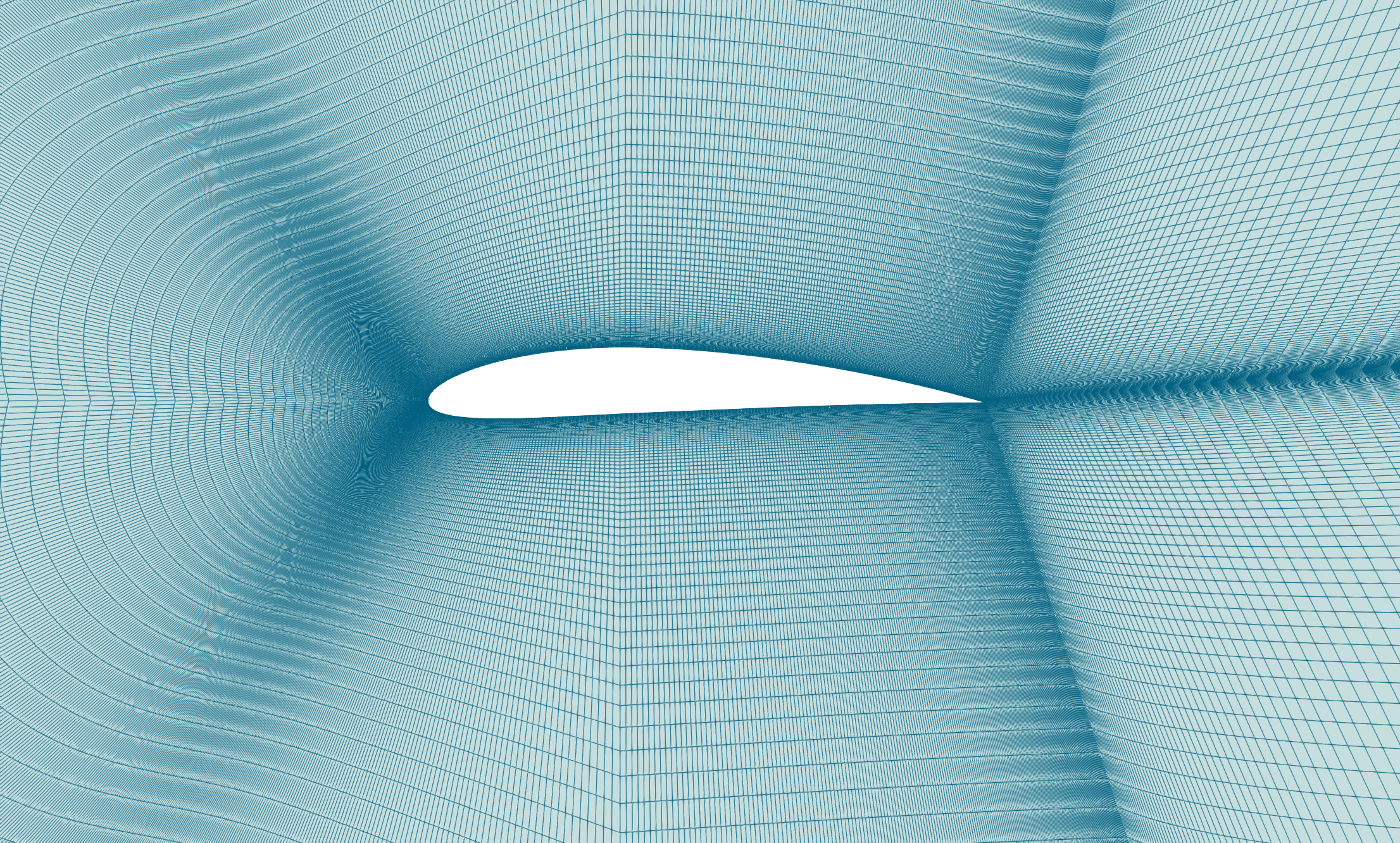}
     \label{case1mesh-zoom}}
    \caption{Domain (\protect \subref{case1domain}) and mesh (\protect \subref{case1mesh}) representation for test case 2, with a zoomed representation for the zone around the airfoil (\protect \subref{case1mesh-zoom}).}
    \label{fig:domain-case1}
\end{figure}

The domain and the mesh considered for this simulation are represented in Figure \ref{fig:domain-case1}. The domain boundaries $\partial \Omega_{\text{freestream}}$ are at a distance $100$ times the airfoil chord length, and so essentially in a freestream condition. Hence, the boundary conditions are the following:
\begin{equation}
\begin{cases}
    \bm{u} = \bm{u}_{\text{freestream}},\\
    p = p_{\text{freestream}},
\end{cases}
\text{on }\partial \Omega_{\text{freestream}}, \quad
\begin{cases}
    \bm{u} = \bm{0},\\
    p = 0,
\end{cases}
\text{on }\partial \Omega_{\text{airfoil}},
\end{equation}
where the freestream condition corresponds to a mixed condition, where the mode of operation switches between fixed value and Neumann, based on the sign of the flux. The fixed velocity value is evaluated starting from the Reynolds number, the dataset parameter, whereas the fixed pressure value is $p=0$.

The pressure-velocity coupling is treated considering the Semi-Implicit Method for Pressure-Linked Equations (SIMPLE) algorithm. For more details, we refer the reader to \cite{patankar1983calculation, patankar2018numerical}.

\begin{figure}[htpb!]
    \centering
      \subfloat[FOM pressure ($Re\approx 118000$)]{\includegraphics[width=0.25\textwidth, trim={2cm 0.85cm 1cm 1.1cm}, clip]{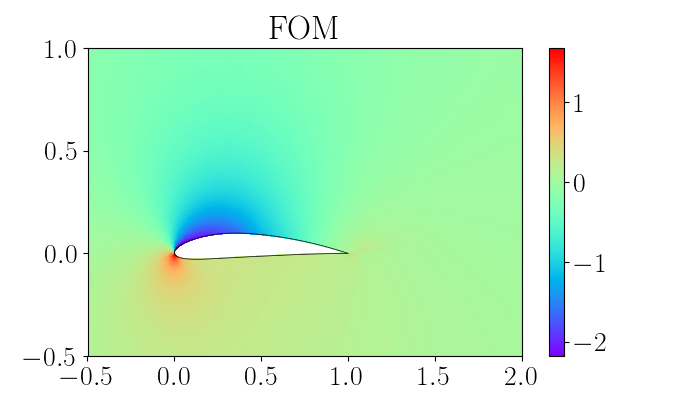}
     \label{case1press118}}
     \subfloat[FOM pressure ($Re\approx \num{518000}$)]
     {\includegraphics[width=0.25\textwidth, trim={2cm 0.85cm 1cm 1.1cm}, clip]{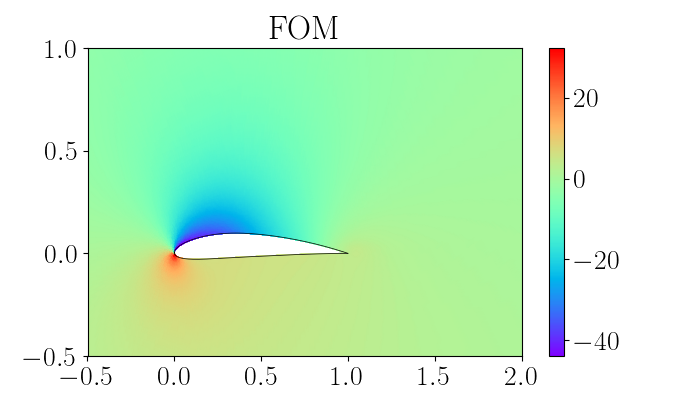}
     \label{case1press518}}
     \subfloat[FOM velocity ($Re\approx \num{118000}$)]{\includegraphics[width=0.25\textwidth, trim={2cm 0.85cm 1cm 1.1cm}, clip]{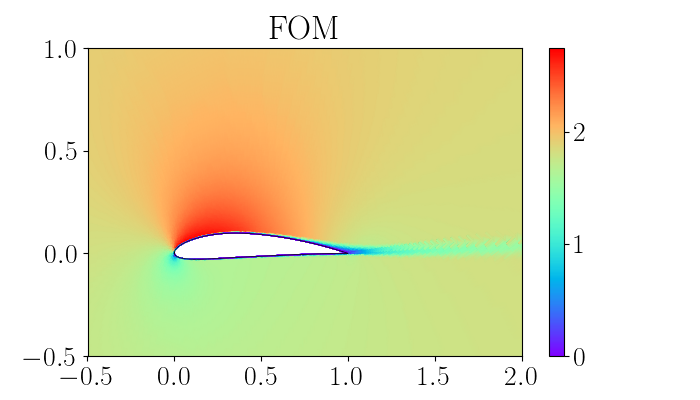}
     \label{case1vel118}}
    \subfloat[FOM velocity ($Re\approx \num{518000}$)]{\includegraphics[width=0.25\textwidth, trim={2cm 0.85cm 1cm 1.1cm}, clip]{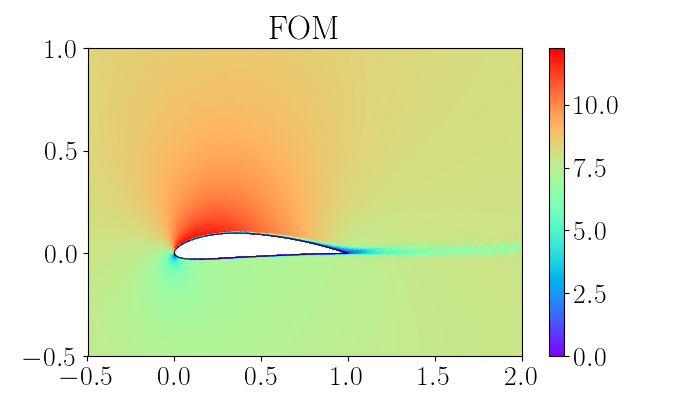}
     \label{case1vel518}}
    \caption{Examples of pressure and velocity magnitude snapshots for the first test case, for $Re\approx \num{118000}$ and $Re\approx \num{518000}$.}
    \label{fig:snapshots-case1}
\end{figure}

\begin{figure}
    \centering
    \subfloat[$-C_p$ for $y>0$]{\includegraphics[width=.5\textwidth]{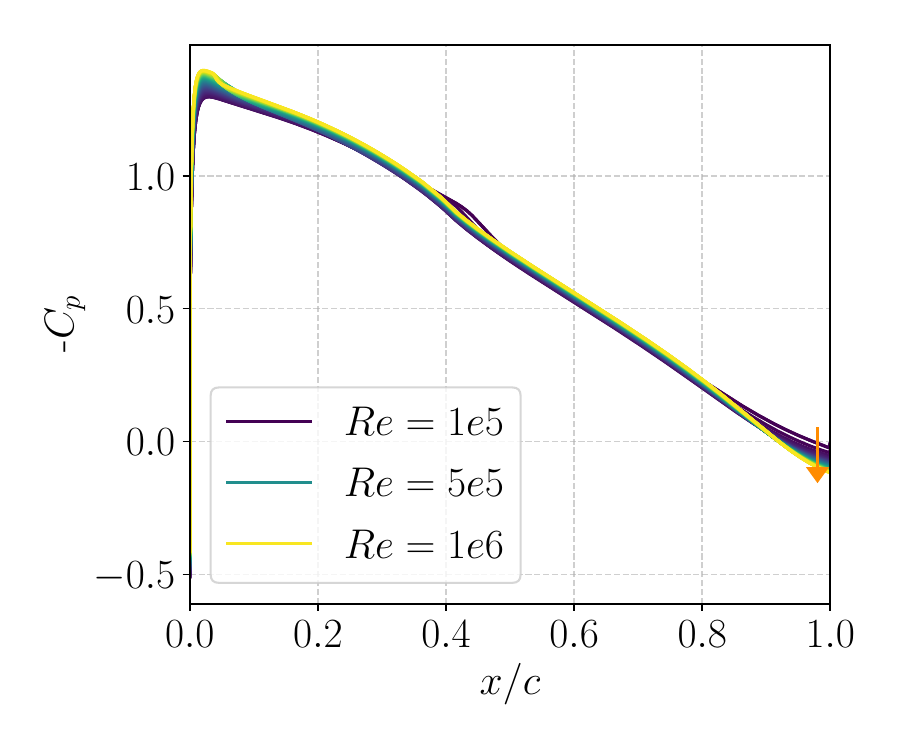}}
    \subfloat[$-C_p$ for $y<0$]{\includegraphics[width=.5\textwidth]{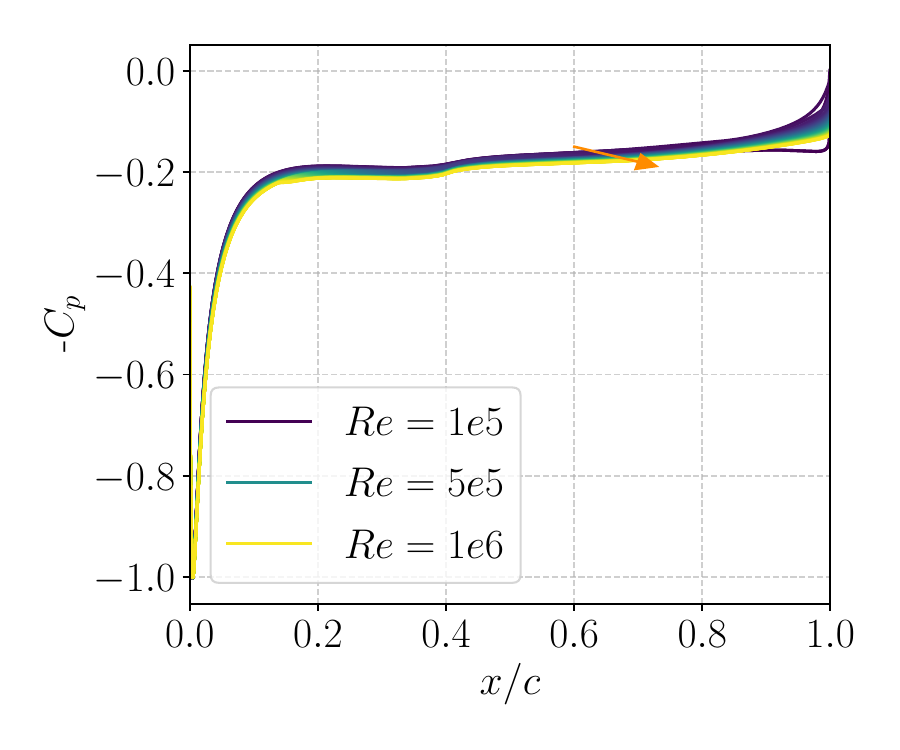}}
    \caption{Pressure coefficient for the elements of the database considered for the first test case, for both suction and pressure sides of the airfoil.}
    \label{fig:cp_snap_case1}
\end{figure}

In order to analyse the flow behavior, in Figure \ref{fig:snapshots-case1} we report a zoomed detail around the airfoil of the pressure and velocity magnitude snapshots, for two values of the Reynolds number. From the snapshots' representation in Figure \ref{fig:snapshots-case1}, we can notice that the solutions represented by the snapshots are all characterized by a similar behavior. Moreover, also the trend of the pressure coefficient in this test case is similar for all the parameters taken into account, as represented in Figure \ref{fig:cp_snap_case1}.

For the above-mentioned reasons, non-intrusive ROMs are particularly suitable to approximate this kind of dataset, as we will see in Section \ref{subsec:rom-result}.


The FOM dataset in this first test case is inspired by the AirfRANS database presented in \cite{bonnet2022airfrans}, where the reader can find more details about the FOM and the turbulence setting.
The above-mentioned FOM has also been validated with the NASA DNS and RANS high-fidelity simulations.

\begin{figure}[htpb!]
    \centering
    \subfloat[]{\includegraphics[width=0.48\textwidth, trim={0cm 0 1cm 0}, clip]{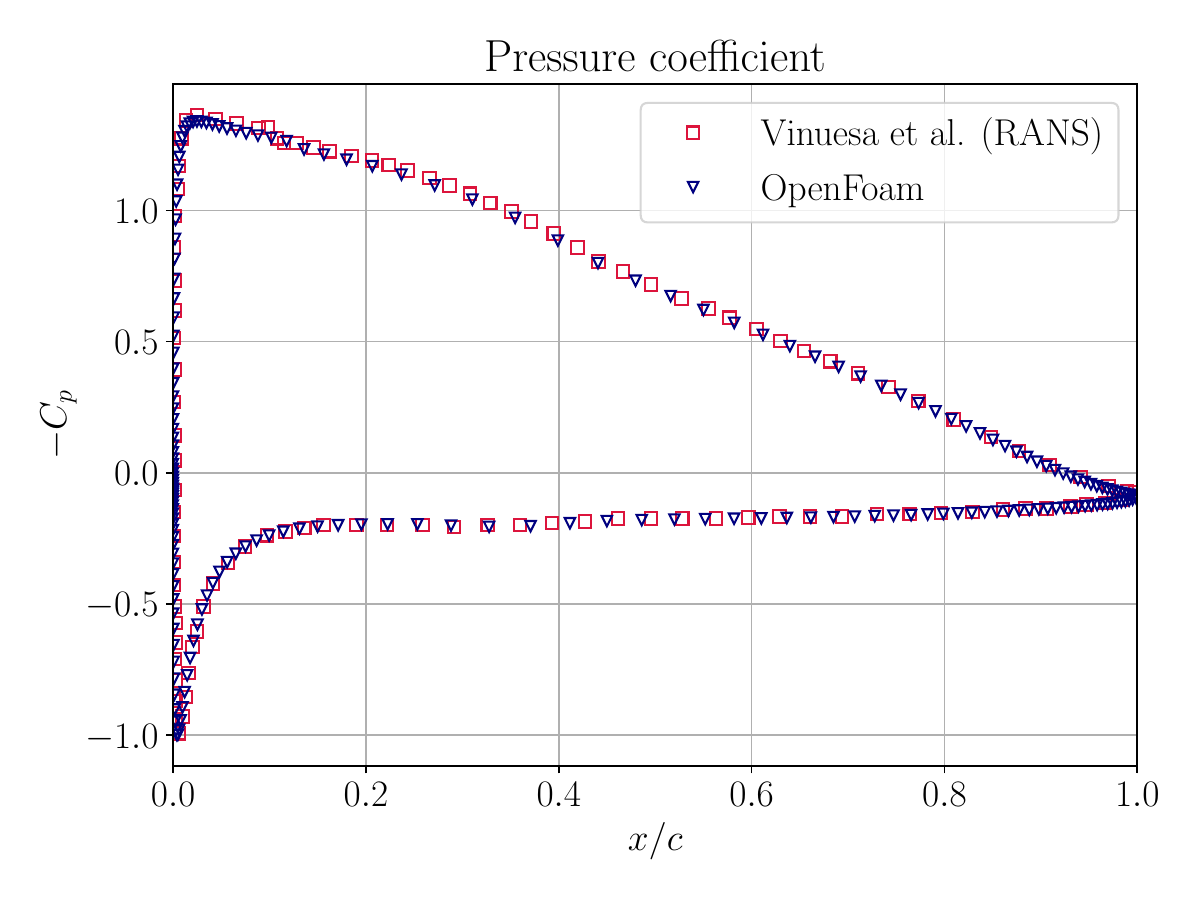} \label{fig1a}}
    \subfloat[]{\includegraphics[width=0.48\textwidth, trim={0cm 0 1cm 0}, clip]{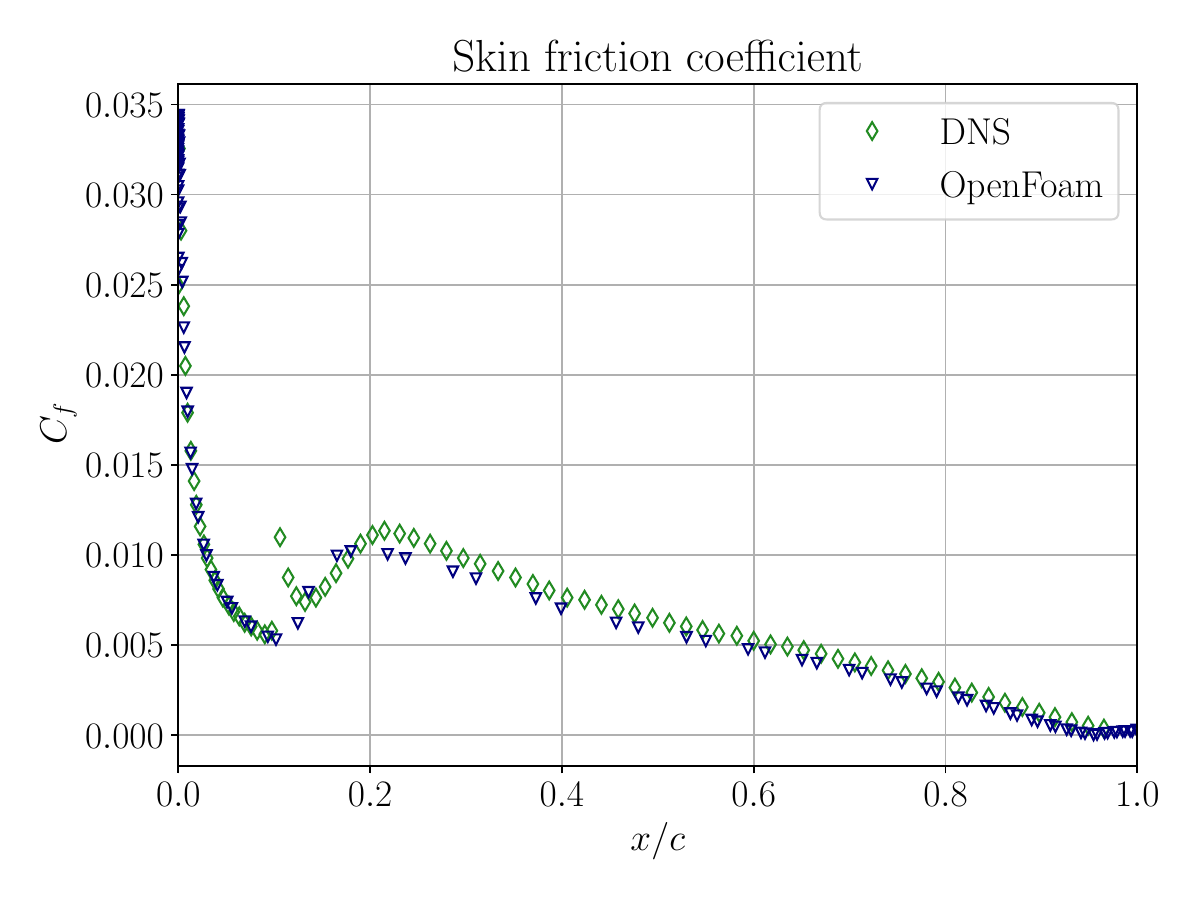}\label{fig1b}}
    \caption{Validation of pressure \protect\subref{fig1a} and upper skin friction \protect\subref{fig1b} coefficients for the first type of full order simulation. The FOM validation is for NACA 4412 at $\mathrm{Re}=4e5$ and angle of attack $\alpha=5^{\circ}$. The reference simulations are taken from \cite{vinuesa2017pressure} (DNS) and \cite{tabatabaei2022rans} (RANS).}
    \label{fig:validation-case-1}
\end{figure}

Indeed, in Figure \ref{fig:validation-case-1} we report the pressure and skin friction coefficients for the NACA $4412$ airfoil at $\mathrm{Re}=4e5$ and angle of attack $\alpha=5^{\circ}$.
In addition to the FOM data, we reported the DNS data~\cite{vinuesa2017pressure} and the RANS simulation data provided in~\cite{tabatabaei2022rans}. Similar validation studies have been performed in other works such as~\cite{tober2018evaluation} and \cite{tabatabaei2022techniques}, where the tripping point is also taken into account in the simulation setup.
We can observe a quite good agreement between the performed RANS and the simulated data. The main difficulty for this specific test case is the considerably high Reynolds number. Due to the predominance of advection effects in this specific set-up, we can expect a particularly challenging test case in terms of model order reduction.

\subsubsection{Test case 2}
\label{subsubsec:fom-case2}
The second test case is the compressible flow past a NACA 0012 airfoil, with fixed inlet velocity at $Re=\num{1e7}$ having as a parameter the angle of attack.

\begin{figure}[htpb!]
    \centering
    \subfloat[Domain and notation]{\includegraphics[width=0.5\textwidth, trim={5cm 0.5cm 3cm 0}, clip]{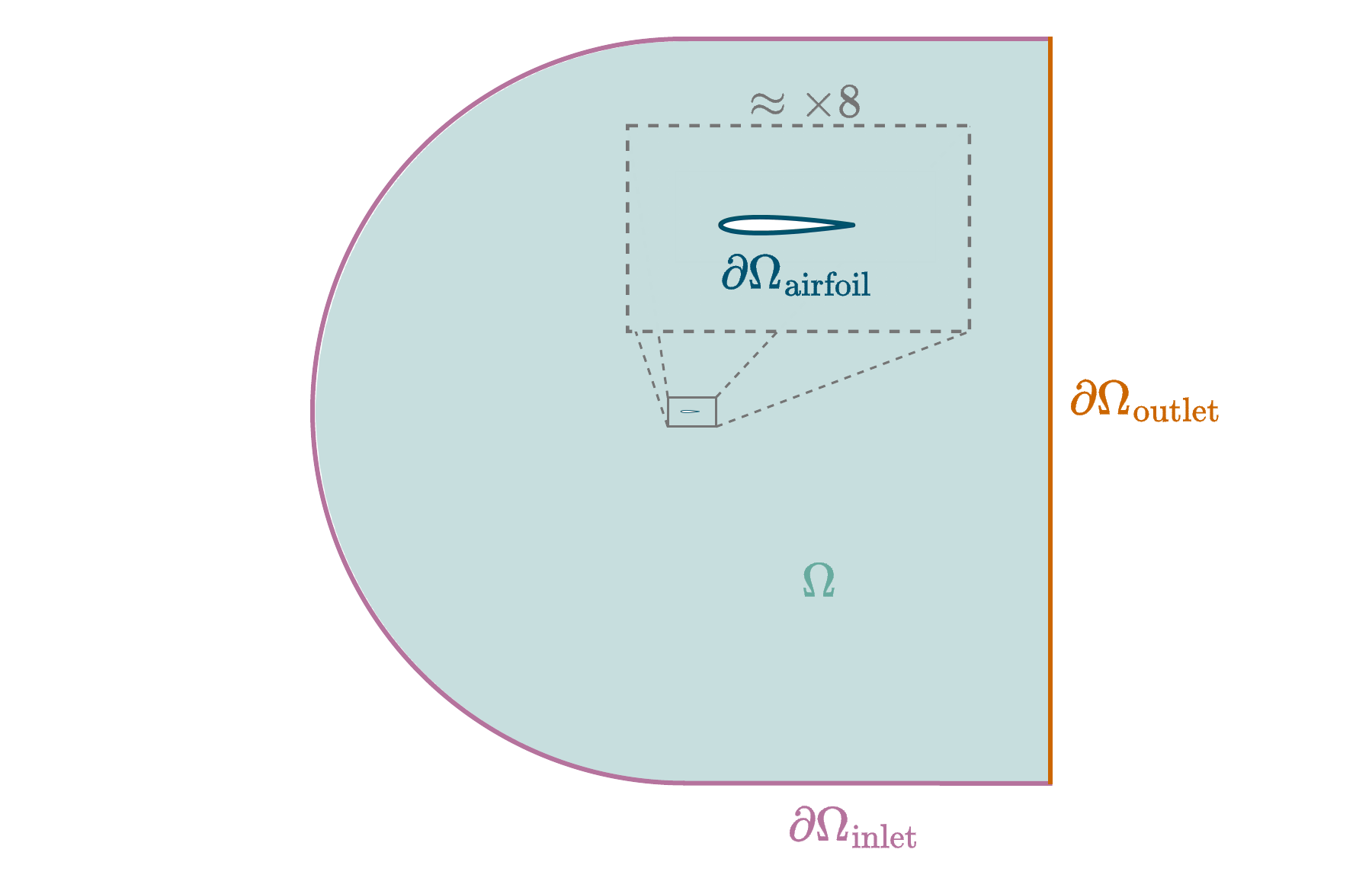}
    \label{case2domain}}
     \subfloat[Mesh representation]{\includegraphics[width=0.5\textwidth, trim={9cm 0cm 13cm 2cm}, clip]{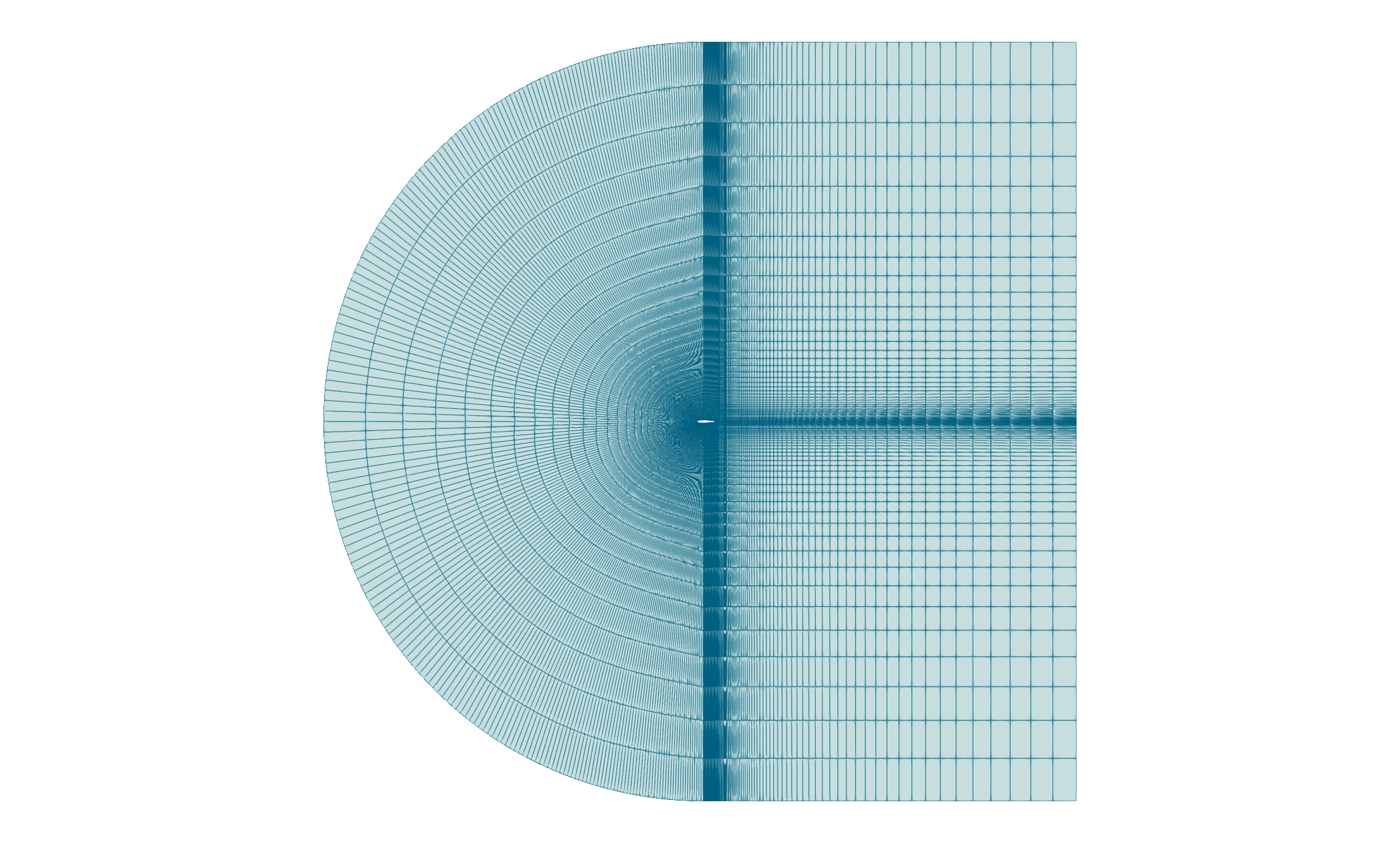}
     \label{case2mesh}}\\
     \subfloat[Zoomed mesh around the airfoil]{\includegraphics[width=0.3\textwidth]{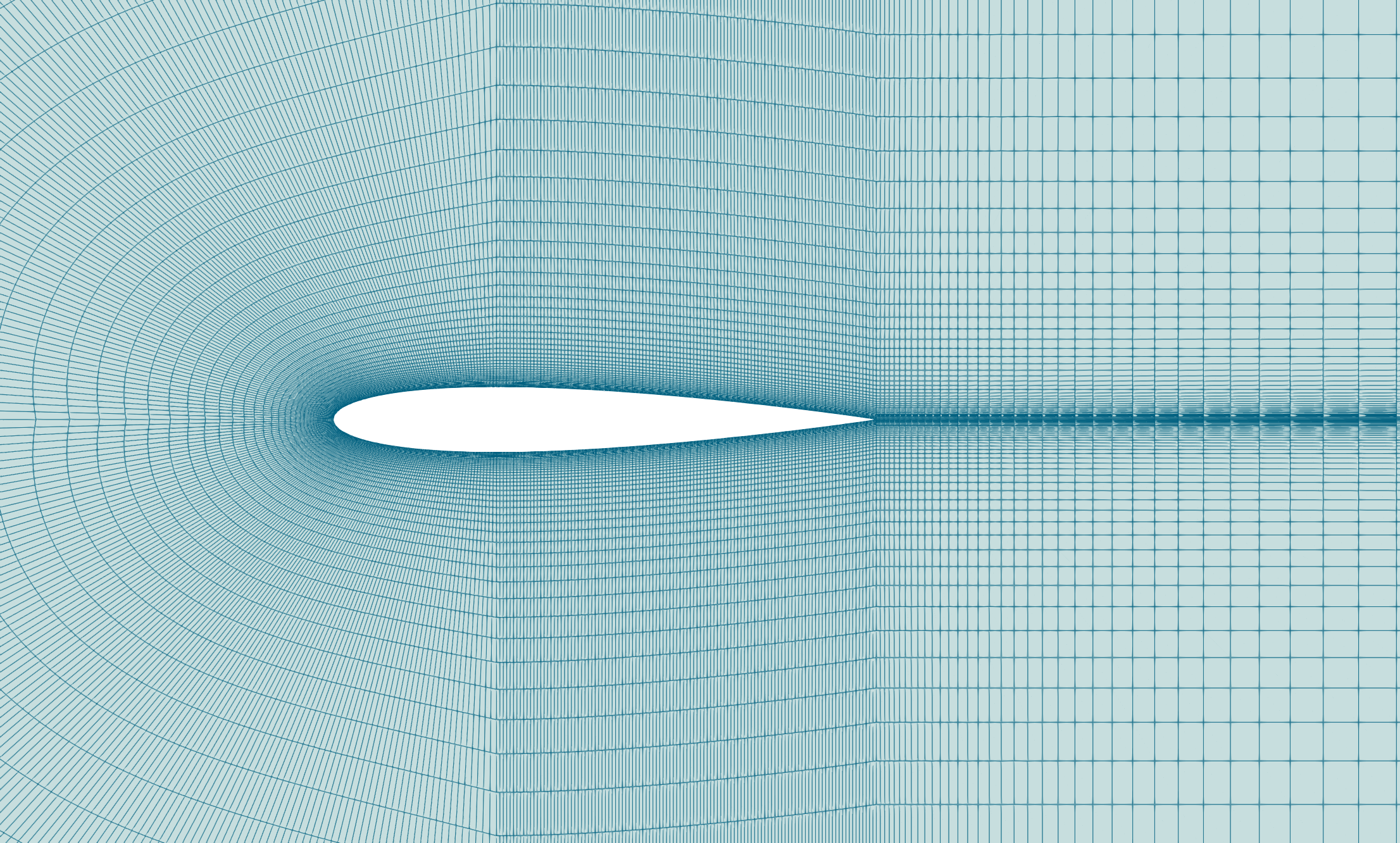}
     \label{case2mesh-zoom}}
    \caption{Domain (\protect \subref{case2domain}) and mesh (\protect \subref{case2mesh}) representation for test case 2, with a zoomed representation for the zone around the airfoil (\protect \subref{case2mesh-zoom}).}
    \label{fig:domain-case2}
\end{figure}

The domain, with the corresponding notation and mesh, is represented in Figure \ref{fig:domain-case2}. It is important to remark here that the mesh is fixed for all snapshots, what changes is the orientation of the inlet velocity, depending on the angle of attack.

In this case, the boundary conditions read as follows:
\begin{equation}
    \begin{cases}
        \bm{u}=\bm{u}_{\text{freestream}},\\
        p=p_{\text{freestream}},
    \end{cases}
    \text{on }\partial \Omega_{\text{inlet}},
    \begin{cases}
        \bm{u}=\bm{u}_{\text{freestream}},\\
        \dfrac{\partial (\int_{\partial \Omega_{\text{outlet}}} w \cdot p dx)}{\partial t}=0,
    \end{cases}
    \text{on }\partial \Omega_{\text{outlet}},
        \begin{cases}
        \bm{u}=\bm{0},\\
        \nabla p \cdot \bm{n} = \bm{0},
    \end{cases}
    \text{on }\partial \Omega_{\text{airfoil}},
\end{equation}
where the pressure condition at the outlet corresponds to the \emph{wave transmission} boundary condition and $w$ is the wave speed.
The freestream conditions are mixed conditions, as in the first test case, but here the pressure fixed value is $\num{1e5}$ whereas the velocity inlet value has a fixed magnitude $\|\bm{u}_{\text{freestream}}\|\approx \SI{277}{\metre \per \second}$.

The solver for the pressure-velocity coupling used in OpenFOAM is \emph{rhoPimpleFoam}, which coincides with a version of the PIMPLE algorithm for compressible flows. The PIMPLE algorithm is a combination of the SIMPLE and PISO (Pressure Implicit with Splitting of Operator) algorithms, and it is used also for steady-state simulations to obtain more robust convergence towards the solution with respect to the SIMPLE counterpart.

Figure \ref{fig:domain-case2} displays two different snapshots, one for $\alpha=0.3^{\circ}$ and one for $\alpha=4.7^{\circ}$.

\begin{figure}[htpb!]
    \centering
    \subfloat[$-C_p$ for $y>0$]{\includegraphics[width=.5\textwidth]{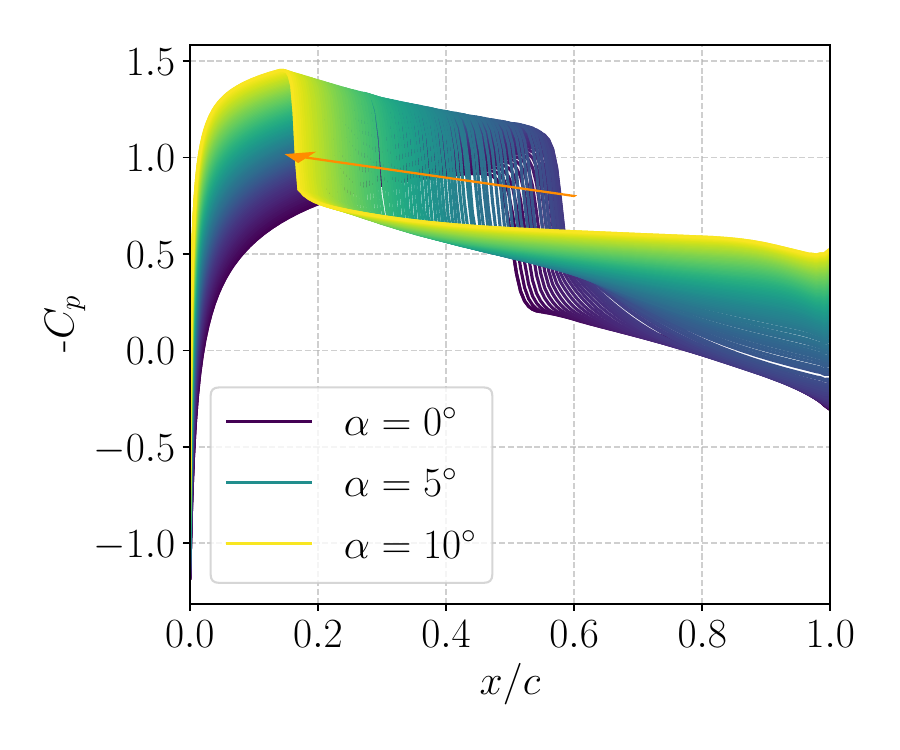}}
    \subfloat[$-C_p$ for $y<0$]{\includegraphics[width=.5\textwidth]{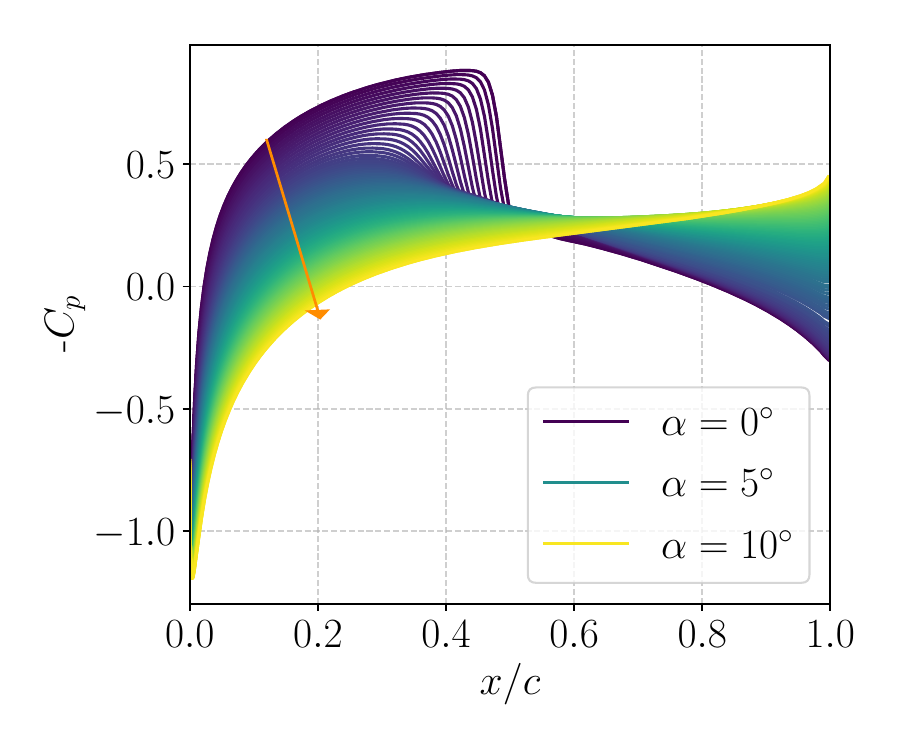}}
    \caption{Pressure coefficient for the elements of the database considered for the second test case, for both suction and pressure sides of the airfoil.}
    \label{fig:cp_snap_case2}
\end{figure}

The difference in the form of the final solution among the snapshots is considerable and can be mainly seen in the shock position, as can be also seen from the trend of the pressure coefficient with respect to the variation of the angle of attack $\alpha$ (Figure \ref{fig:cp_snap_case2}). 
In the comparison between test case 1 (Figure \ref{fig:cp_snap_case1}) and 2 (Figure \ref{fig:cp_snap_case2}), the data variability appears significantly wider in test case 2.
Indeed, we face here an additional challenge due to the presence of shock waves which are well-known complex features to be captured using reduced order modeling. For this particular application, furthermore, the shock location varies significantly within the range of angles of attack herein considered.
In particular, we considered an angle of attack varying between $0^{\circ}$ to $10^{\circ}$. Such a range offers many different behaviors both in terms of boundary layer separation and shock location.
All the previous considerations make this test case more challenging than test case 1 for the characterization through purely data-based ROMs. 

\begin{figure}[htpb!]
    \centering
     \subfloat[FOM pressure ($\alpha=4.7^{\circ}$)]{\includegraphics[width=0.25\textwidth, trim={2cm 0.85cm 1cm 1.1cm}, clip]{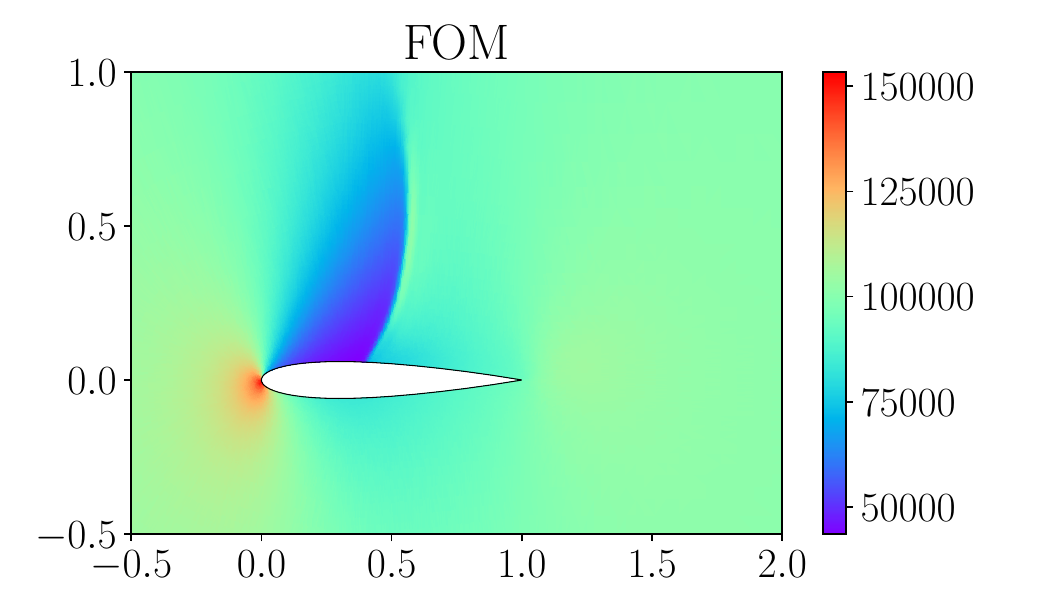}
     \label{case2press-alpha47}}
     \subfloat[FOM pressure ($\alpha=0.3^{\circ}$)]{\includegraphics[width=0.25\textwidth, trim={2cm 0.9cm 1cm 1.1cm}, clip]{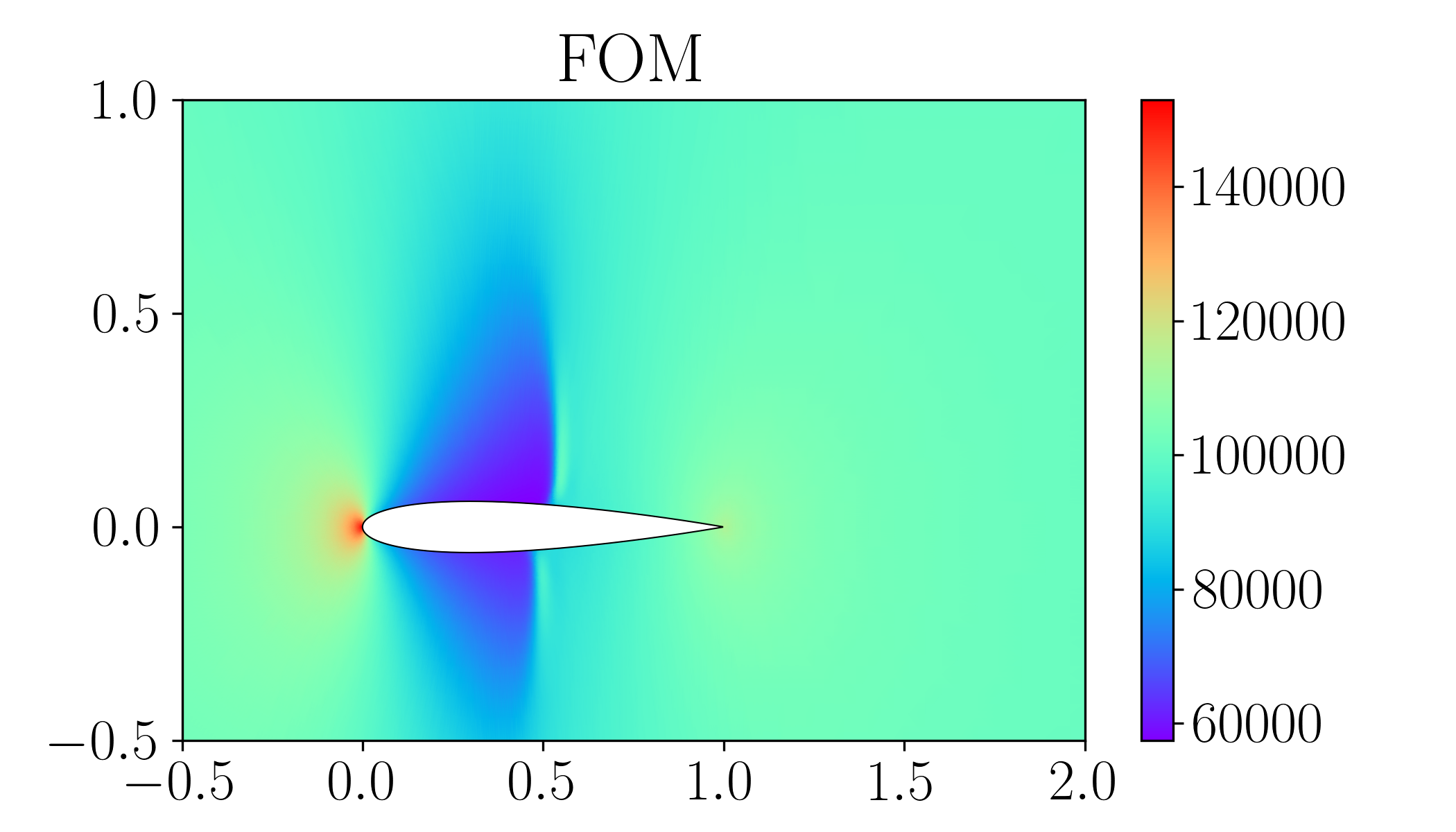}
     \label{case2press-alpha03}}
     \subfloat[FOM velocity ($\alpha=4.7^{\circ}$)]{\includegraphics[width=0.25\textwidth, trim={2cm 0.85cm 1cm 1.1cm}, clip]{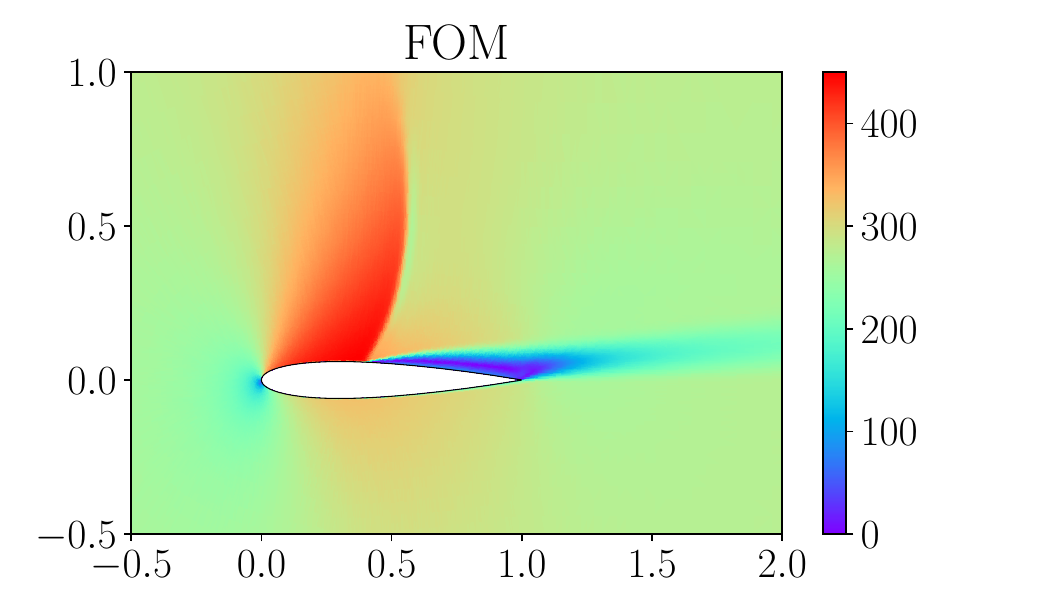}
     \label{case2vel-alpha47}}
          \subfloat[FOM velocity ($\alpha=0.3^{\circ}$)]{\includegraphics[width=0.25\textwidth, trim={2cm 0.9cm 1cm 1.1cm}, clip]{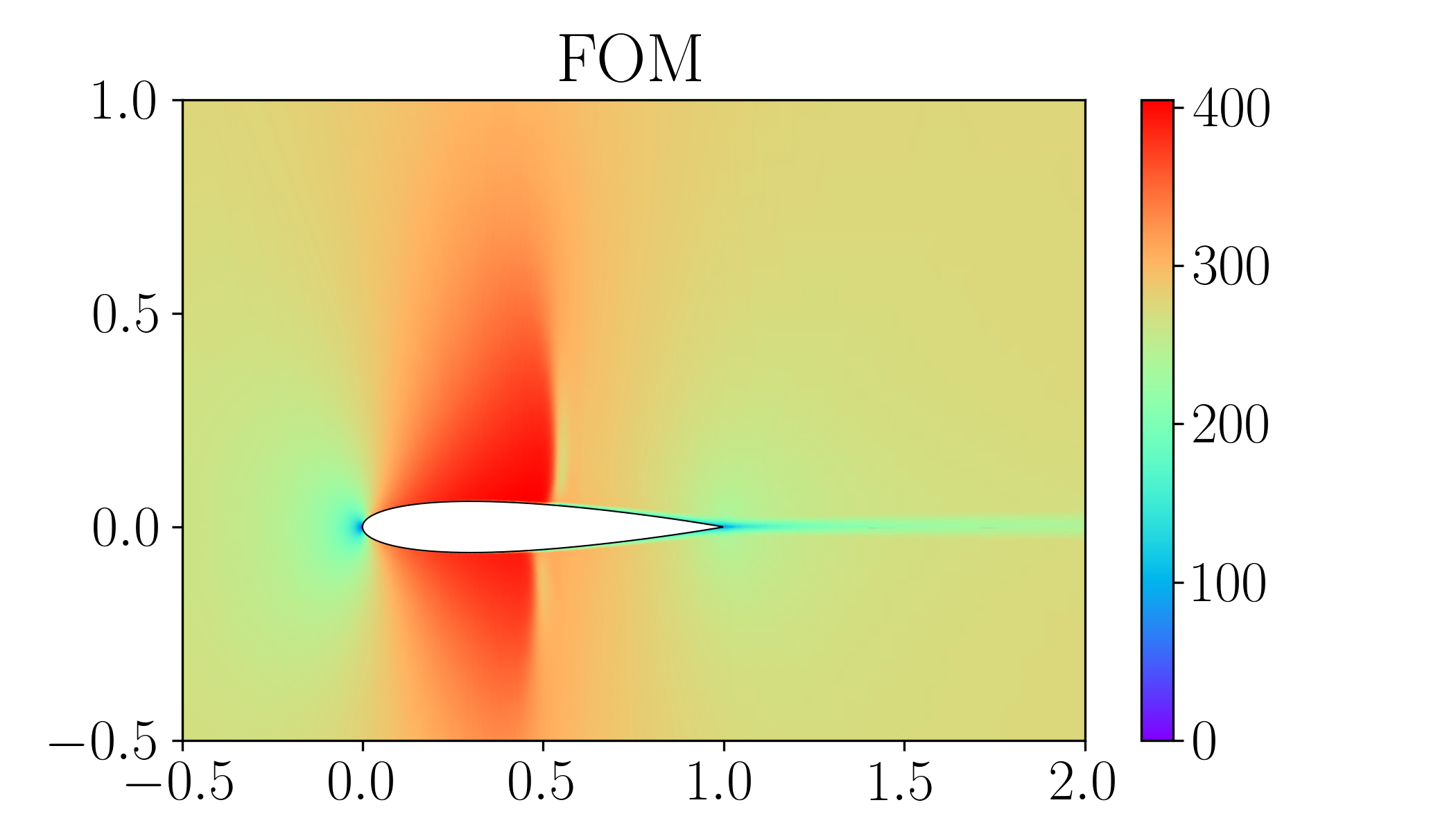}
     \label{case2vel-alpha03}}
    \caption{Examples of FOM snapshots for the second test case, for two values of the angle of attack, $\alpha=0.3$ and $\alpha=4.7$.}
    \label{fig:snapshots-case2}
\end{figure}
\begin{figure}[htpb!]
    \centering
    \subfloat[]{\includegraphics[width=0.48\textwidth, trim={0cm 0 1cm 0}, clip]{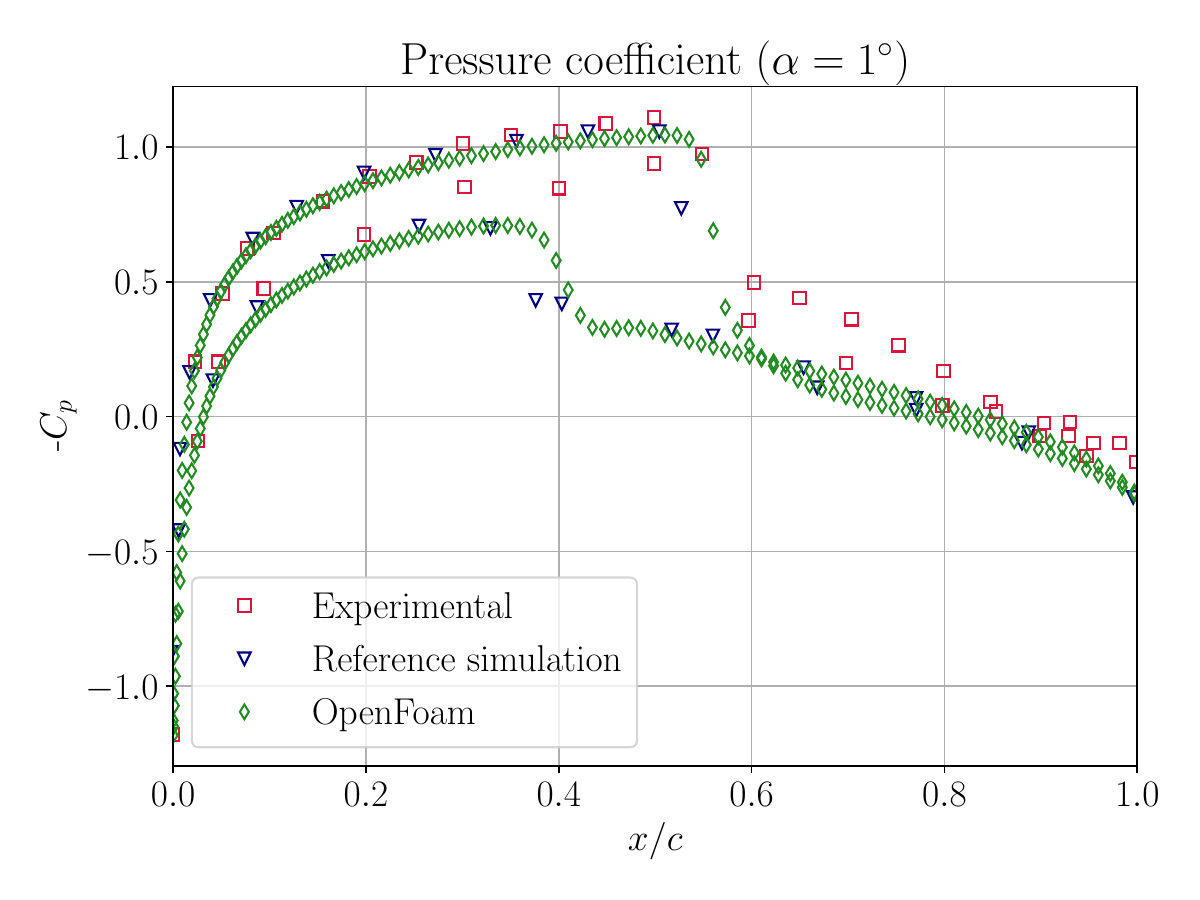} \label{fig2a}}
    \subfloat[]{\includegraphics[width=0.48\textwidth, trim={0cm 0 1cm 0}, clip]{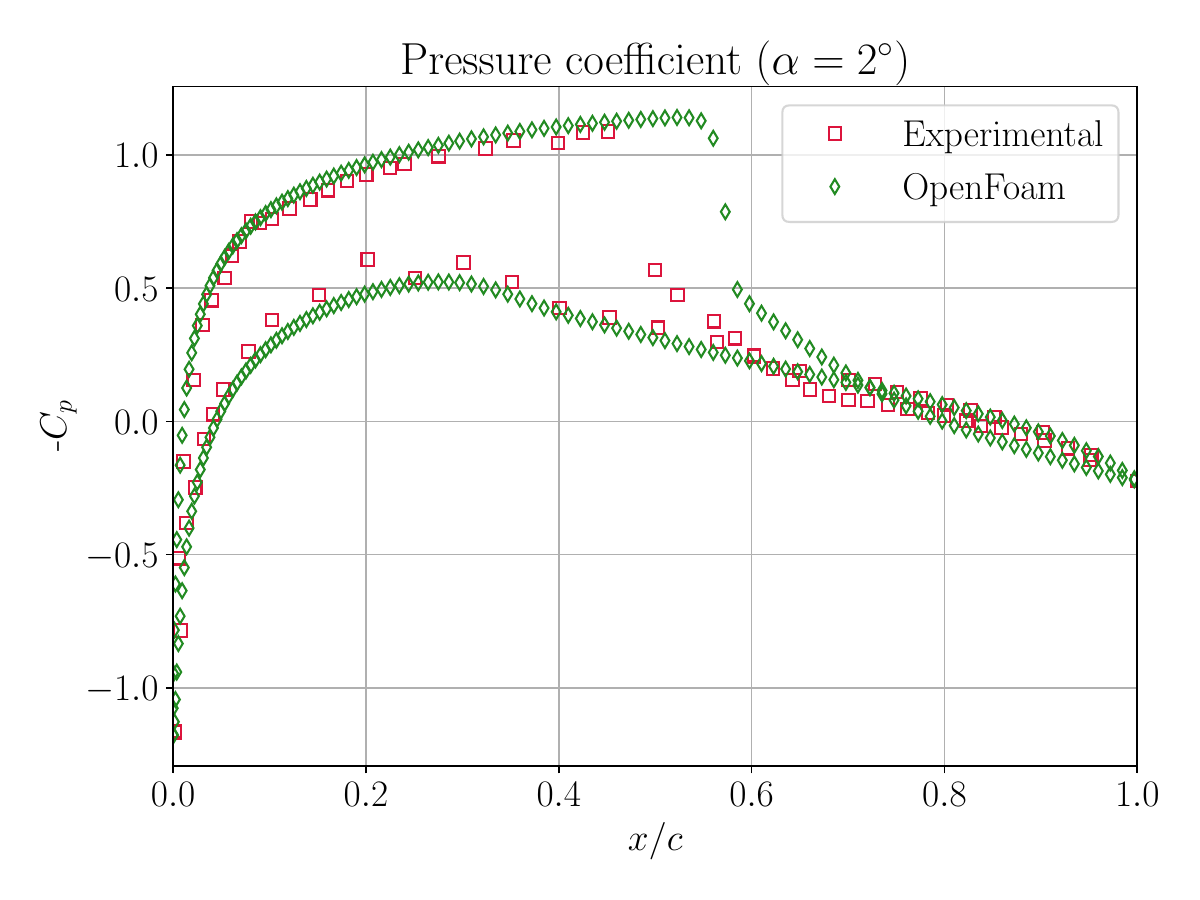}\label{fig2b}}
    \caption{Validation of pressure coefficient of the full order model for two angles of attack: $\alpha=1^{\circ}$ \protect\subref{fig2a} and $\alpha=2^{\circ}$\protect\subref{fig2b}. The FOM validation is for NACA 0012 at Mach $\mathrm{Ma}=0.8$ and the reference results are presented in~\cite{mcdevitt1985static} (experimental) and ~\cite{iovnovich2012reynolds} (wall-resolved reference simulation).}
    \label{fig:validation-case-2}
\end{figure}

A validation analysis has also been conducted for the second test case. Figure \ref{fig:validation-case-2} represents the pressure coefficient in two different operating conditions, namely at $\alpha=1^{\circ}$ and $\alpha=2^{\circ}$.
Here, the results obtained by OpenFoam are compared with the experimental data presented in~\cite{mcdevitt1985static} and a wall-resolved reference simulation~\cite{iovnovich2012reynolds}, which used the same turbulence model.
The agreement is quite good, especially with respect to the simulation taken as reference. The deviations from the experimental data are intrinsically related to the RANS models employed. This set-up has been used precisely because different models can provide quite different outcomes. In this scenario, it would be useful to apply the same type of RANS aggregation proposed by de Zordo-Banliat et al. \cite{de2024space}. However, such improvements in terms of RANS modelling is out of the scope of the present work and future research will be focused on improving the RANS modelisation and combine it with the reduced-order modelling strategies herein proposed.

\subsubsection{Computational time}
In the first test case, the CPU time required for running each full order simulation is in the range $[2 \text{ h }30 \text{ min }, 7 \text{ h}]$, depending on the Reynolds value we consider. In the second test case, the CPU time is in the interval $[1 \text{ h }45 \text{ min }, 3 \text{ h}]$, depending on the angle of attack considered. The FOM simulations are all performed using 4 processor cores on SISSA HPC cluster Ulysses (200 TFLOPS, 2TB RAM, 7000 cores).

\subsection{Individual ROM analysis}
\label{subsec:rom-result}
This part of the paper is dedicated to the analysis of the individual ROMs considered for the model mixture.

The fields considered for the model reduction can be divided into two types: (i) one-dimensional (1D) wall variables on the airfoil (pressure and wall shear stress); (ii) two-dimensional (2D) fields in the internal mesh (pressure and velocity magnitude).

In both the test cases taken into account, we collect the predictions of the following non-intrusive models:
\begin{itemize}
    \item POD-RBF;
    \item POD-GPR;
    \item POD-ANN;
    \item AE-RBF and PODAE-RBF, for the 1D airfoil fields and the 2D mesh fields, respectively;
    \item AE-GPR and PODAE-GPR, for the 1D airfoil fields and the 2D mesh fields, respectively;
    \item AE-ANN and PODAE-ANN, for the 1D airfoil fields and the 2D mesh fields, respectively.
\end{itemize}
As above-mentioned, the nonlinear reduction approach considered for the 2D fields is the mixed technique PODAE, instead of the pure AE. The reason is that the PODAE allows for a valuable reduction in the computational effort in the offline part, and especially in the training of the autoencoder.

Following the notation introduced in Section \ref{subsec:met-aggregation}, the training dataset considered for all the data-based ROMs is composed of $\nsnap=70$ snapshots, the evaluation set of $\neval=20$ elements, and the test set of $\ntest=10$ elements.

\begin{figure}[htpb!]
    \centering
    \subfloat[1D fields on airfoil]{\includegraphics[width=0.5\textwidth]{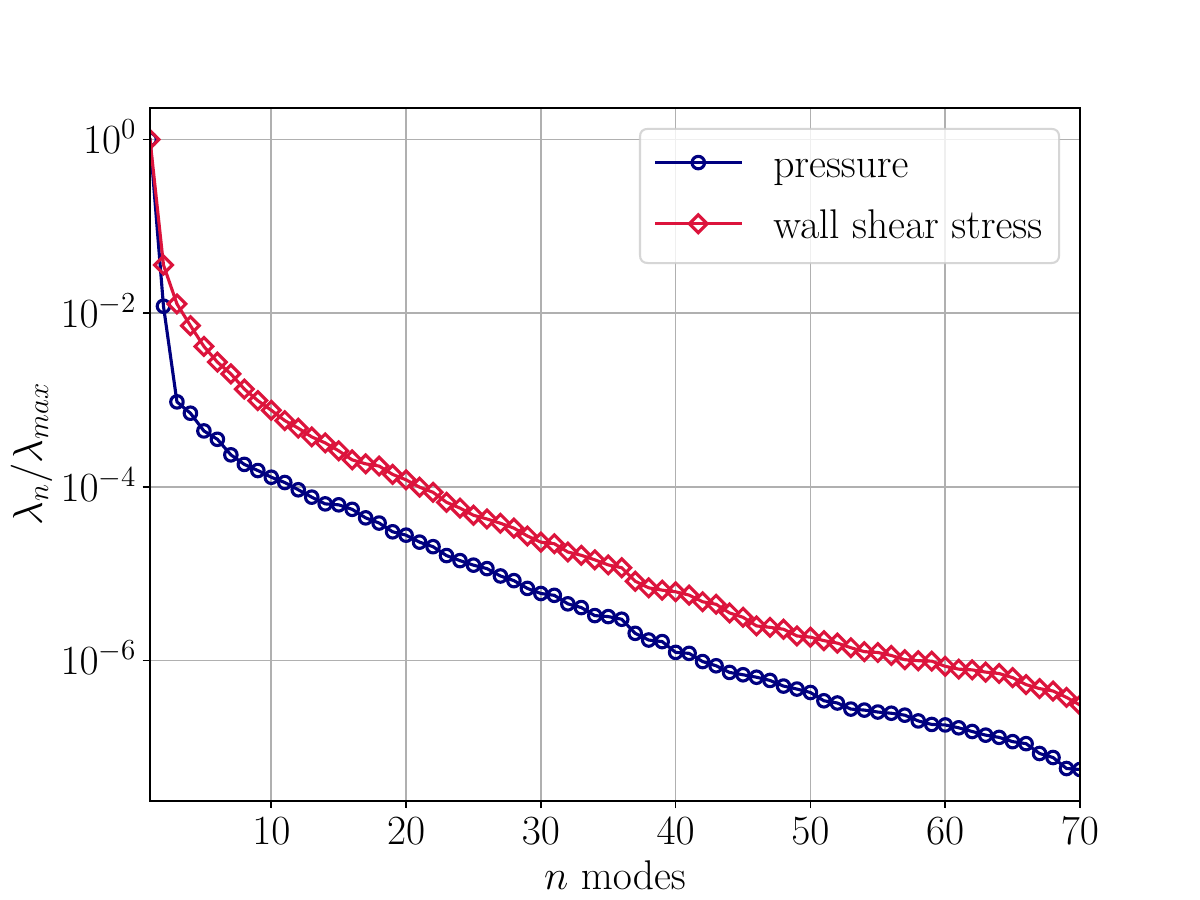}}
    \subfloat[2D fields on mesh]{\includegraphics[width=0.5\textwidth]{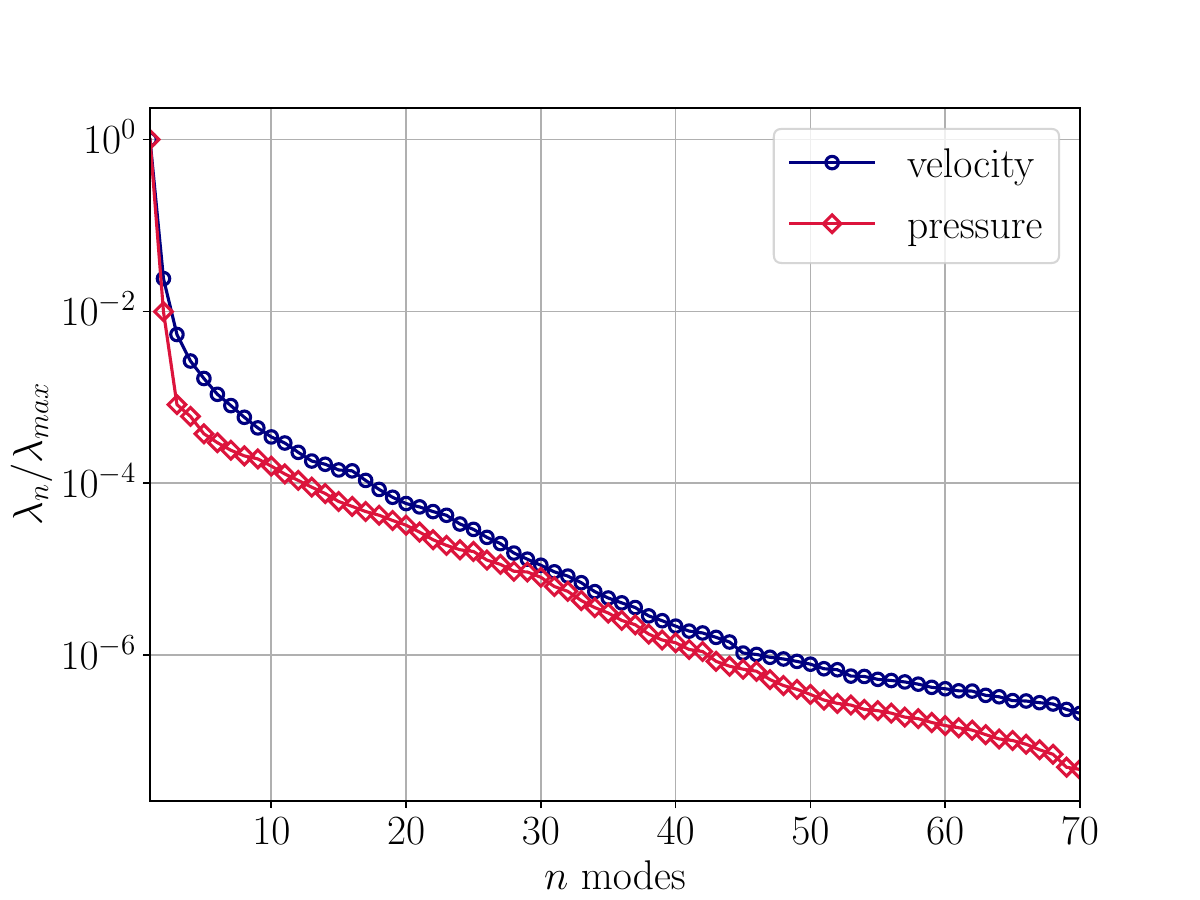}}
    \caption{Decay of POD eigenvalues for the 1D fields on the airfoil and for 2D fields on the computational mesh around the airfoil for test case 1.}
    \label{fig:decay-eig-case1}
\end{figure}

\begin{figure}[htpb!]
    \centering
    \subfloat[1D fields on airfoil]{\includegraphics[width=0.5\textwidth]{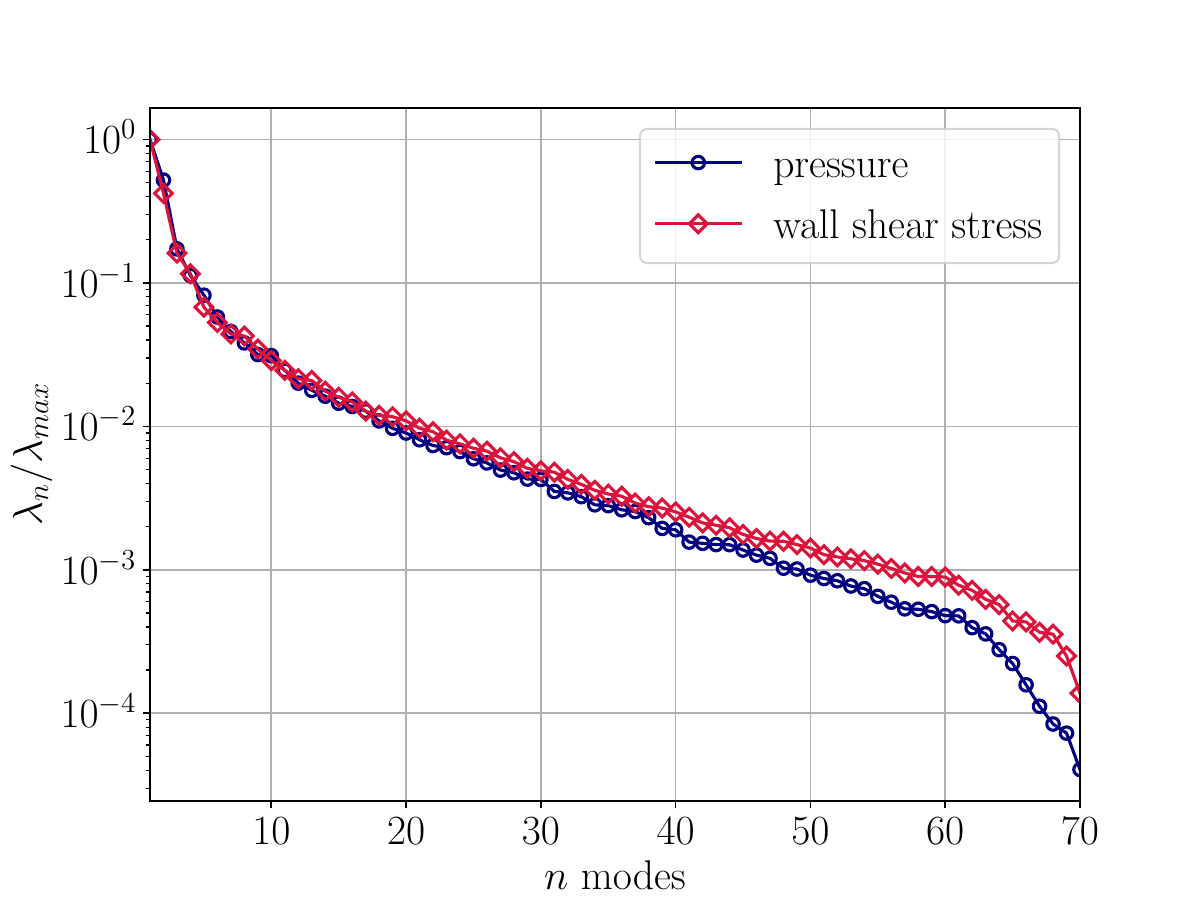}}
    \subfloat[2D fields on mesh]{\includegraphics[width=0.5\textwidth]{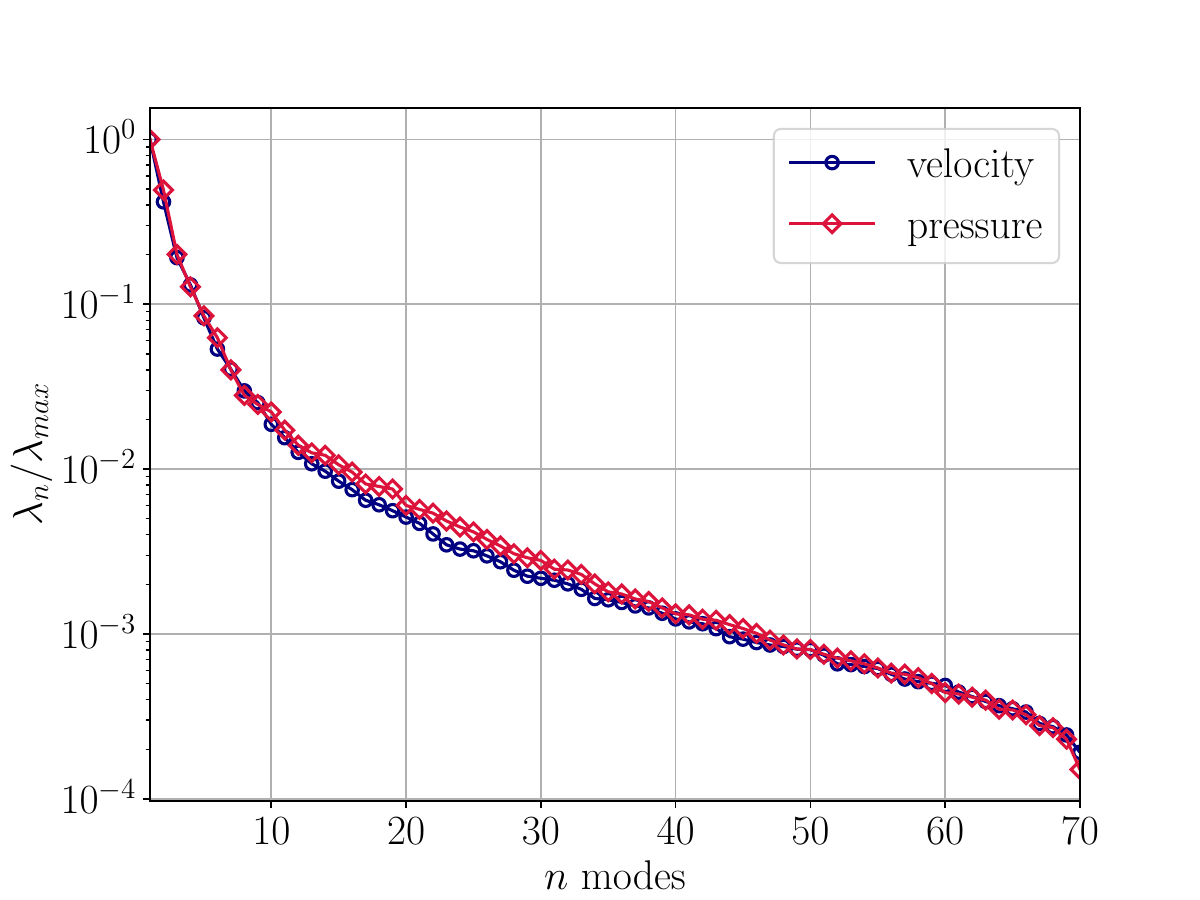}}
    \caption{Decay of POD eigenvalues for the 1D fields on the airfoil and for 2D fields on the computational mesh around the airfoil for test case 2.}
    \label{fig:decay-eig-case2}
\end{figure}

Figures \ref{fig:decay-eig-case1} and \ref{fig:decay-eig-case2} represent the POD eigenvalues' decays for all the variables considered and for the two first cases, respectively.

It can be seen that the eigenvalues' decays for all variables in the second test case is considerably slower than in the corresponding variables in the first test case. This fact is reflected on the accuracy of the ROMs and of the model mixture, that are presented in Sections \ref{subsubsec:aggr-result-1} and \ref{subsubsec:aggr-result-2} for the first and second test case, respectively. The problems herein considered are consequently presented for increasing complexity in terms of model order reduction. 

At this point, we address the accuracy of the individual ROMs on the \emph{ROM test set}, namely the union of the evaluation and test datasets. The two test cases are separately considered in \ref{subsubsec:rom-case1} and \ref{subsubsec:rom-case2}.

\subsubsection{Test case 1}
\label{subsubsec:rom-case1}

\begin{figure}[htpb!]
    \centering
    \subfloat[Latent dimension$=3$]{\includegraphics[width=.5\textwidth]{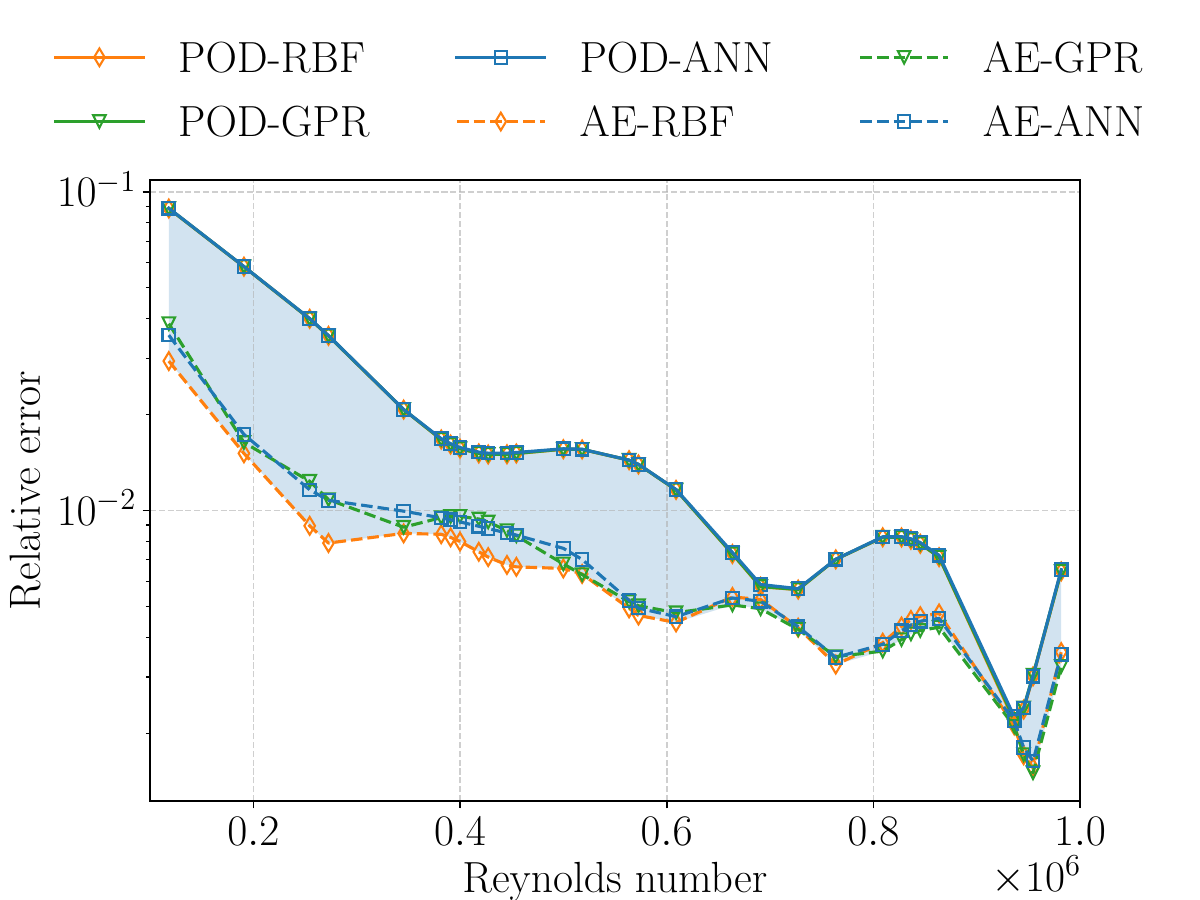}}
    \subfloat[Latent dimension$=10$]{\includegraphics[width=.5\textwidth]{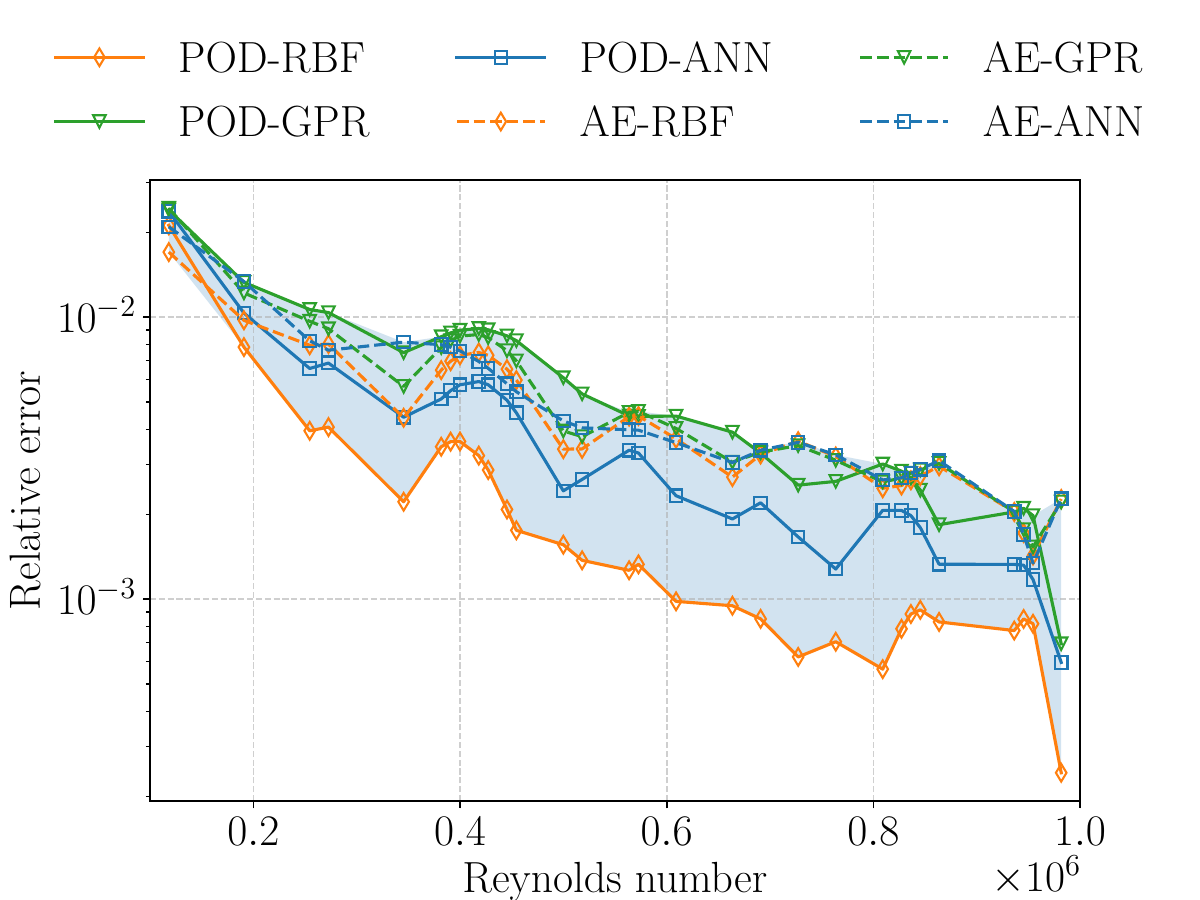}}
    \caption{Relative test errors of individual ROMs for the 1D pressure field, for two different latent dimensions, namely 3 and 10. The results are for test case 1.}
    \label{fig:test-w-err-case-1}
\end{figure}

\begin{figure}[htpb!]
    \centering
    \subfloat[Latent dimension$=3$]{\includegraphics[width=.5\textwidth]{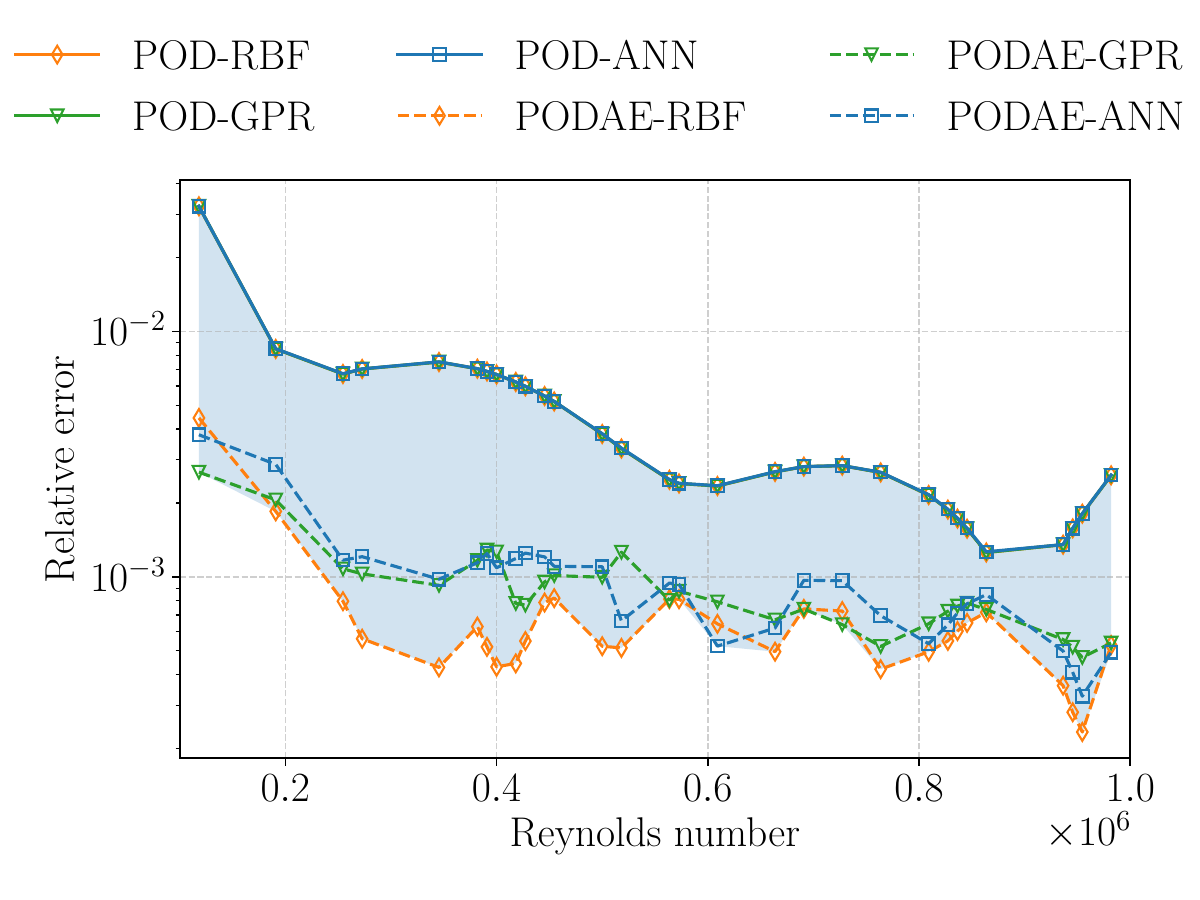}}
    \subfloat[Latent dimension$=10$]{\includegraphics[width=.5\textwidth]{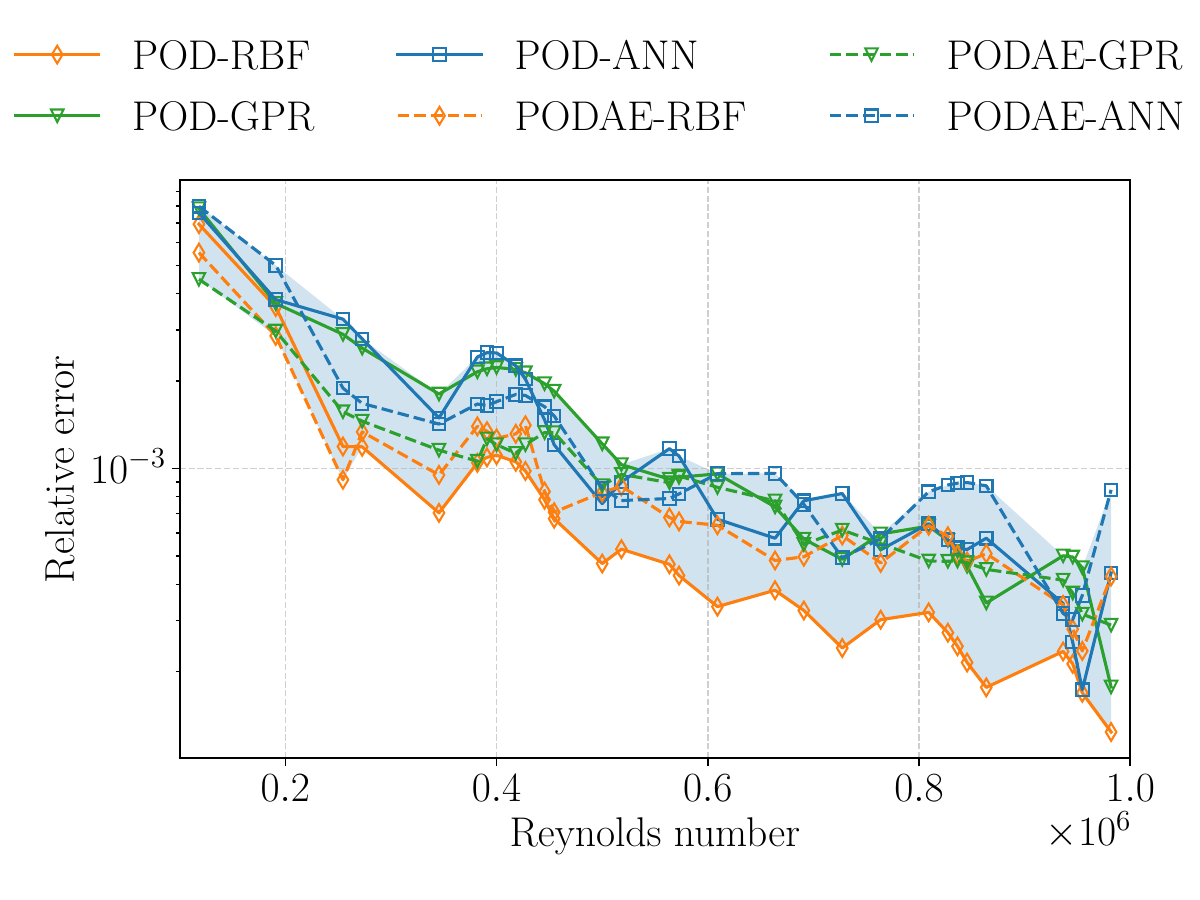}}
    \caption{Relative test errors of individual ROMs for the 2D velocity field, for two different latent dimensions, namely 3 and 10. The results are for test case 1.}
    \label{fig:test-vint-err-case-1}
\end{figure}

In Figures \ref{fig:test-w-err-case-1} and \ref{fig:test-vint-err-case-1} we represent the relative errors on the ROM test set with respect to the FOM counterpart, on the 1D wall shear stress field and on the 2D velocity magnitude field, respectively.

From both the Figures the following considerations can be drawn:
\begin{itemize}
    \item all the individual ROMs in both the latent dimensions considered are characterized by better performances for high values of the Reynolds parameter;
    \item the AE-based ROMs have similar accuracy for the two latent dimensions;
    \item the AE outperforms the results of the POD as a reduction approach almost everywhere in the parametric domain when the latent dimension is small. Indeed, in this case we are in the so-called \emph{under-resolved} modal regime, namely the POD modes are not able to fully characterized the dynamics of the system.
\end{itemize}

As a comparison between the two fields (1D wall shear and 2D velocity), the results are less accurate for the wall shear stress field, as can be seen from the order of magnitude of the errors. Indeed, it is usually a complex field to characterize in the flow past an airfoil.

For what concerns the computational time, all the ROMs are performed on an Intel(R) Core(TM) i5-4570 CPU @ 3.20GHz 16GB RAM on only one processor core

\subsubsection{Test case 2}
\label{subsubsec:rom-case2}

\begin{figure}[htpb!]
    \centering
    \subfloat[Latent dimension$=3$]{\includegraphics[width=.5\textwidth]{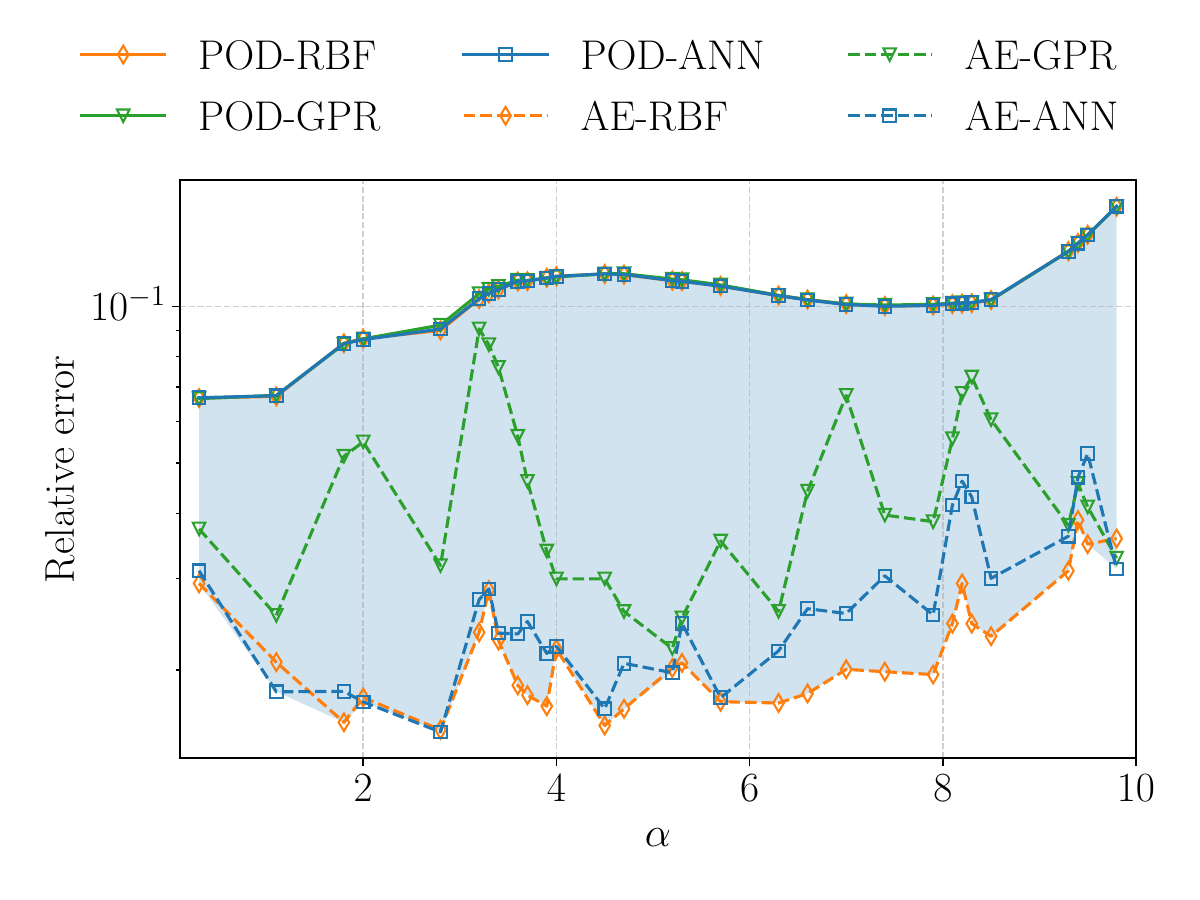}}
    \subfloat[Latent dimension$=10$]{\includegraphics[width=.5\textwidth]{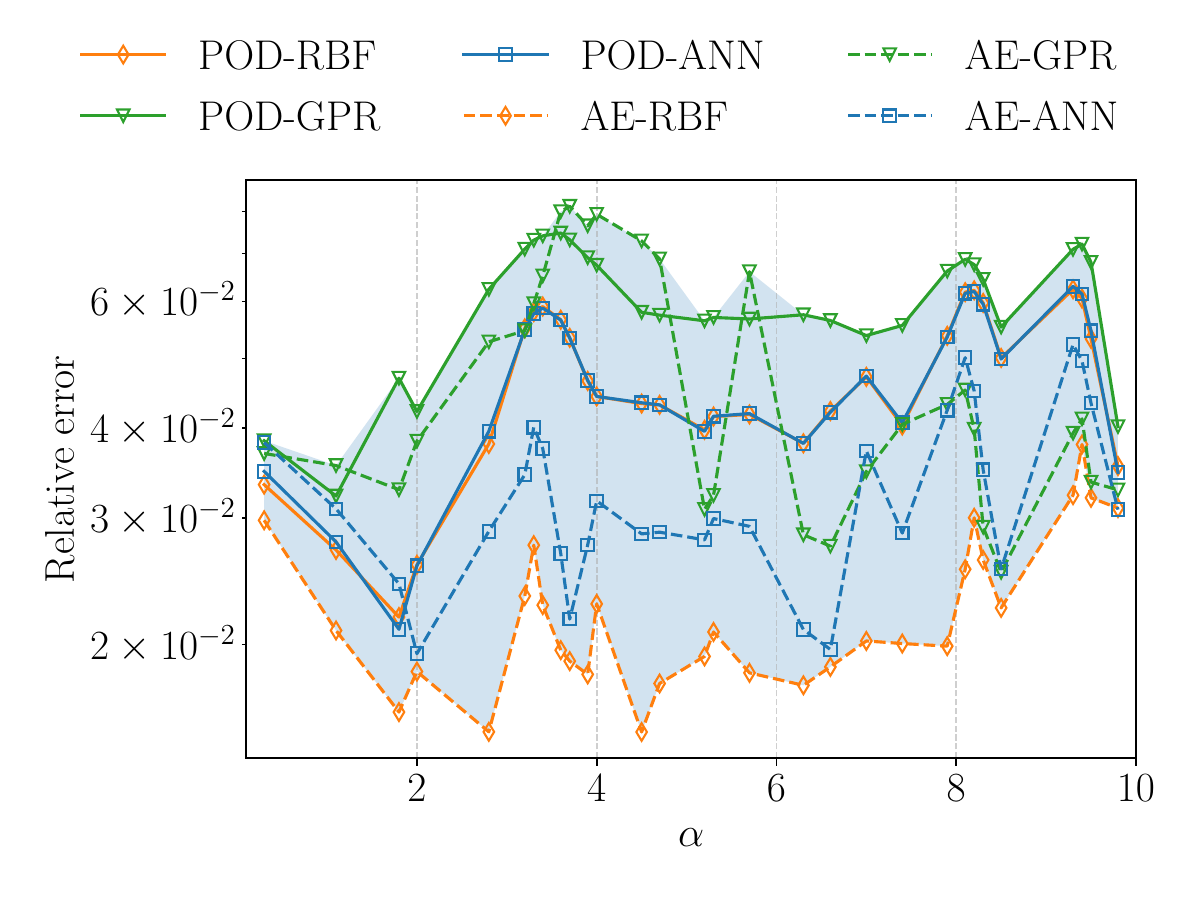}}
    \caption{Relative test errors of individual ROMs for the 2D pressure field, for two different latent dimensions, namely 3 and 10. The results are for test case 2.}
    \label{fig:test-w-err-case-2}
\end{figure}

\begin{figure}[htpb!]
    \centering
    \subfloat[Latent dimension$=3$]{\includegraphics[width=.5\textwidth]{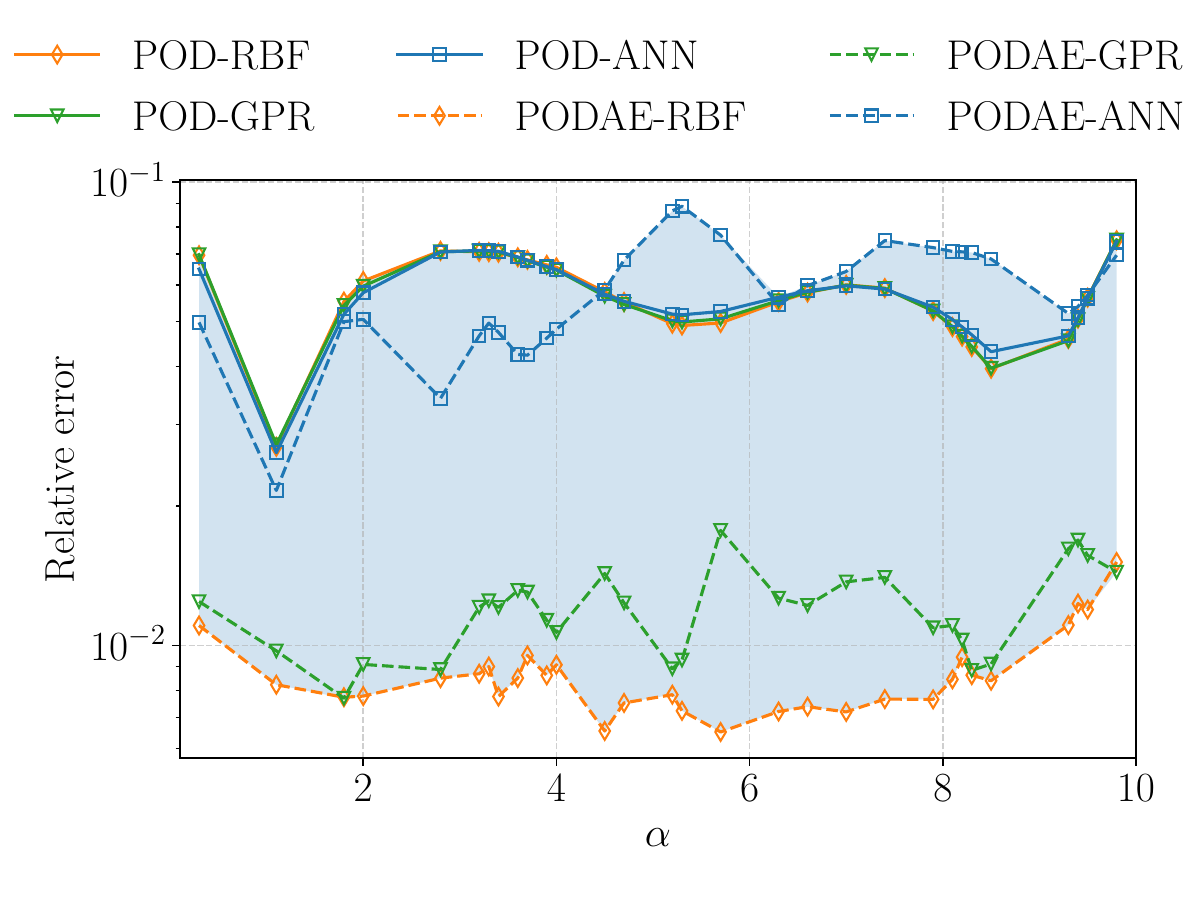}}
    \subfloat[Latent dimension$=10$]{\includegraphics[width=.5\textwidth]{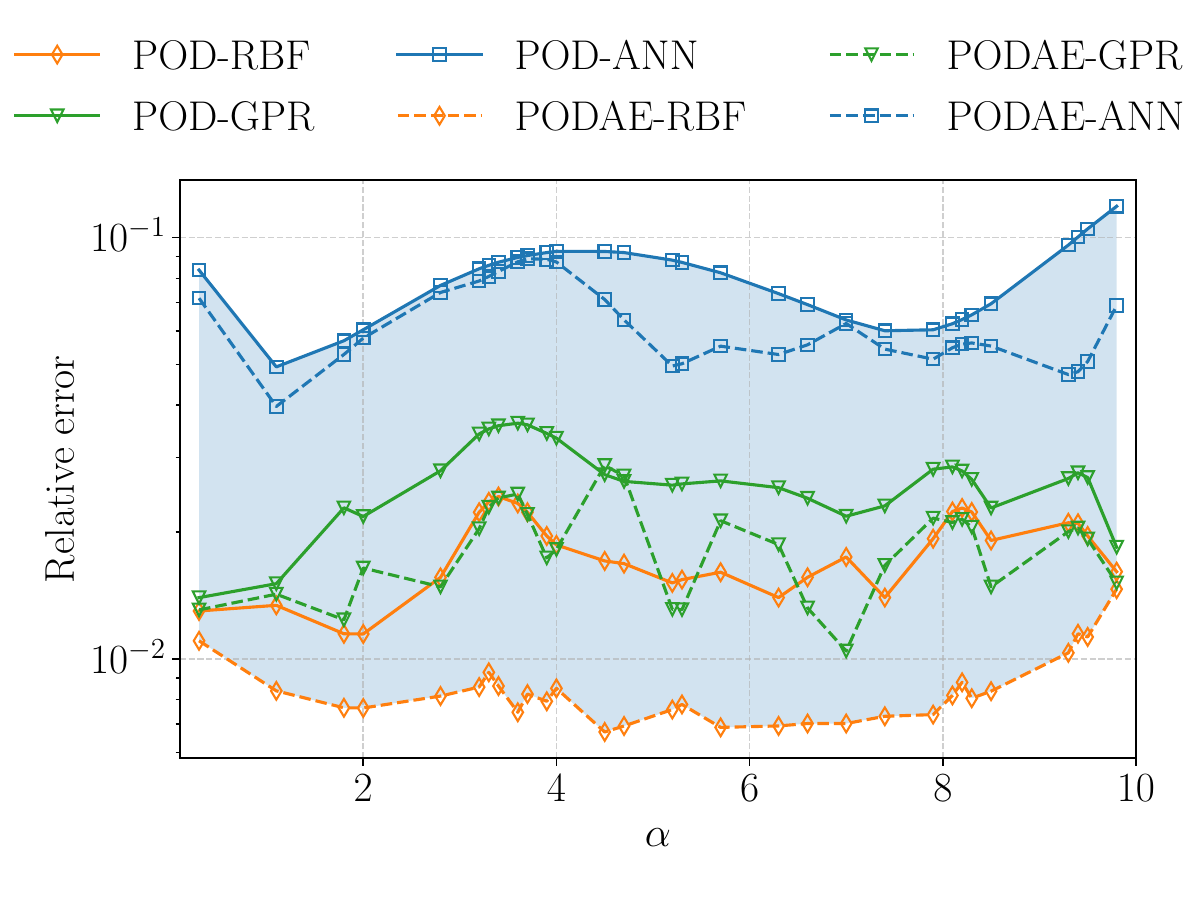}}
    \caption{Relative test errors of individual ROMs for the 2D velocity field, for two different latent dimensions, namely 3 and 10. The results are for test case 2.}
    \label{fig:test-vint-err-case-2}
\end{figure}

In this second test case, we represent the relative test errors for the 1D wall shear stress and the 2D pressure fields, in Figures \ref{fig:test-w-err-case-2} and \ref{fig:test-vint-err-case-2} respectively.

As expected, the performances of the individual ROMs are remarkably different from the first test case.
We can, hence, notice that:
\begin{itemize}
    \item the performances of the single ROMs are more differentiated than in the first test case, and, therefore, the accessible region by all the ROMs is bigger;
    \item there is not a clear trend of the errors as the parameter (the angle of attack in this case) increases;
    \item the AE (or PODAE)-RBF model outperforms all the other techniques in most of the parameter's values;
    \item while in the first test case, the choice of the reduction is the most relevant in the ROM performance, in this test case also the approximation type plays an important role.
\end{itemize}

The complexity of this problem can be quantified by observing the decay of the singular values of the snapshots' matrix for pressure coefficient and wall shear stress already depicted in figure \ref{fig:decay-eig-case2}. In order to achieve an error of $1\%$, in fact, at least $20$ modes are needed.

In order to have a large variability between the different reduced-order models herein considered, we decided to consider a relatively small number of latent variables (namely $3$ and $10$). The motivation for this choice is twofold: first, we obviously want the latent space to be small so that we can achieve a significant speed-up by the use of ROMs. Secondly, we want to position ourselves in a case where the different ROMs provide significantly different results. Only in this scenario, the aggregation can be beneficial.

Subsection \ref{subsec:aggr-result} will thoroughly discuss how the different decays of the two cases will affect both the individual ROMs and their aggregation on the overall accuracy of the proposed approach.

\medskip

\subsubsection{Computational time}
All the ROMs are performed on an Intel(R) Core(TM) i5-4570 CPU @ 3.20GHz 16GB RAM on only one processor core. Tables \ref{tab:rom-time} and \ref{tab:rom-time-online} report the execution offline and online times for the non-intrusive ROM models, respectively. It is worth to specify that the offline times do not include the FOM times, which are reported in Sections \ref{subsubsec:fom-case1} and \ref{subsubsec:fom-case2}, for the two test cases respectively.
The online execution time is here the time required by the ROM prediction in the test set. The difference in the execution times between 1D airfoil and 2D internal fields for the reduction techniques is caused by the different initial dimension. Indeed, the compression of the 2D fields' snapshots' matrix takes more time than in the 1D fields' case.

We also specify that the initial dimension of the snapshots is $\ndof=965$ for the 1D airfoil fields, $\ndof\sim 2.5e5$ for 2D fields, in the first test case. In the second test case we have $\ndof=448$ for 1D fields, and $\ndof\sim2.5e5$ for 2D fields.
The results in terms of computational times are similar for the two test cases, so we report in the tables the execution time only for test case 1.

\begin{table}[htpb!]
    \centering
 \caption{ROM execution times in offline stage, for all the reduction and approximation techniques for $\reddim=3$.}
    \label{tab:rom-time}
    \begin{tabular}{>{\centering\arraybackslash}p{0.1\linewidth}
    >{\centering\arraybackslash}p{0.2\linewidth}
    >{\centering\arraybackslash}p{0.3\linewidth}
    >{\centering\arraybackslash}p{0.3\linewidth}
    }
    \bottomrule
    {\textbf{Model}}&Field& {Offline time - Test case 1 [s]}\\
    \midrule
\multirow{2}{*}{\textbf{POD}}&airfoil& $\sim \num{4e-3}$\\
 &internal&$\sim \num{2.18}$\\
   \midrule
  {\textbf{AE}}&airfoil& $\sim \num{40}-\num{50}$\\
  \midrule
 {\textbf{PODAE}}&internal&$\sim \num{25}-\num{35}$\\
 \midrule
 {\textbf{RBF}}&(all cases)&$\sim \num{3.2e-2}$\\
 \midrule{\textbf{GPR}}&(all cases)&$\sim \num{3.5e-2}$\\
   \midrule
   \multirow{2}{*}{\textbf{ANN}}&airfoil&$\sim \num{7}-\num{10}$\\
 &internal&$\sim \num{80}-\num{100}$\\
    \bottomrule
    \end{tabular}
\end{table}

\begin{table}[htpb!]
    \centering
 \caption{ROM execution times in online stage, for all the reduction and approximation techniques for $\reddim=3$.}
    \label{tab:rom-time-online}
    \begin{tabular}{>{\centering\arraybackslash}p{0.2\linewidth}
    >{\centering\arraybackslash}p{0.1\linewidth}
    >{\centering\arraybackslash}p{0.3\linewidth}
    }
    \bottomrule
    {\textbf{Model}}&Field& {Online time - Test case 1 [s]} \\
    \midrule
\multirow{2}{*}{\textbf{POD-RBF}}&airfoil& $\sim \num{1e-4}$\\
 &internal&$\sim \num{3.3e-2}$\\
   \midrule
  {\textbf{AE-RBF}}&airfoil& $\sim \num{5e-3}$\\
  \midrule
 {\textbf{PODAE-RBF}}&internal&$\sim \num{4.3e-2}$\\
 \midrule
 \multirow{2}{*}{\textbf{POD-GPR}}&airfoil& $\sim \num{2e-4}$\\
 &internal&$\sim \num{4.5e-2}$\\
   \midrule
  {\textbf{AE-GPR}}&airfoil& $\sim \num{6e-3}$\\
  \midrule
 {\textbf{PODAE-GPR}}&internal&$\sim \num{7.6e-2}$\\
  \midrule
 \multirow{2}{*}{\textbf{POD-ANN}}&airfoil& $\sim \num{2e-4}$\\
 &internal&$\sim \num{2e-3}$\\
   \midrule
  {\textbf{AE-ANN}}&airfoil& $\sim \num{4e-4}$\\
  \midrule
 {\textbf{PODAE-ANN}}&internal&$\sim \num{4.3e-2}$\\
    \bottomrule
    \end{tabular}
\end{table}

\subsection{Aggregation results}
\label{subsec:aggr-result}

This Section is dedicated to the results of the aggregated ROM models for the first and second test cases, addressed in \ref{subsubsec:aggr-result-1} and \ref{subsubsec:aggr-result-2}, respectively.

\subsubsection{Test case 1}
\label{subsubsec:aggr-result-1}

\begin{figure}[htpb!]
    \centering
    \subfloat[]{\includegraphics[width=1.\textwidth, trim={2.2cm 0 3cm 0}, clip]{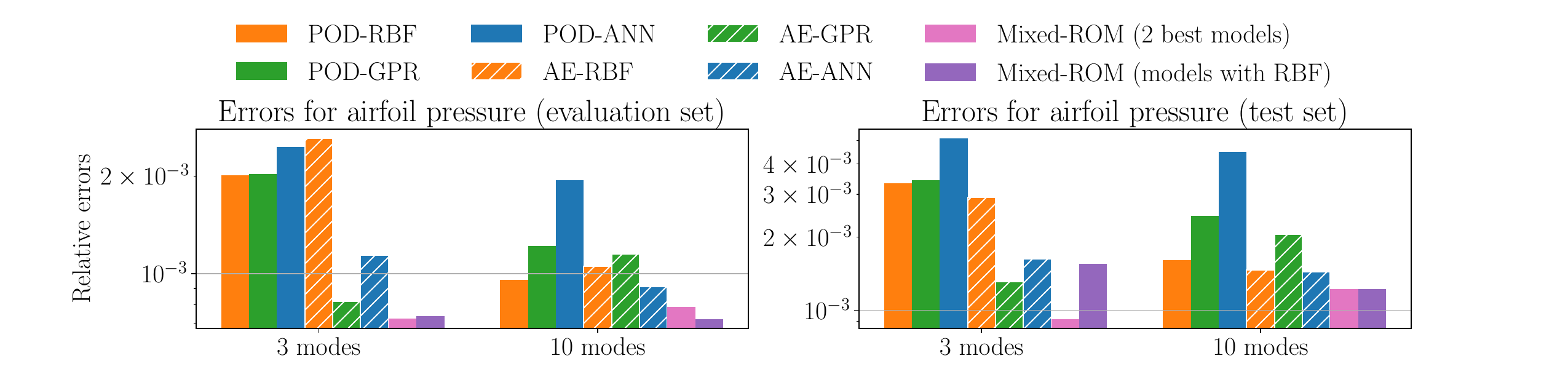}}\\
    \subfloat[]{\includegraphics[width=1.\textwidth, trim={2.2cm 0 3cm 0}, clip]{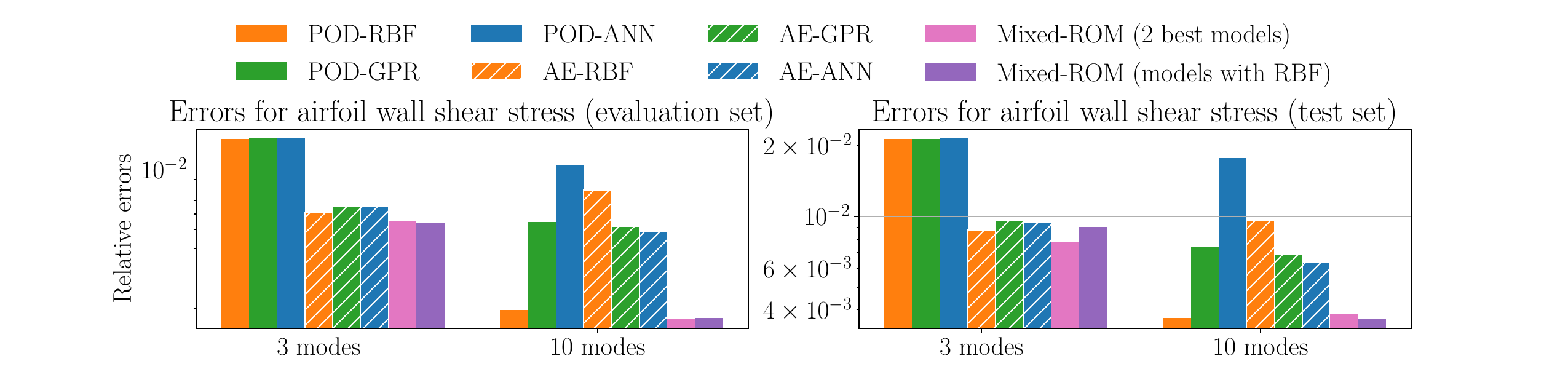}}
    \caption{Relative errors on evaluation and test set for 1D fields on the airfoil for all ROMs and for aggregated models in test case 1.}
    \label{fig:histograms-airfoil-case1}
\end{figure}

\begin{figure}[htpb!]
    \centering
    \subfloat[]{\includegraphics[width=1.\textwidth, trim={2.2cm 0 3cm 0}, clip]{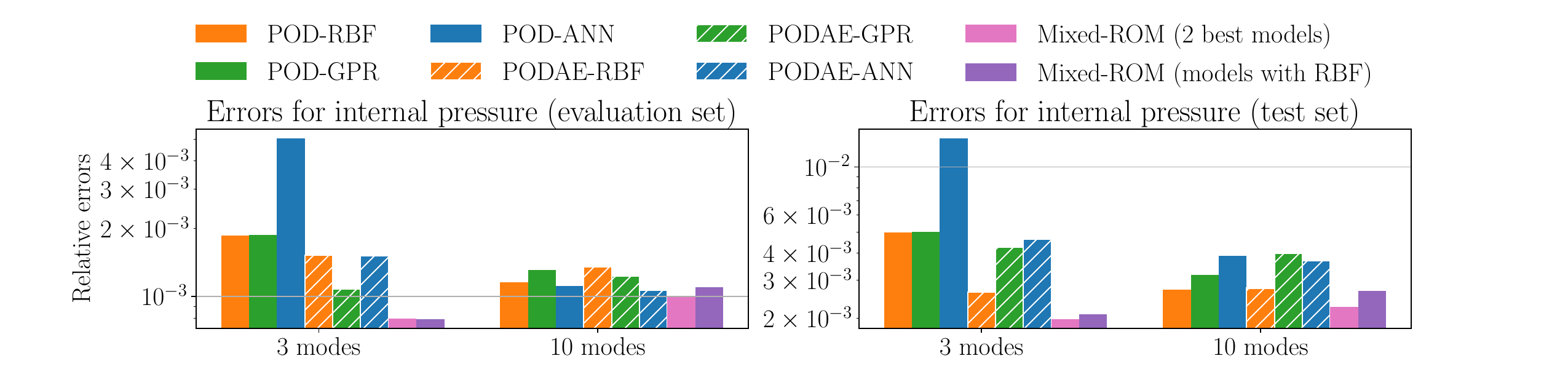}}\\
    \subfloat[]{\includegraphics[width=1.\textwidth, trim={2.2cm 0 3cm 0}, clip]{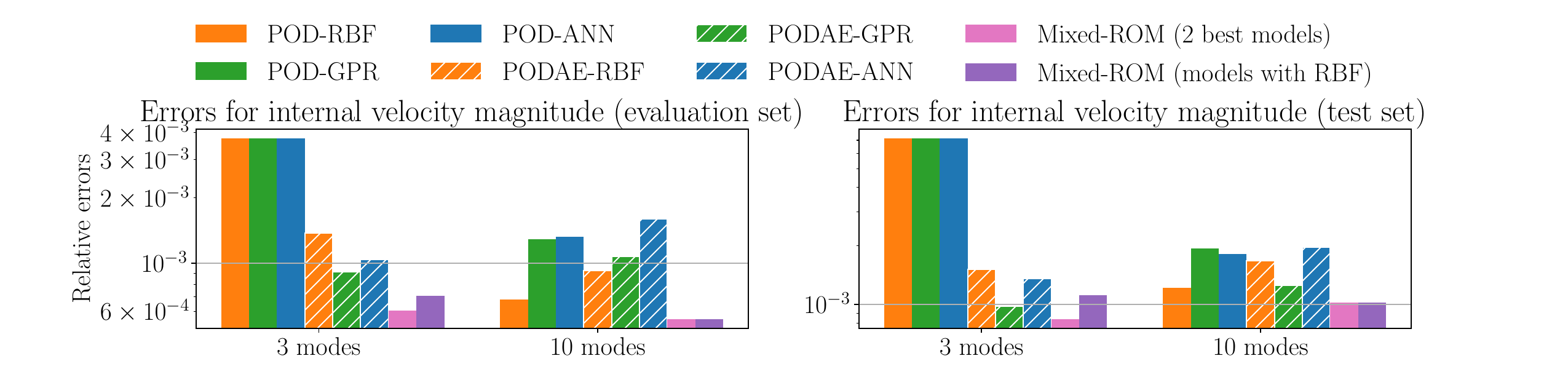}}
    \caption{Relative errors on evaluation and test set for 2D internal fields for all ROMs and for aggregated models in test case 1.}
    \label{fig:histograms-internal-case1}
\end{figure}

As a first quantitative analysis, we include Figures \ref{fig:histograms-airfoil-case1} and \ref{fig:histograms-internal-case1}, showing a comparison in the relative errors among the individual ROM models and the aggregated models (or \emph{mixed-ROMs}).
In particular, in both the test cases we built two different aggregation models, considering:
\begin{itemize}
    \item POD-RBF and AE(or PODAE)-RBF, namely two methods with different reduction techniques but with the same approximation method;
    \item the two ``\emph{best}" ROMs, based on the relative errors of the individual ROMs in the training set.
\end{itemize}

Based on Figures \ref{fig:histograms-airfoil-case1} and \ref{fig:histograms-internal-case1}, we can make different considerations. First of all, as already noticed in the preliminary analysis in \ref{subsubsec:rom-case1}, the performance of POD-based methods improves in general as the number of modes increases, whereas AE-based methods are not affected as much by the latent dimension.
Secondly, the aggregated models outperform always the individual ROMs originating the aggregations in the evaluation set. That is obvious by mathematical law, since the weights are calculated to minimize the discrepancy with the full-order counterpart in that specific set.
However, this may not happen in the test set, where the weights are computed through regression law, namely a Random Forest regression. In fact, the performance of the aggregation in the test set highly depends on the database used to fit and test the regression.
In particular, we have in all cases an improvement of the results in the test set for the airfoil pressure field and for the internal fields, but not for the wall shear field. In that case, we have comparable results with respect to the elementary non-intrusive ROMs. Moreover, the errors for the wall shear stress field are about one order of magnitude higher than for the other fields. This may lead to an increased complexity when combining the individual models.


\begin{figure}[htpb!]
    \centering{\includegraphics[width=0.6\textwidth]{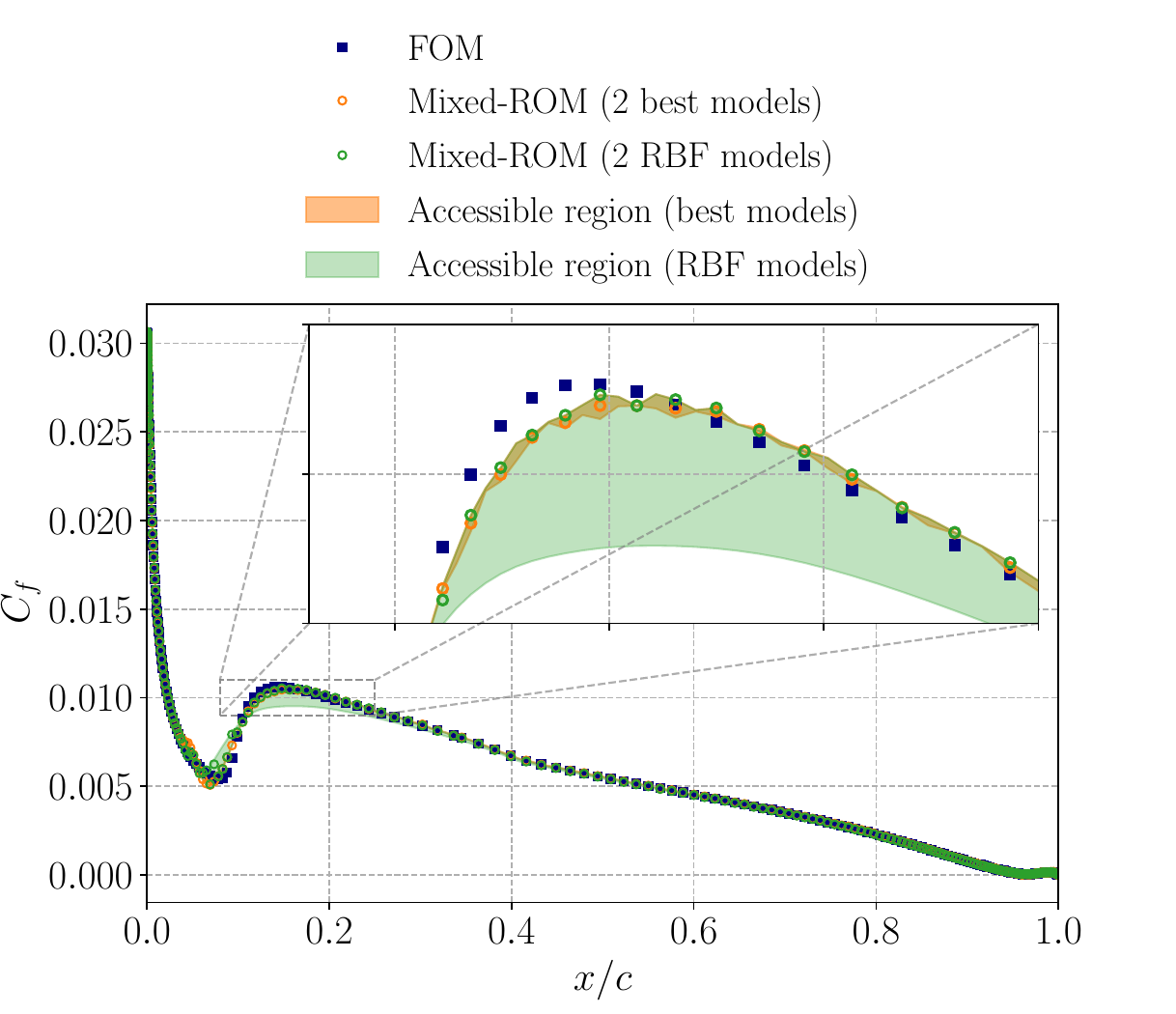}}
    \caption{Representation of the skin friction coefficient $C_f$ on the suction side for mixture models and FOM for a test parameter $Re \simeq 518000$. We consider 3 as reduced dimension. The \emph{accessible region}, namely the region accessed by the ROMs predictions considered in the aggregation, is also here represented.}
    \label{fig:cf-case1}
\end{figure}

As an example of how the aggregation performs on the chords coordinates, we include here Figure \ref{fig:cf-case1}, displaying the skin friction coefficient for a test parameter. In particular, the Figure shows the predictions of the \emph{mixed-ROM}, the FOM ground truth, and the regions accessed by the ROMs considered for the model mixtures. In addition, the light orange/green region is representative of the accessible area for the aggregated ROMs in the two types of aggregation here considered. It essentially quantifies the envelope of all the possible linear combinations of the individual ROMs aggregated.

The \emph{mixed-ROM} provides a prediction that is the closest as possible to the FOM snapshot, while remaining inside the accessible ROMs region, as can be seen in the zoomed region in Figure \ref{fig:cf-case1}.

\begin{figure}[htpb!]
    \centering{\includegraphics[width=1.\textwidth]{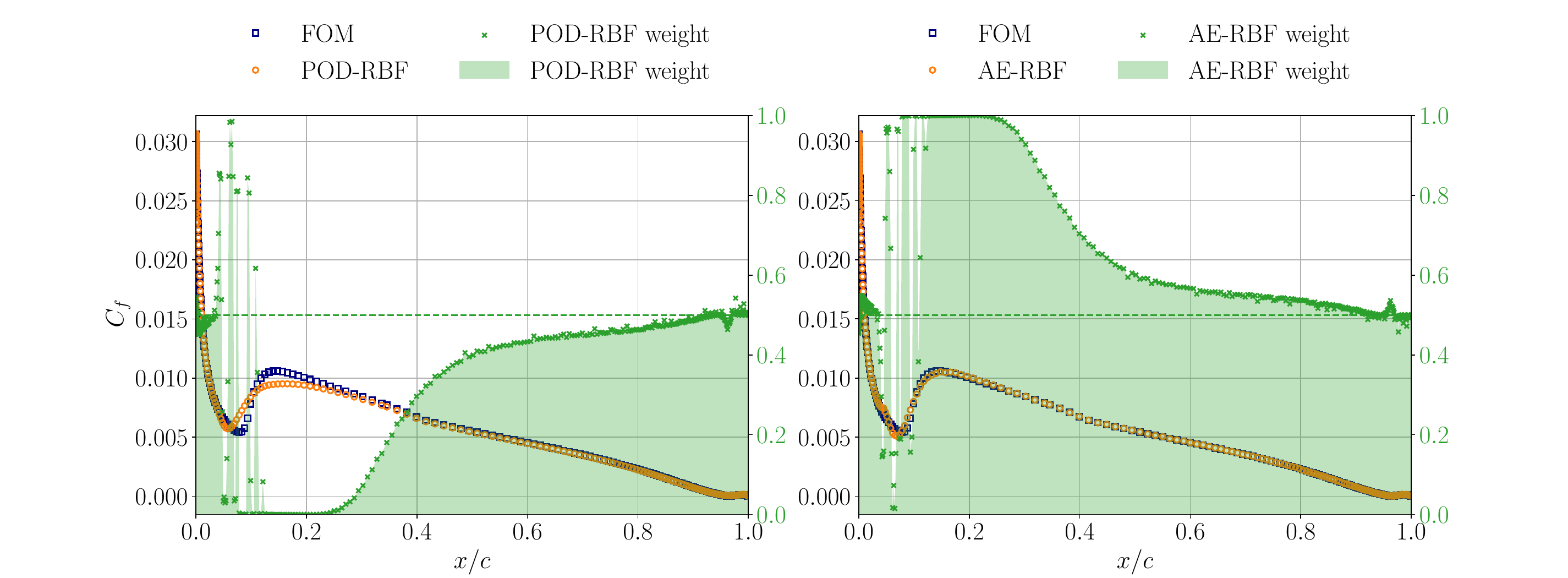}}
    \caption{Representation of the skin friction coefficient $C_f$ on the suction side for two ROMs (POD-RBF and AE-RBF) and of the corresponding weights for the aggregation with only RBF models. We consider as test parameter $Re \simeq 518000$ and as reduced dimension 3.}
    \label{fig:cf-weights-case1}
\end{figure}

Moreover, Figure \ref{fig:cf-weights-case1} provides the weights spatial distributions for the mixture of RBF-based models. The Figure shows that the AE-RBF is mostly activated where $x/c \in [0.15, 0.3]$, where it provides the most accurate approximation. Moreover, in the region $x/c \in [0.6, 1]$, since the two ROMs have similar predictions, the weights are close to 0.5, i.e., both are activated in the same percentage.

\begin{figure}[htpb!]
    \centering
\includegraphics[width=1.\textwidth, trim={4cm 0 3cm 0}, clip]{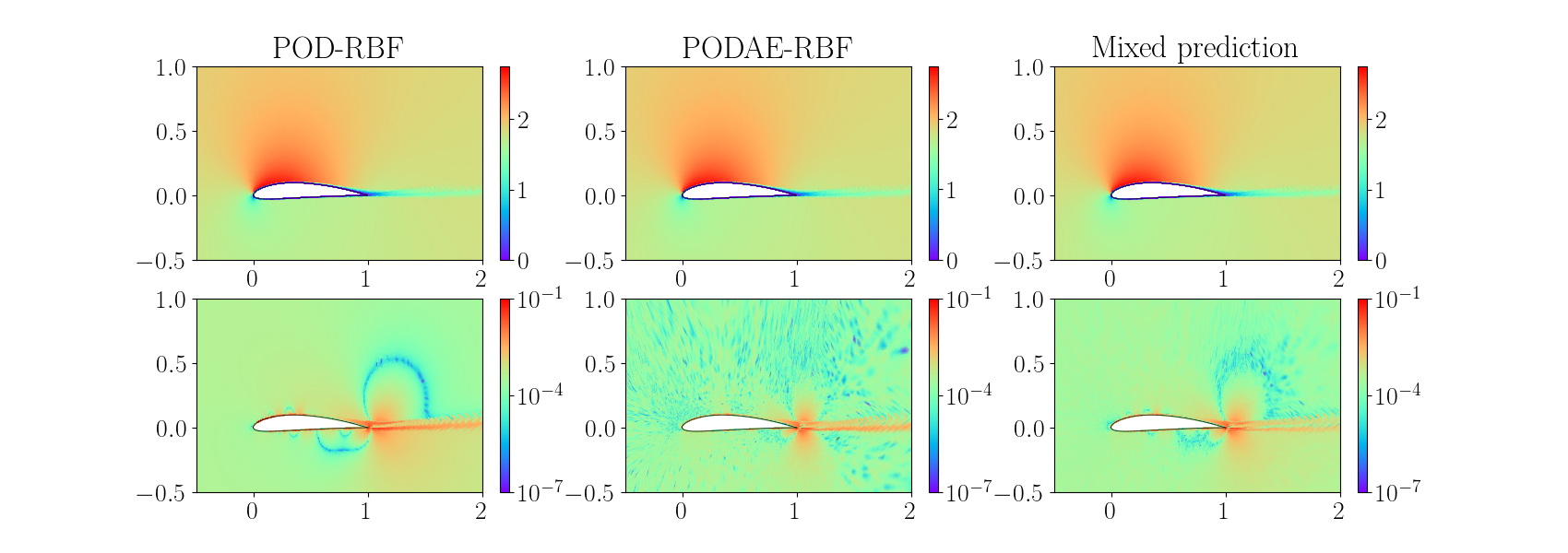}
    \caption{Velocity magnitude predictions for two ROMs and for the aggregation model (first row), and corresponding relative errors with respect to FOM (second row). The test parameter considered is $Re\approx 118000$ and the reduced dimension is 10.}
    \label{fig:v-10-internal-case1-Re118}
\end{figure}

\begin{figure}[htpb!]
    \centering
    \includegraphics[width=0.8\textwidth]{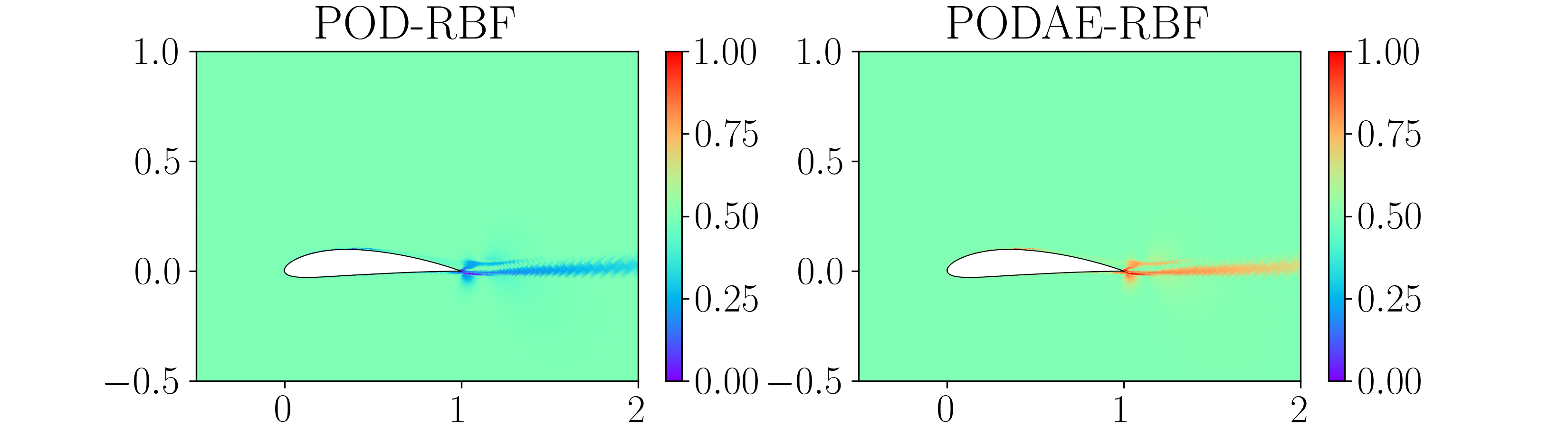}
    \caption{Weights for RBF models for the velocity magnitude field when considering the test parameter $Re \approx 118000$. The reduced dimension is 10.}
    \label{fig:v-10-weights-case1-Re118}
\end{figure}

If we consider a velocity internal field (Figures \ref{fig:v-10-internal-case1-Re118} and \ref{fig:v-10-weights-case1-Re118}), we can notice that the two ROMs have a similar performance and are really close to the FOM reference.

However, the AE-RBF model is closer to the FOM prediction in the wake, as can be seen from the spatial relative error in Figure \ref{fig:v-10-internal-case1-Re118}, and, hence, the corresponding weights are close to 1 (i.e., full activation) in that specific zone of the domain (Figure \ref{fig:v-10-weights-case1-Re118}).


\subsubsection{Test case 2}
\label{subsubsec:aggr-result-2}
Based on the results obtained in the previous section, we wanted to consider a test case that is more challenging in terms of reduction. Consequently, we moved towards a transonic airfoil where the angle of attack has been chosen as the parameter of interest. In this scenario, due to the significant movement of the shock wave for different parametric states, the reduction is expected to be more demanding.

\begin{figure}[htpb!]
    \centering
    \subfloat[]{\includegraphics[width=1.\textwidth, trim={2.2cm 0 3cm 0}, clip]{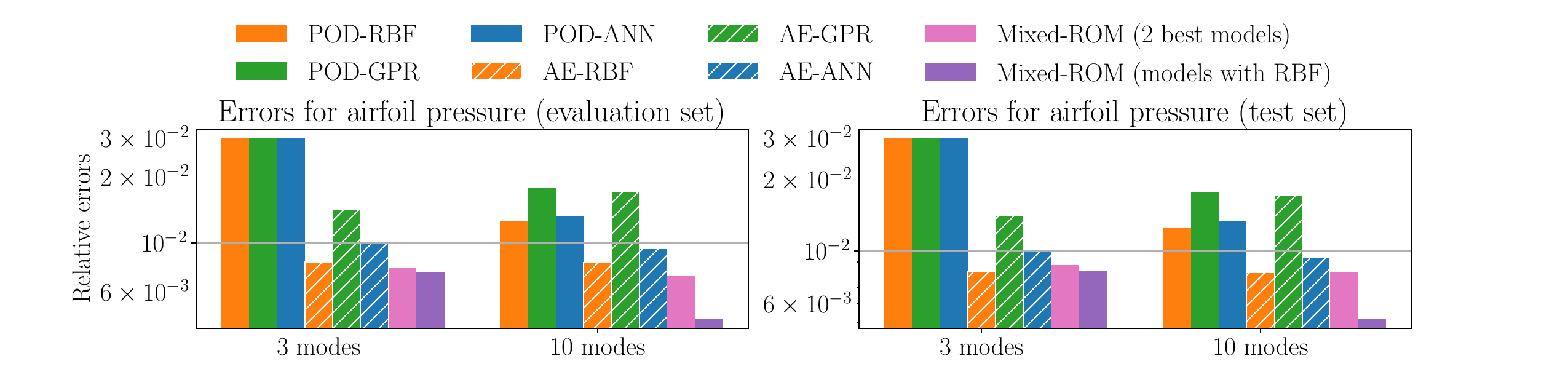}}\\
    \subfloat[]{\includegraphics[width=1.\textwidth, trim={2.2cm 0 3cm 0}, clip]{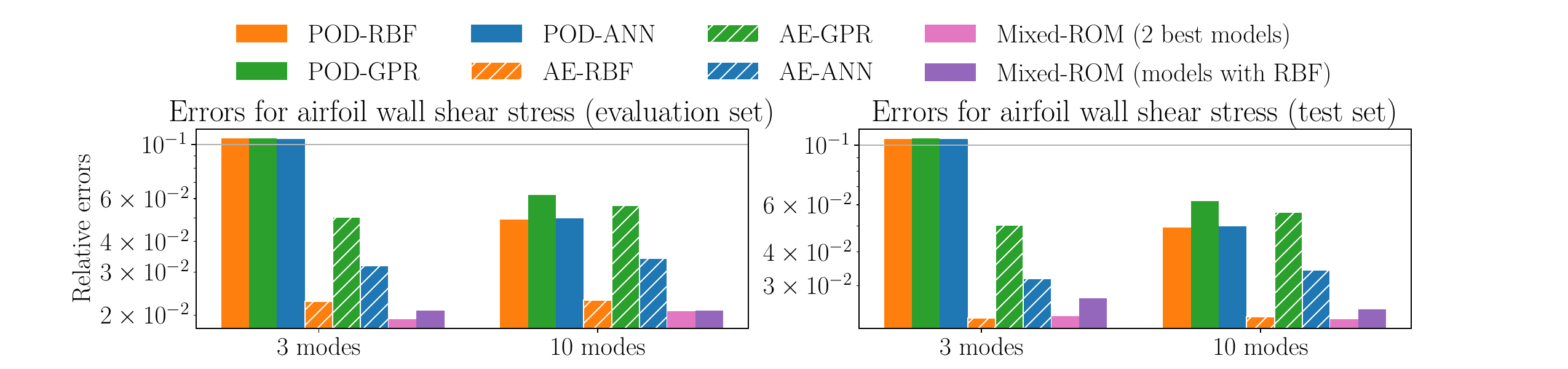}}
    \caption{Relative errors on evaluation and test set for the 1D fields on the airfoil for all ROMs and for aggregated models.}
    \label{fig:histograms-airfoil-case2}
\end{figure}

\begin{figure}[htpb!]
    \centering
    \subfloat[]{\includegraphics[width=1.\textwidth, trim={2.2cm 0 3cm 0}, clip]{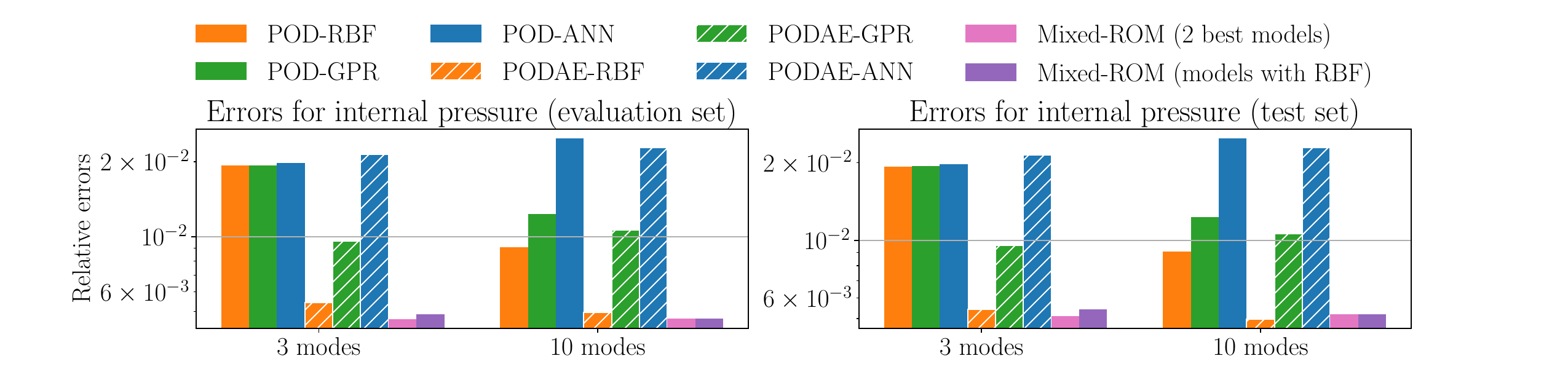}}\\
    \subfloat[]{\includegraphics[width=1.\textwidth, trim={2.2cm 0 3cm 0}, clip]{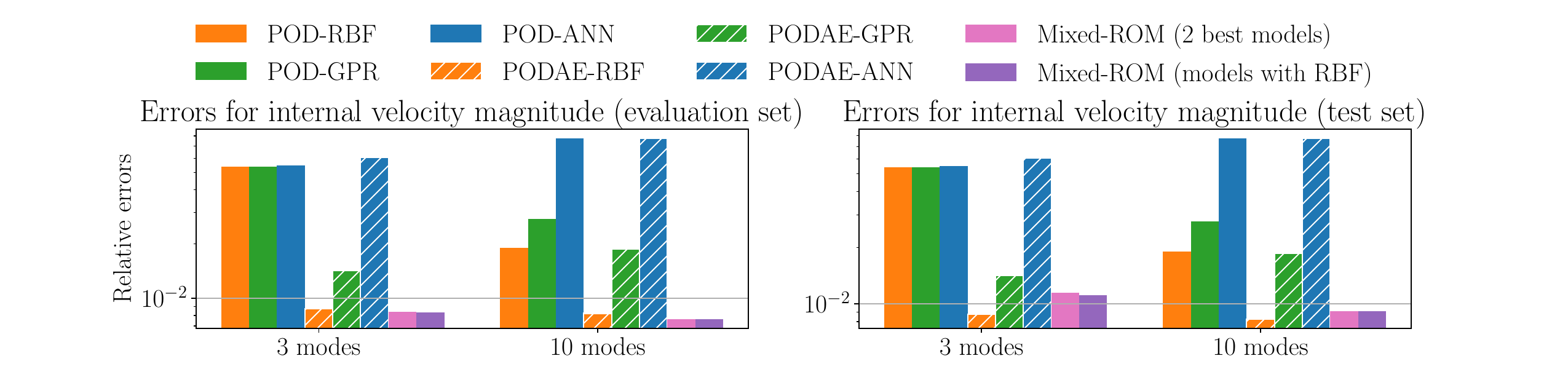}}
    \caption{Relative errors on evaluation and test set for the 2D fields in the internal mesh for all ROMs and for aggregated models.}
    \label{fig:histograms-internal-case2}
\end{figure}

In Figure \ref{fig:histograms-airfoil-case2} and \ref{fig:histograms-internal-case2} we show the relative errors of each individual ROM (varying the type of reduction and approximation in the parametric space) and with the different type of mixing combinations.

As a general observation, we can notice that the relative error of the aggregated models is always significantly small.
The same type of analysis is carried out also on the variables on the internal mesh with respect to the airfoil. The relative errors are thereby shown in Figure \ref{fig:histograms-internal-case2}. Also in this case we can observe that the aggregated model can significantly improve the best possible combination of individual ROMs, or, at least, identify the best possible model for each configuration. 

\begin{figure}[htpb!]
    \centering{\includegraphics[width=1\textwidth]{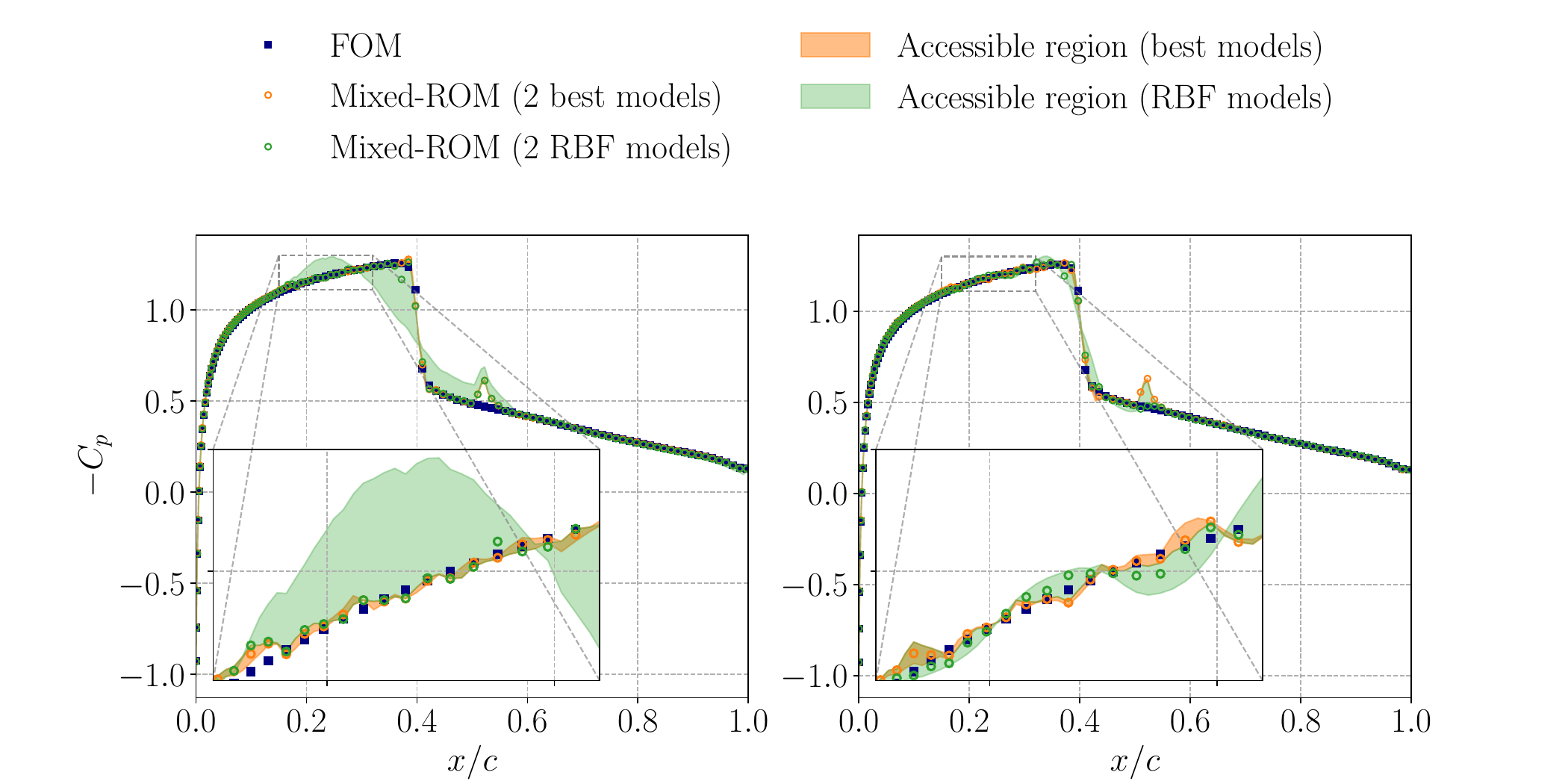}}
    \caption{Representation of the pressure coefficient $C_p$ for mixture models and FOMs for a test parameter $\alpha=4.7^{\circ}$. We consider two different reduced dimensions, 3 (on the left) and 10 (on the right). The \emph{accessible region}, namely the region accessed by the ROMs predictions considered in the aggregation, is also here represented.}
    \label{fig:cp-case2}
\end{figure}

In Figure \ref{fig:cp-case2}, we show the results for a specific solution of the test set. In particular, we depict the pressure coefficient at the wall given by the full-order model and the two aggregated ROMs. On the left, we show the results obtained with a latent space of dimension $3$ and on the right with a latent space of dimension $10$. 
Clearly, the accessible area may be significantly larger when the dimension of the latent space is small: in this scenario we can expect much different results from each individual ROM. As the dimension of the latent space increases, all of the different ROM will tend to converge to the exact solution, limiting the variability among them. 
We can notice that within this region, the aggregated model correctly predicts a mixed solution that is as close as possible to the reference.

It is important to stress that having a significantly large accessible region is fundamental in the success of the proposed strategy. Such a feature, in fact, indicates that the large differences between the models can be exploited in order to obtain a better aggregated solution. It can be noticed, in fact, that considering a latent space of dimension $10$, the prediction of the aggregated model is better when using the two models using RBF interpolation rather than the two best models. In this latter case, in fact, the differences would be so small between the two best models, that there would be practically no benefit in their aggregation. 
\begin{figure}[htpb!]
    \centering{\includegraphics[width=0.6\textwidth]{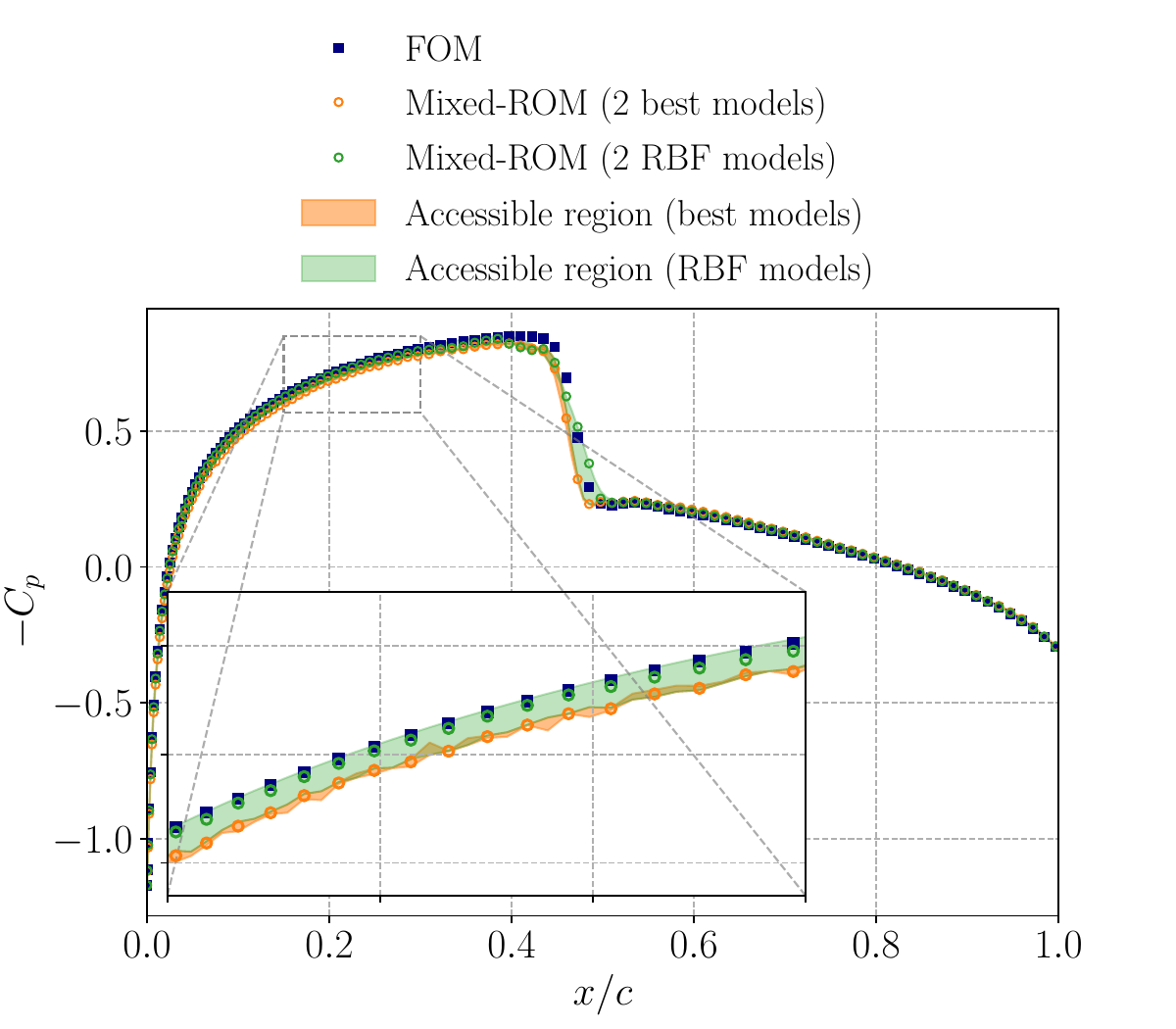}}
    \caption{Representation of the pressure coefficient $C_p$ for mixture models and FOMs for a test parameter $\alpha=0.3^{\circ}$. We consider as reduced dimension $\reddim=10$. The \emph{accessible region}, namely the region accessed by the ROMs predictions considered in the aggregation, is also here represented.}
    \label{fig:cp-case2-alphalow-10modes}
\end{figure}

\begin{figure}[htpb!]
    \centering
    \subfloat[Latent dimension$=3$]{\includegraphics[width=1.\textwidth]{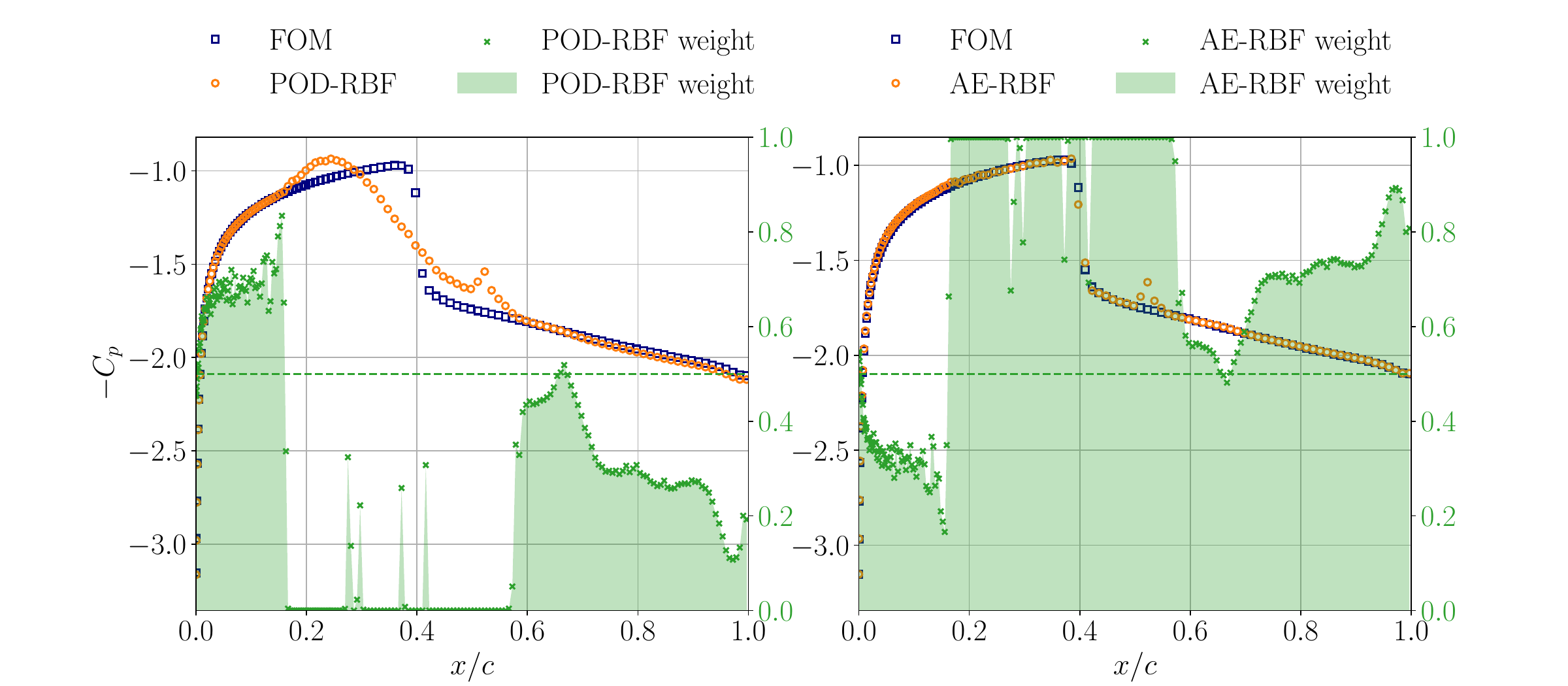}}\\
    \subfloat[Latent dimension$=10$]{\includegraphics[width=1.\textwidth]{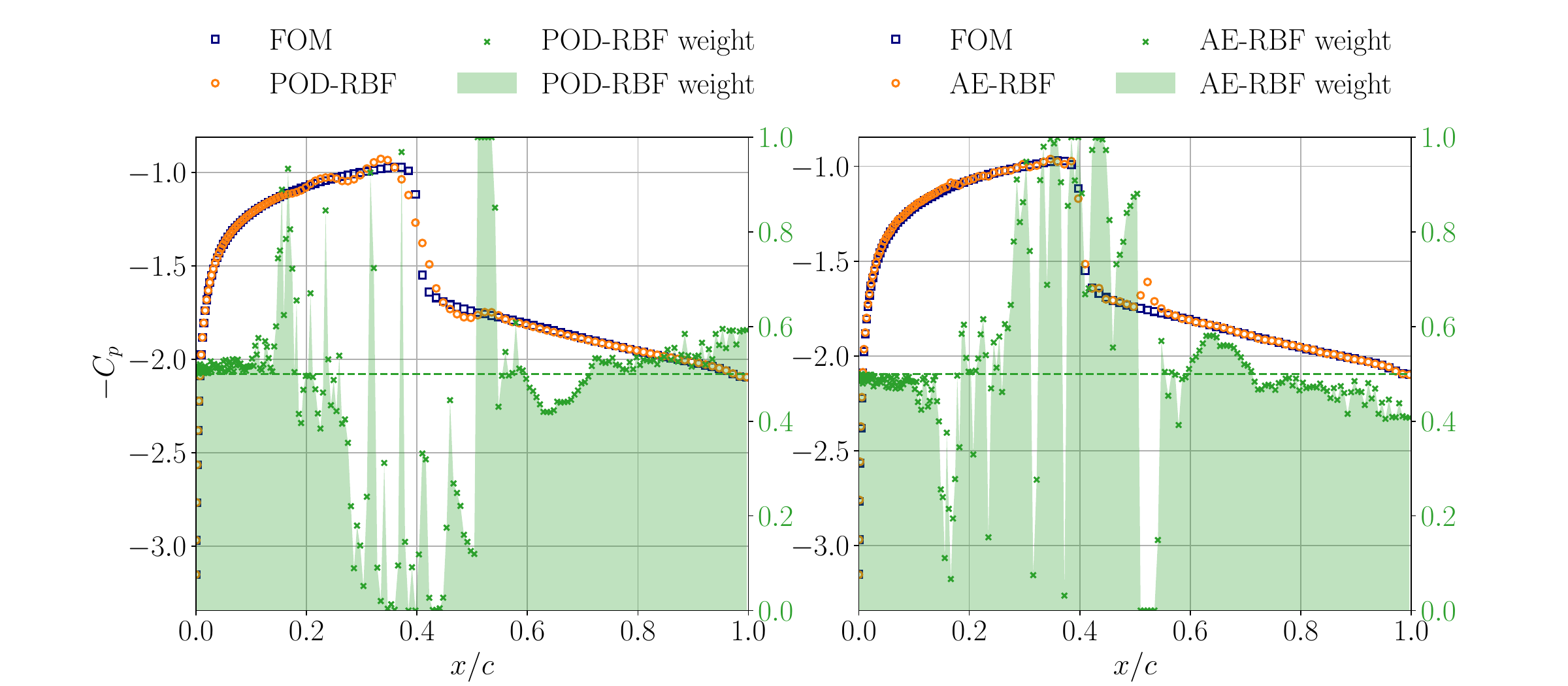}}
    \caption{Representation of the pressure coefficient $C_p$ for two ROMs (POD-RBF and AE-RBF) and of the corresponding weights for the aggregation with only RBF models. We consider as test parameter $\alpha=4.7^{\circ}$ and two different reduced dimensions, 3 and 10.}
    \label{fig:cp-weights-case2}
\end{figure}

To better understand the interplay between the different models and how they work together, in Figure \ref{fig:cp-weights-case2} we show also the values of the weights of the elementary ROMs. We can observe that as one might expect, the intrinsic non-linearity of the autoencoder helps in predicting correctly the sharp feature represented by the shock wave. This automatically implies a full activation of the AE at this location, whereas in other regions, such as in proximity of the leading edge, the POD is more active. This is more evident when the number of modes is small. In fact, in this scenario, the POD prediction does not provide a lot of information to properly reconstruct the manifold in the parametric space. By increasing the number of modes (right figure), we can still observe similar trends, although the POD is reasonably good in proximity of the shock and the activation of the AE is less evident.

\begin{figure}[htpb!]
    \centering
    \includegraphics[width=1.\textwidth, trim={4cm 0 3cm 0}, clip]{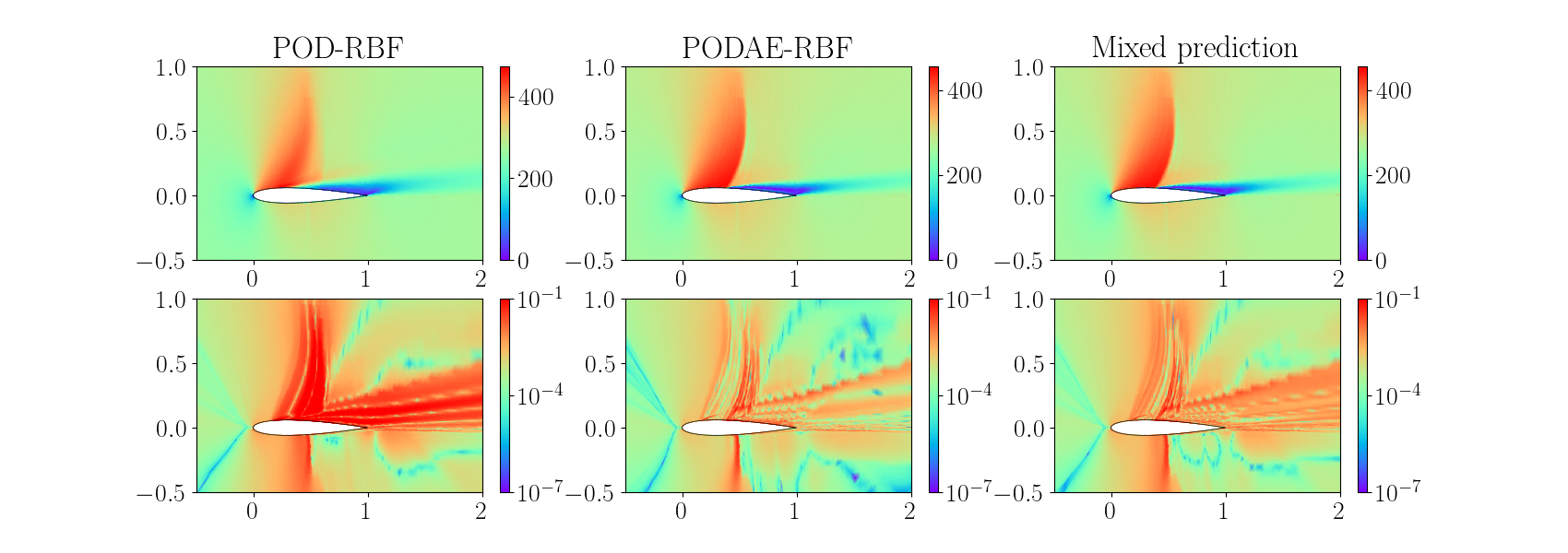}
    \caption{Velocity magnitude predictions for two ROMs and for the aggregation model (first row), and corresponding absolute errors with respect to FOM (second row). The test parameter considered is $\alpha=4.7^{\circ}$ and the reduced dimension is 3.}
    \label{fig:v-3-internal-case2}
\end{figure}

The velocity magnitude on the surroundings of the airfoil is shown in Figure \ref{fig:v-3-internal-case2}. In this figure, the POD and POD-AE reductions are compared and aggregated together.

\begin{figure}[htpb!]
    \centering
    \includegraphics[width=0.8\textwidth]{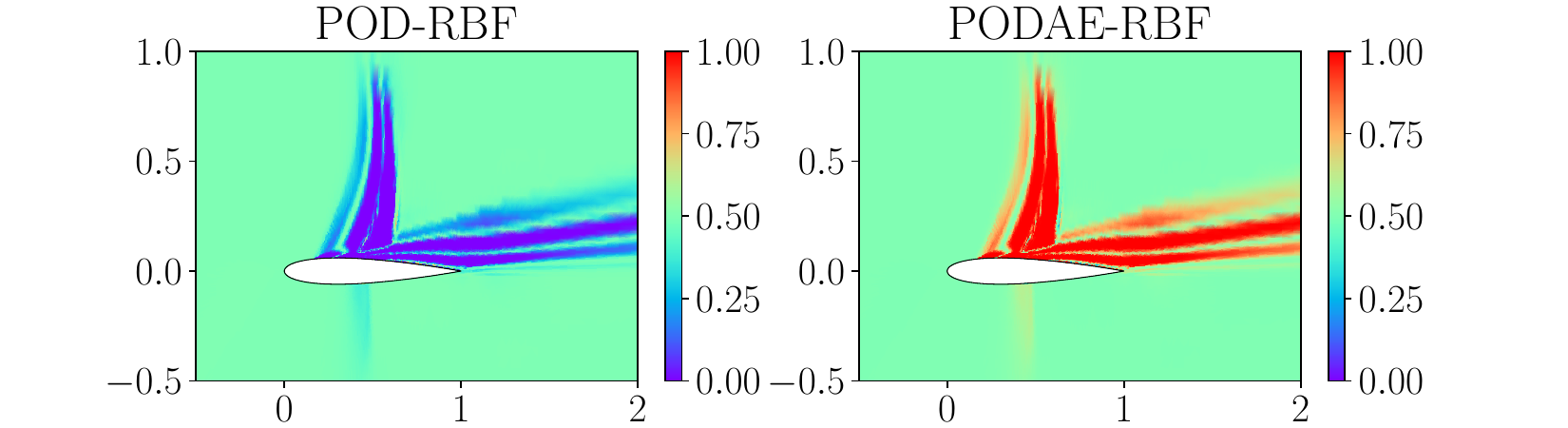}
    \caption{Weights representation for the aggregation model when considering the POD-RBF and PODAE-RBF models. The test parameter considered is $\alpha=4.7^{\circ}$ and the reduced dimension is 3.}
    \label{fig:v-weights-3-internal-case2}
\end{figure}

Both models provide a reasonably good approximation of the unseen parametric solution. However, if we observe more closely, the weights of the two models (shown in Figure \ref{fig:v-weights-3-internal-case2}), we can clearly observe that regions which are particularly challenging in terms of reduction are correctly detected in the spatial domain and the non-linear technique based on PODAE is much more active with respect to the linear POD.

If we increase the dimension of the latent space it is even more interesting to notice that the sheet behind the airfoil is not a challenging feature anymore and both models are correctly predicting its behavior. However, at the shock, the PODAE is much more active than the POD. This observation is a further indication that even if turbulent boundary layers can be challenging in terms of reduction, sharp features moving significantly in the parametric space (such as shock-waves) are still the dominant difficulty in terms of reduced-order modelling. This aspect reiterates how the second test case is more complex to reduce and represents a more challenging scenario for the aggregation of ROMs.

\subsubsection{Computational time}
\label{subsubsec:time-aggr}
For what concerns the computational time needed for the Random Forest Regression, most of the execution time is due to the fitting process to the evaluation data. In particular, the execution time for fitting the regression is $\sim 1$ s for the 1D airfoil fields, and $\sim 40-45$ s for 2D fields, because of the larger space dimensionality.

The execution times for prediction are much smaller and comparable with the prediction times for the individual ROMs, i.e., $\sim \num{1e-2}$ for the airfoil fields and $\sim \num{1e-1}$ for the internal fields.
\newpage

\section{Conclusions}\label{sec:conclusions}
In this manuscript we proposed a strategy to enhance the performance of individual data-driven ROMs, which are previously trained in the so-called \emph{training set}. The ROM methodology is extensively discussed in Section \ref{sec:methods}.

In particular, each individual ROM is mainly characterized by two parts, a \emph{reduction} stage and an \emph{approximation} step, which have been discussed in Section \ref{met-roms}.

The space-dependent aggregation approach, described in Section \ref{subsec:met-aggregation}, builds convex combination of a set of alternative ROMs, with space-dependent weights evaluated on an \emph{a-priori-selected} \emph{evaluation} set. A machine learning technique (Random Forest Regression) is then trained with the given weights' data, and used to predict the weights in the \emph{test set}.

The model mixture, namely the \emph{mixed ROM}, is tested on two different test cases. Both the cases refer to the flow past an airfoil, but use different FOMs and present a different parameter, the Reynolds number in the first test case, and the angle of attack in the second test case.
The results of the aggregation approach are presented in Section \ref{sec:results}.

In both the test cases the proposed methodology has proved to improve the accuracy of the individual ROMs. It is indeed able to identify and activate the ROM with the best performance in different regions of the space and parameter domains.
For instance, in the second test case, that is more challenging than the first one, the aggregation approach was successful in fully activating the nonlinear reduction techniques in the region close to the shock position.

Therefore, we can definitely say that it provides an automatic detection of the need between linear and nonlinear reduction without any a-priori knowledge of the problem of interest.

Moreover, the computational time needed for computing the aggregation is small if compared to the FOM simulation time. In particular, most of the execution time in the aggregation is employed for fitting the RF regression, but it takes only a few seconds, as pointed out in section \ref{subsubsec:time-aggr}.

It is important to remark that the performance of the technique in the test set strongly depends on the type of regression used to predict the weights' space distributions.
One of the possible extensions of this work will be to test the method with more advanced machine learning techniques for the weights' regression.

A possible extension is to change the features space in order to allow for different space coordinates, and hence, different geometries.

Another extension that the authors will address in the future is related to more complex test cases, and, in particular, to the mixture of ROMs built on the top of FOMs characterized by different turbulence models. In this case, a further level of variability at the ROM level would be introduced by the different performances of the turbulence models at the FOM level.





\section*{Acknowledgements}
This study was funded by the European Union - NextGenerationEU, in the framework of the iNEST -
Interconnected Nord-Est Innovation Ecosystem (iNEST ECS00000043 – CUP G93C22000610007). The views and
opinions expressed are solely those of the authors and do not necessarily reflect those of the
European Union, nor can the European Union be held responsible for them. In addition, the authors would like to acknowledge INdAM–GNCS.

\bibliographystyle{plain}
\bibliography{main}

\section*{Supplementary material}
In this supplementary Section, we report the specifications of the hyperparameters of the techniques exploited in non-intrusive ROMs.

In particular, in the RBF interpolation, we consider a \emph{thin plate spline
} kernel, which reads $\varphi(r)=r^2 \text{log}(r)$, and the degree of the added polynomial is $0$.

In the GPR approximation, we consider a squared exponential form:
\begin{equation}
\mathcal{K}(\bm{\mu_i}, \bm{\mu_j}) = \sigma^2 \exp{\left( -\dfrac{\|\bm{\mu_i}-\bm{\mu_j}\|^2}{2l} \right)}.
    \label{eq:GP_form}
\end{equation}

Once the hyperparameters $\sigma$ and $l$ of the covariance kernel are computed in order to fit our data, we can exploit the distribution \eqref{eq:GP} to evaluate $\mathcal{M}(\bm{\mu^*})$ in \eqref{eq:GP} and predict the new modal coefficients. The hyperparameters are automatically fitted by the algorithm itself.

\begin{table}[h!]
    \centering
 \caption{Neural networks setting in ROMs.}
    \label{tab:networks_first}
    \begin{tabular}{>{\centering\arraybackslash}p{0.13\linewidth}
    >{\centering\arraybackslash}p{0.13\linewidth}
    >{\centering\arraybackslash}p{0.12\linewidth}
    >{\centering\arraybackslash}p{0.14\linewidth}
    >{\centering\arraybackslash}p{0.06\linewidth}
    >{\centering\arraybackslash}p{0.06\linewidth}
    >{\centering\arraybackslash}p{0.12\linewidth}
    }
    \bottomrule
    \multirow{3}{*}{\textbf{Network}}& \multirow{3}{*}{Hidden layers} & \multirow{3}{*}{Non-linearity} &\multirow{3}{*}{Learning rate} & \multicolumn{2}{c}{\multirow{2}{*}{Stop criteria}}& \multirow{3}{*}{Weight decay}\\
    &&&&&\\
    &&&& epochs & final loss &\\
    \midrule
\textbf{ANN} &{$[20,20,20]$}&{Softplus}& {$\num{5e-3}$} & {100000} & {$\num{1e-4}$} & {$\num{1e-7}$}\\ 
   \midrule
  \multirow{2}{*}{\textbf{AE}} &{$[50, 20, \reddim,$}&\multirow{2}{*}{Softplus}& \multirow{2}{*}{$\num{5e-4}$} & \multirow{2}{*}{20000} & \multirow{2}{*}{$\num{5e-6}$} & \multirow{2}{*}{0}\\
   &{$,\reddim, 20, 50]$}&&&&&\\
    
   \midrule
  \multirow{2}{*}{\textbf{AE} (\textbf{PODAE})} &
  {$[\mediumdim, 50, 20, \reddim,$} &\multirow{2}{*}{Softplus}& \multirow{2}{*}{$\num{5e-4}$} & \multirow{2}{*}{20000} & \multirow{2}{*}{$\num{5e-6}$} & \multirow{2}{*}{0}\\
   &{$,\reddim, 20, 50, \mediumdim]$} &&&&&\\
  
    \bottomrule
    \end{tabular}
\end{table}

\end{document}